\documentclass[12pt]{amsart}

\usepackage{amsmath, amssymb, amsfonts, amsthm, amscd}
\usepackage{manfnt}
\usepackage{mathtools}
\usepackage{fullpage}
\usepackage{url}
\usepackage[all]{xy}
   \SelectTips{cm}{10}
\usepackage{txfonts}
\usepackage{enumitem}
\newtheorem{thm}{Theorem}[section]
\newtheorem{lemma}[thm]{Lemma}
\newtheorem{prop}[thm]{Proposition}
\newtheorem{cor}[thm]{Corollary}

\newtheorem{prop-conj}[thm]{Proposition-Conjecture}

\theoremstyle{definition}
\newtheorem{defn}[thm]{Definition}

\theoremstyle{remark}
\newtheorem{rmk}[thm]{Remark}

\newtheorem{assumption}[thm]{Assumption}

\theoremstyle{remark}

\theoremstyle{remark}

\theoremstyle{remark}
\newtheorem{eg}[thm]{Example}

\newcommand{\Q}{\mathbb{Q}}
\newcommand{\Qb}{\overline{\mathbb{Q}}}
\newcommand{\Z}{\mathbb{Z}}
\newcommand{\CC}{\mathbb{C}}

\newcommand{\Ql}{\mathbb{Q}_{\ell}}
\newcommand{\Qlb}{\overline{\mathbb{Q}}_\ell}

\newcommand{\fl}{\mathbb{F}_{\ell}}

\DeclareMathOperator{\Hom}{Hom}

\DeclareMathOperator{\Aut}{Aut}
\DeclareMathOperator{\Out}{Out}

\DeclareMathOperator{\Sym}{Sym}
\DeclareMathOperator{\Ad}{Ad}
\DeclareMathOperator{\ad}{ad}

\DeclareMathOperator{\im}{im}

\DeclareMathOperator{\Fil}{Fil}

\DeclareMathOperator{\rk}{rk}

\DeclareMathOperator{\Rep}{Rep}
\DeclareMathOperator{\Spec}{Spec}
\DeclareMathOperator{\id}{id}

\DeclareMathOperator{\Frac}{Frac}
\DeclareMathOperator{\Lie}{Lie}

\newcommand{\gal}[1]{\Gamma_{#1}} 
\newcommand{\Gal}{\mathrm{Gal}} 

\newcommand{\into}{\hookrightarrow}
\newcommand{\onto}{\twoheadrightarrow}

\newcommand{\mc}{\mathcal}
\newcommand{\mf}{\mathfrak}
\newcommand{\mr}{\mathrm}
\newcommand{\mbf}{\mathbf}

\newcommand{\tG}{\widetilde{G}}

\newcommand{\br}{\bar{\rho}} 
\newcommand{\fg}{\mathfrak{g}}
\newcommand{\fb}{\mathfrak{b}}
\newcommand{\fn}{\mathfrak{n}}
\DeclareMathOperator{\Lift}{\mathrm{Lift}}
\DeclareMathOperator{\Def}{\mathrm{Def}}
\newcommand{\Sets}{\mathbf{Sets}}
\newcommand{\CO}{\mathcal{C}_{\mathcal{O}}}
\newcommand{\fr}{\mathrm{fr}}
\newcommand{\kbar}{\overline{\kappa}}
\newcommand{\Ram}{\mathrm{Ram}}
\newcommand{\unr}{\mathrm{unr}}
\newcommand{\brgm}{\bar{\rho}(\mathfrak{g}_{\mu})}
\newcommand{\St}{\mathrm{St}}
\newcommand{\tv}{\tilde{v}}
\newcommand{\brgd}{\bar{\rho}(\mathfrak{g}^\vee)}
\newcommand{\tF}{\widetilde{F}}

\begin{document}

\author{Stefan Patrikis}
\email{patrikis@math.utah.edu}
\address{Department of Mathematics\\University of Utah\\Salt Lake City, UT 84112}
\date{June 2015}

\title{Deformations of Galois representations and exceptional monodromy}
\thanks{This paper has a long history, and it gives me great pleasure to extend thanks both old and new. In 2006 with Richard Taylor's guidance I proved a version of Ramakrishna's lifting theorem for symplectic groups, as my Harvard undergraduate thesis. Without Richard's encouragement and singular generosity, I would likely not have continued studying mathematics, much less been equipped to write the present paper; I am enormously grateful to him. At the time Frank Calegari and Brian Conrad both read and helpfully commented on the thesis, and I thank them as well. Moving toward the present day, I am grateful to Dick Gross, from whom I learned an appreciation of the principal $\mr{SL}_2$, which turned out to be crucial for this paper; and to Shekhar Khare, for his comments both on this paper and on other aspects of Ramakrishna's work. Finally, I am greatly indebted to Rajender Adibhatla, whose proposal of the joint project \cite{adibhatla-patrikis} spurred me to revisit Ramakrishna's methods.}
\begin{abstract}
For any simple algebraic group $G$ of exceptional type, we construct geometric $\ell$-adic Galois representations with algebraic monodromy group equal to $G$, in particular producing the first such examples in types $\mr{F}_4$ and $\mr{E}_6$. To do this, we extend to general reductive groups Ravi Ramakrishna's techniques for lifting odd two-dimensional Galois representations to geometric $\ell$-adic representations.
\end{abstract}
\maketitle

\section{Introduction}
Prior to the proof of Serre's conjecture on the modularity of odd representations 
\[
\br \colon \Gal(\Qb/\Q) \to \mr{GL}_2(\overline{\mathbb{F}}_{\ell}),
\]
one of the principal pieces of evidence for the conjecture was a striking theorem of Ravi Ramakrishna (\cite{ramakrishna02}) showing that under mild hypotheses such a $\br$ could be lifted to a geometric (in the sense of Fontaine-Mazur) characteristic zero representation. A generalization of Ramakrishna's techniques to certain $n$-dimensional representations then played a key part in the original proof of the Sato-Tate conjecture for rational elliptic curves (\cite{hsbt:dwork}). Meanwhile, dramatic advances in potential automorphy theorems, culminating in \cite{blggt:potaut}, have drawn attention away from Ramakrishna's method, since for suitably odd, regular, and self-dual representations
\[
\br \colon \Gal(\Qb/\Q) \to \mr{GL}_N(\overline{\mathbb{F}}_{\ell}),
\]
the potential automorphy technology, combined with an argument of Khare-Wintenberger (see eg \cite[\S 4.2]{kisin:cdm}) now yields remarkably robust results on the existence of characteristic zero lifts of prescribed inertial type (eg, \cite[Theorem E]{blggt:potaut}). Put another way, there are connected reductive groups $G/\Q$ of classical type in the Dynkin classification, for which one can apply potential automorphy theorems to find geometric characteristic zero lifts of certain
\[
\br \colon \Gal(\Qb/\Q) \to {}^L G(\overline{\mathbb{F}}_{\ell}),
\]
where ${}^L G$ denotes as usual an L-group of $G$. But these potential automorphy techniques are currently quite limited outside of classical types,
for neither the existence of automorphic Galois representations, nor any hoped-for potential automorphy theorems, have been demonstrated in more general settings. The first aim of this paper is to prove a generalization, modulo some local analysis, of Ramakrishna's lifting result for essentially arbitrary reductive groups. Such results have some intrinsic interest, and among other things, they provide evidence for generalized Serre-type conjectures. 

At first pass, we do this under very generous assumptions on the image of $\br$, analogous to assuming the image of a two-dimensional representation contains $\mr{SL}_2(\mathbb{F}_{\ell})$. Here is a special case: see Theorem \ref{maximaltheorem} and Theorem \ref{maximalimageout} for more general versions.
\begin{thm}\label{maximalintro}
Let $\mc{O}$ be the ring of integers of a finite extension of $\Ql$, and let $k$ denote the residue field of $\mc{O}$. Let $G$ be an adjoint Chevalley group, and let $\br \colon \Gal(\Qb/\Q) \to G(k)$ be a continuous representation unramified outside a finite set of primes $\Sigma$ containing $\ell$ and $\infty$. Impose the following conditions on $\br$ and $\ell$:
\begin{enumerate}
\item Let $G^{\mr{sc}} \to G$ denote the simply-connected cover of $G$. Assume there is a subfield $k' \subset k$ such that the image $\br(\Gal(\Qb/\Q))$ contains $\im \left( G^{\mr{sc}}(k') \to G(k') \right)$.
\item $\ell > 1+ \max(8 , 2h-2)\cdot \#Z_{G^{\mr{sc}}}$, where $h$ denotes the Coxeter number of $G$.
\item $\br$ is odd, i.e. for a choice of complex conjugation $c$, $\Ad(\br(c))$ is a split Cartan involution of $G$ (see \S \ref{archcond}).
\item For all places $v \in \Sigma$ not dividing $\ell \cdot \infty$, $\br|_{\gal{F_v}}$ satisfies a liftable local deformation condition $\mc{P}_v$ with tangent space of dimension $\dim_k \left( \mf{g}^{\br(\gal{F_v})}\right)$ (eg, the conditions of \S \ref{steinberg} or \S \ref{minimal}).
\item $\br|_{\gal{\Ql}}$ is ordinary, satisfying the conditions (REG) and (REG*), in the sense of \S \ref{ordsection}.
\end{enumerate}
Then there exists a finite set of primes $Q$ disjoint from $\Sigma$, and a lift 
\[
\xymatrix{
& G(\mc{O}) \ar[d] \\
\Gal(\Qb/\Q) \ar[r]_{\br} \ar@{-->}[ur]^{\rho} & G(k) 
}
\]
such that $\rho$ is type $\mc{P}_v$ at all $v \in \Sigma$ (taking $\mc{P}_v$ to be an appropriate ordinary condition at $v \vert \ell$), of Ramakrishna type (see \S \ref{ramdefsection}) at all $v \in Q$, and unramified outside $\Sigma \cup Q$.  Moreover, we can arrange that $\rho|_{\gal{\Ql}}$ is de Rham, and hence that $\rho$ is geometric in the sense of Fontaine-Mazur.
\end{thm}
The statement of this theorem takes for granted the study of certain local deformation conditions; in \S \ref{localsection}, we study a few of the possibilities, but we have largely ignored this problem as unnecessary for our eventual application. We should note, however, especially for a reader familiar with Kisin's improvement of the Taylor-Wiles method, that Ramakrishna's method, in contrast, seems to require the local deformation rings to be formally smooth. Regarding the global hypothesis in Theorem \ref{maximalintro},
note that for two-dimensional $\br$, a case-by-case treatment of the possible images $\br(\gal{\Q})$ can be undertaken to establish the general case. Such an approach would be prohibitive in general, so after establishing Theorem \ref{maximalintro}, we focus on $\br$ suited to this paper's principal application, the construction of geometric Galois representations
\[
\rho \colon \Gal(\Qb/\Q) \to {}^L G(\Qlb)
\]
with Zariski-dense image in ${}^L G$, where $G$ is one of the exceptional groups $\mr{G}_2$, $\mr{F}_4$, $\mr{E}_6$, $\mr{E}_7$, or $\mr{E}_8$. Let us recall some of the history of this problem.

In his article \cite{serre:motivicgalois}, Serre raised the question of whether there are motives (over number fields, say), whose motivic Galois groups are equal to the exceptional group $\mr{G}_2$, of course implicitly raising the question for other exceptional types as well. If `motive' is taken to mean either
\begin{itemize}
\item pure homological motive in the sense of Grothendieck, but assuming the Standard Conjectures (\cite{kleiman:algcycles}); or
\item motivated motive in the sense of Andr\'{e} (\cite{andre:motivated}),
\end{itemize}
then Dettweiler and Reiter (\cite{dettweiler-reiter:rigidG2}) answered Serre's question affirmatively for the group $\mr{G}_2$, using Katz's theory of rigid local systems on punctured $\mathbb{P}^1$ (\cite{katz:rls})--in particular Katz's remarkable result that all irreducible rigid local systems are suitably motivic. Then, in an astounding development, Yun (\cite{yun:exceptional}) answered a somewhat weaker form of Serre's question for the exceptional types $\mr{G}_2$, $\mr{E}_7$, and $\mr{E}_8$. Namely, he showed that there are motives in the above sense whose associated $\ell$-adic Galois representations have algebraic monodromy group (i.e. the Zariski-closure of the image) equal to these exceptional groups. Yun's work is also deeply connected to the subject of rigid local systems, but the relevant local systems are constructed not via Katz's work, but as the eigen-local systems of suitable `rigid' Hecke eigensheaves on a moduli space of $G$-bundles with carefully-chosen level structure on $\mathbb{P}^1-\{0, 1, \infty\}$. 

In particular, Yun produces the first examples of geometric Galois representations with exceptional monodromy groups $\mr{E}_7$ and $\mr{E}_8$. The main theorem of the present paper is the construction of geometric Galois representations with monodromy group equal to any of the exceptional types: 
\begin{thm}[see Theorem \ref{notE6theorem}, Theorem \ref{E6theorem}]\label{mainintro}
There is a density one set of rational primes $\ell$ such that for each exceptional Lie type $\Phi$, and for a suitably-chosen rational form $G/\Q$ of $\Phi$, there are $\ell$-adic representations
\[
\rho_{\ell} \colon \Gal(\Qb/\Q) \to {}^L G(\Qlb)
\] 
with Zariski-dense image. For $G$ of types $\mr{G}_2$, $\mr{F}_4$, $\mr{E}_7$, or $\mr{E}_8$, i.e. the exceptional groups whose Weyl groups contain $-1$, we can replace ${}^L G$ simply by the dual group $G^\vee$. For $G$ of type $\mr{E}_6$, the algebraic monodromy group of $\rho_{\ell}|_{\Gal(\Qb/K)}$ is $G^\vee$ for an appropriate quadratic imaginary extension $K/\Q$.
\end{thm}
We achieve this via a quite novel method, and indeed the examples we construct are disjoint from those of Dettweiler-Reiter and Yun. The case of $\mr{E}_6$ should stand out here, as it has proven especially elusive: for example, in the paper \cite{heinloth-ngo-yun:kloosterman}, which served as much of the inspiration for Yun's work, certain ${}^L G$-valued $\ell$-adic representations of the absolute Galois group of the function field $\mathbb{F}_q(\mathbb{P}^1)$ were constructed for $G$ of any simple type; their monodromy groups were all computed, and, crucially, the monodromy group turns out to be `only' $\mr{F}_4$ in the case $G= \mr{E}_6$.\footnote{The same thing happens in \cite{frenkel-gross:rigid}, which served as inspiration for \cite{heinloth-ngo-yun:kloosterman}.} I believe the present paper contains the first sighting of the group $\mr{E}_6$ as any sort of arithmetic monodromy group. 

In the rest of this introduction, we will sketch the strategy of Theorem \ref{mainintro}. Let $\gal{\Q}= \Gal(\Qb/\Q)$ (in what follows, $\Q$ can be replaced by any totally real field $F$ for which $[F(\zeta_{\ell}):F]= \ell -1$). The essential content of the argument is already present in the following somewhat simpler case, where as in Theorem \ref{maximalintro} we restrict to the consideration of  
\[
\br \colon \gal{\Q} \to G(k),
\]
for $G$ a simple Chevalley group of adjoint type.\footnote{So what follows will literally apply except in type $\mr{E}_6$; type $\mr{E}_6$ turns out to require a minor, merely technical, modification, carried out in \S \ref{defL}-\ref{E6section}.} The hope would be to start with an appropriate $\br$, and to use Ramakrishna's method to deform it to characteristic zero. But it is already difficult to produce such $\br$, and finding one with image containing $G(\mathbb{F}_{\ell})$, so that the initial version of our lifting theorem might apply, would pose a problem already as difficult\footnote{In fact, substantially more so, because we would need not only a Galois extension with group $G(\mathbb{F}_{\ell})$, but also to know that the associated representation $\br$ satisfied the various technical hypotheses of the lifting theorem, eg ordinarity at $\ell$.} as the inverse Galois problem for the group $G(\mathbb{F}_{\ell})$.

But any (adjoint, for simplicity) $G$ admits a principal homomorphism $\varphi \colon \mr{PGL}_2 \to G$ (see \cite{gross:principalsl2}, \cite{serre:principalsl2}); for instance, for the classical groups, these are the usual symmetric power representations of $\mr{PGL}_2$. We have at our disposal a very well-understood collection of 2-dimensional representations $\bar{r} \colon \gal{\Q} \to \mr{GL}_2(k)$, those associated to holomorphic modular forms, and for our $\br$ we consider composites
\[
\xymatrix{
\gal{\Q} \ar[r]^-{\bar{r}} \ar@/^2 pc/[rr]^{\br} & \mr{PGL}_2(k) \ar[r]^{\varphi} & G(k).
} 
\]
We therefore undertake to prove a version of Ramakrishna's lifting theorem that applies when $\br$ factors through a principal $\mr{PGL}_2$. This is carried out in \S \ref{principalSL2}, buttressed by an axiomatic version of Ramakrishna's argument (\S \ref{formalglobal}). Choosing $\bar{r}$ so that $\br= \varphi \circ \bar{r}$ satisfies the hypotheses of the lifting Theorem \ref{principalsl2lift}, and so that the resulting lift $\rho \colon \gal{\Q} \to G(\mc{O})$ can be guaranteed to have maximal algebraic monodromy group, is rather delicate, however: in particular, we don't want to lift $\bar{r}$ to $r \colon \gal{F} \to \mr{GL}_2(\mc{O})$ and then take $\rho= \varphi \circ r$! The rough idea for ensuring that the monodromy group $G_{\rho}= \overline{\rho(\gal{\Q})}^{\mr{Zar}}$ is equal to $G$ is the following:
\begin{itemize}
\item Ensure $G_{\rho}$ is reductive and contains a regular unipotent element of $G$. There is a straightforward (to use, not to prove) classification of such subgroups of $G$ (\cite[Theorem A]{saxl-seitz:dynkin}). Namely, when $G$ is of type $\mr{G}_2$, $\mr{F}_4$, $\mr{E}_7$, or $\mr{E}_8$, such a $G_{\rho}$ is either a principal $\mr{PGL}_2$ or all of $G$; and when $G$ is of type $\mr{E}_6$, $G_{\rho}$ might also be $\mr{F}_4$.
\item Arrange that the Hodge-Tate weights of $\rho$ are `sufficiently generic' that $G_{\rho}$ must be $G$ itself. This is where our examples veer away from those of \cite{dettweiler-reiter:rigidG2} and \cite{yun:exceptional}.
\end{itemize}
If $\bar{r}(\gal{\Q})$ contains $\mr{PSL}_2(\mathbb{F}_{\ell})$, it is essentially formal that $G_{\rho}$ is reductive: see Lemma \ref{reductivemono}. It is not necessarily true, however, that $G_{\rho}$ contains a principal $\mr{PGL}_2$: indeed, a beautiful result of Serre (deforming the finite groups $\mr{PSL}_2(\mathbb{F}_{\ell})$ themselves to characteristic zero: see Example \ref{serreeg}) shows it is sometimes possible for $G_{\rho}$ to be finite! 

We now explain the subtleties in choosing a modular form $f$, with associated compatible system 
\[
r_{f, \lambda} \colon \gal{\Q} \to \mr{GL}_2(\mc{O}_{\lambda}),
\]
so that the reductions $\bar{r}_{f, \lambda}$, at least for $\lambda$ lying over a density one set of rational primes $\ell$, can be used as our $\bar{r}$. If we were content with establishing Theorem \ref{mainintro} for only a single $\ell$, it seems likely that a significant computer calculation of a well-chosen residual Galois representation for a particular modular form might suffice; but any result for infinitely many $\ell$ (or in our case, a density one set) seems to require hurdling further theoretical difficulties.
\begin{enumerate}
\item Needing $\rho$ to have `sufficiently generic' Hodge-Tate weights forces us to work with ordinary deformations; in particular $\bar{r}_{f, \lambda}$ must be ordinary. Except for $f$ of weights 2 and 3, establishing ordinarity of $f$ even for infinitely many $\ell$ is a totally open problem, so we are forced into weights 2 or 3. Requiring ordinarity restricts the conclusion of Theorem \ref{mainintro} to a density one set of primes rather than almost all primes.
\item While it would be possible to work with $f$ of weight two, the ordinary deformation rings in this case are more often singular, and some extra work is needed to avoid this possibility; in the spirit of minimizing the amount of local work, we therefore take $f$ to have weight 3. 
\item Keeping in mind Serre's cautionary example (Example \ref{serreeg}), we use a local analysis to ensure that $G_{\rho}$ contains a regular unipotent element. The idea is to choose $f$ with $\Gamma_0(p)$ level at some prime $p$, and to consider deformations of `Steinberg' type of $\varphi \circ \bar{r}_{f, \lambda}|_{\gal{\Q_p}}$. As long as the resulting deformations to characteristic zero are as ramified as possible, they will provide us with regular unipotent elements in $G_{\rho}$. But it is quite difficult purely using deformation theory to guarantee that these characteristic zero lifts have regular unipotent ramification unless the residual representations $\bar{r}_{f, \lambda}|_{\gal{\Q_p}}$ are themselves ramified: the analogous issue in studying `lifts of prescribed type' using potential automorphy theorems is only treated by invoking settled cases of the Ramanujan conjecture.\footnote{Note that in the case $G= \mr{G}_2$, the composition of the principal $\mr{SL}_2$ with the quasi-minuscule representation $\mr{G}_2 \into \mr{GL}_7$ remains irreducible, so that potential automorphy techniques could be applied in this case. This approach does not work for the other exceptional groups.} Thus, we need $\bar{r}_{f, \lambda}|_{\gal{\Q_p}}$ to be ramified for almost all $\lambda$. If $f$ were associated to an elliptic curve over $\Q$ with multiplicative reduction at $p$, then the theory of the Tate curve would imply this, but it seems to be quite a deep result in general: we establish it (Proposition \ref{levellower}, which Khare has pointed out was previously proven by Weston in \cite{weston:unobstructed}) using essentially the full strength of level-lowering results for classical modular forms.
\end{enumerate}
Having juggled the demands of the lifting theorem as just described, we can then look in tables of modular forms to find plenty of $f$ that do in fact serve our purpose, or even show by theoretical means that infinitely many can be found. For the final steps, see Theorem \ref{notE6theorem}.

In this sketch we have omitted the case of $\mr{E}_6$. When the Weyl group of $G$ does not contain $-1$, $G$ contains no order 2 element inducing a split Cartan involution, and so there are no `odd' representations $\gal{\Q} \to G(k)$. We instead deform odd representations valued in the L-group of a suitable outer form of $G$. We develop, only in the degree of generality needed for our application, the basics of deformation theory for L-groups in \S \ref{defL}; our task is made easy by the template provided in \cite[\S 2]{clozel-harris-taylor}, which treats the case of type $\mr{A}_n$. With these foundations in place, there are no new difficulties in extending the arguments of earlier sections; we explain the very minor modifications needed in \S \ref{E6section}. The proof of Theorem \ref{mainintro} is then completed in Theorem \ref{E6theorem}. We hope the reader does not object to this expository decision: it would have been possible to work throughout with (possibly non-connected) L-groups, but I think as written the argument will be easier to digest, and in the end it only costs us a few extra pages. 
\section{Notation}\label{notation}
For a field $F$ (always a number field or local field), we let $\gal{F}$ denote $\Gal(\overline{F}/F)$ for some fixed choice of algebraic closure $\overline{F}$ of $F$. When $F$ is a number field, for each place $v$ of $F$ we fix once and for all embeddings $\overline{F} \into \overline{F}_v$, giving rise to inclusions $\gal{F_v} \into \gal{F}$. If $\Sigma$ is a (finite) set of places of $F$, we let $\Gamma_{F, \Sigma}$ denote $\Gal(F_{\Sigma}/F)$, where $F_{\Sigma}$ is the maximal extension of $F$ in $\overline{F}$ unramified outside of $\Sigma$. In \S \ref{defL} and \S \ref{E6section}, for an extension $\tF/F$ we will also write $\gal{\tF, \Sigma}$, where $\Sigma$ is implicitly interpreted as the set of place of $\tF$ above $\Sigma$. If $v$ is a place of $F$ outside $\Sigma$, we write $\fr_v$ for the corresponding geometric frobenius element in $\gal{F, \Sigma}$. If $F$ is a local field, with no reference to a global field, we write $\fr_F$ for a choice of geometric frobenius element in $\gal{F}$. For a representation $\rho$ of $\gal{F}$, we let $F(\rho)$ denote the fixed field of the kernel of $\rho$.

Consider a group $\Gamma$, a ring $R$, an affine group scheme $H$ over $\Spec A$, and a homomorphism $\rho \colon \Gamma \to H(A)$. Then for any algebraic representation $V$ of $H$, we let $\rho(V)$ denote the $A[\Gamma]$-module with underlying $A$-module $V$ induced by $\rho$. This will typically be applied to the adjoint representation of $H$.

Let $\mc{O}$ be a finite totally ramified extension of the ring $W(k)$ of Witt vectors of an algebraic extension $k$ of $\fl$, and let $E$ be the fraction field $\Frac \mc{O}$. We let $\CO^f$ denote the category of artinian local $\mc{O}$-algebras for which the structure map $\mc{O} \to R$ induces an isomorphism on residue fields, and let $\CO$ denote the category of complete local noetherian $\mc{O}$-algebras with residue field $k$. Let $\varpi$ denote a uniformizer of $\mc{O}$.

All the (Galois) cohomology groups we consider will be $k$-vector spaces, and we will always abbreviate $\dim_k H^i(\bullet)$ by $h^i(\bullet)$.

We write $\kappa$ for the $\ell$-adic cyclotomic character, and $\kbar$ for its mod $\ell$ reduction.

Until \S \ref{defL}, $G$ will be a Chevalley group scheme over $\mc{O}$; we mean this in the sense of \cite{conrad:luminy}, so $G$ is not necessarily semi-simple. The reader will not lose anything essential by taking $G$ to be adjoint and considering the classical construction of Chevalley groups as in \cite{steinberg:chevalley}. We refer to \cite{conrad:luminy} for a thorough and accessible treatment of the theory of reductive group schemes; this reference is vastly more general than we require, but it still seems to be the most convenient. Throughout the present paper, we will give as needed more specific pointers to results in \cite{conrad:luminy}, but a reader with additional questions will surely find them answered there as well.
\section{Review of deformation theory}\label{review}
In this section we establish our conventions and notation for the deformation theory of Galois representations. Although we could work much more generally, considering representations valued in an arbitrary connected reductive $\mc{O}$-group scheme $G$ (compare \cite[\S 2]{tilouine:defs}), for simplicity we restrict as in \S \ref{notation} and require that $G$ be a connected reductive Chevalley group scheme. In \S \ref{defL}, following the example of \cite[\S 2.2]{clozel-harris-taylor}, we will recast this background to allow representations valued in certain non-connected L-groups, but for the bulk of the paper, and most cases of the main Theorem \ref{mainintro}, the present discussion suffices.
We write $\fg$ for the Lie algebra of $G$; we will abuse notation and continue to write $\mf{g}$ for the base-change to various coefficient rings, most notably the special fiber. 

Some preliminary hypotheses on the prime $\ell$ are also needed. In the central results of this paper (eg Theorem \ref{principalsl2lift}), we will impose somewhat stricter requirements, but for now it suffices to take $\ell \neq 2$ to be a `very good prime' for all simple factors of $G$: see \cite[\S 1.14]{carter:finitelie}.\footnote{In particular, taking $\ell  \nmid n+1$ in type $\mr{A}_n$ and $\ell \geq 7$ in all other cases suffices.} Here are the relevant consequences:
\begin{itemize}
\item The isogeny theorem for root data (\cite[Theorem 6.1.16]{conrad:luminy}) yields two canonical central isogenies, $Z_G^0 \times G^{\mr{der}} \to G$ and $G \to G/G^{\mr{der}} \times G/Z_G^0$, with $Z_G^0$ the maximal central torus and $G^{\mr{der}}$ the derived group of $G$. These are isomorphisms at the level of Lie algebras, and note that $\mf{g}^{\mr{der}}$ does not depend up to isomorphism on the isogeny class of $G^{\mr{der}}$: if $\ell$ is very good, then it does not divide the determinant of the Cartan matrix. From now on, we write
\begin{equation}\label{multiplier}
\mu \colon G \to S= G/G^{\mr{der}}
\end{equation}
for this map onto the maximal quotient torus of $G$ and denote $G^{\mr{der}}$ by $G_{\mu}$. Likewise, the Lie algebra $\fg^{\mr{der}}$ will be denoted $\fg_{\mu}$.\footnote{In some of the discussion that follows, the reader could replace $\mu \colon G \to S$ by some other map to an $\mc{O}$-torus, whose kernel may be bigger than the derived group.}
\item In particular, $\fg= \mf{z}(\fg) \oplus \fg_{\mu}$, where $\mf{z}(\fg)$ is the center of the Lie algebra of $G$ (and the Lie algebra of $Z_G$), and $\fg_{\mu}$ is an irreducible representation; the latter statement is checked case-by-case on simple types. 
\item There is a non-degenerate $G$-invariant form $\fg_{\mu} \times \fg_{\mu} \to k$: use \cite[\S 1.16]{carter:finitelie} for types other than $\mr{A}_n$ and the Killing form in type $\mr{A}_n$ (whose discriminant is divisible only by primes dividing $2(n+1)$). We fix such a pairing and throughout make the resulting identification $\fg_{\mu}^* \cong \fg_{\mu}$ (this will come up when studying dual Selmer groups). Note that the Killing form of a simple Lie algebra (in very good characteristic, so we don't have to specify the isogeny class) is non-degenerate as long as $\ell$ does not divide 2, the discriminant of the Cartan matrix, the dual Coxeter number, and the ratio of long and short roots (\cite[Proposition I.4.8]{springer-steinberg:conjugacy}). Since in our main theorems we will require $\ell$ to be even larger, most readers will lose nothing by making this stronger assumption from the outset.
\end{itemize}

\subsection{The basics}\label{defbasics}
Let $\Gamma$ be a profinite group, and fix a continuous homomorphism $\br \colon \Gamma \to G(k)$. For simplicity (allowing the restriction to noetherian coefficient algebras), we assume that $\Gamma$ satisfies the $\ell$-finiteness condition of \cite{mazur:deforming}: for all open subgroups $\Gamma_0 \subset \Gamma$, there are only finitely many continuous homomorphisms $\Gamma_0 \to \Z/\ell \Z$. Noetherian hypotheses can be avoided, as in \cite[\S 2.2]{clozel-harris-taylor}, but we at least would gain nothing, and can save some work, by imposing them. We recall the basic definitions of Mazur's deformation theory:
\begin{defn}\label{liftfunctor}
\begin{itemize}
\item Denote by 
\[
\Lift_{\br} \colon \mc{C}_{\mc{O}} \to \Sets
\]
the functor whose $R$-points is the set of lifts of $\br$ to a continuous homomorphism $\Gamma \to G(R)$. It is easy to see that $\Lift_{\br}$ is representable, and we denote its representing object, the universal lifting ring, by $R^{\square}_{\br}$.
\item We say that two lifts $\rho_1, \rho_2 \colon \Gamma \to G(R)$ of $\br$ are strictly equivalent if they are conjugate by an element of
\[
\widehat{G}(R)= \ker\left( G(R) \to G(k) \right).
\]
The functor $\widehat{G}$ is represented by a smooth group scheme over $\mc{O}$.
\item Denote by
\[
\Def_{\br} \colon \CO \to \Sets
\]
the functor assigning to $R$ the set of strict equivalence classes of elements of $\Lift_{\br}(R)$. A \textit{deformation} of $\br$ is an element of $\Def_{\br}(R)$.
\end{itemize}
\end{defn}
As usual, we will need to study certain representable sub-functors of $\Lift_{\br}$. We can always initially define these sub-functors only on $\CO^f$ and then, having proven (pro-)representability, extend them to $\CO$ by `continuity;' but often one wants a `moduli-theoretic' description on all of $\CO$, in which case one might define a functor on $\CO$ and verify directly that it commutes with filtered limits. We will allow ourselves a minor abuse of terminology and allow a `sub-functor of $\Lift_{\br}$' to refer to either of these cases. Similarly, we will allow ourselves to write `representable' when what is strictly speaking meant is `pro-representable.'
\begin{defn}\label{defcond}
A \textit{deformation condition} is a representable sub-functor $\Lift_{\br}^{\mc{P}}$ of $\Lift_{\br}$ that is closed under strict equivalence. We denote the representing object by $R_{\br}^{\square, \mc{P}}$; it is canonically the quotient of $R_{\br}^\square$ by some $\widehat{G}(R_{\br}^\square)$-invariant ideal $J^{\mc{P}}$ of $R_{\br}^\square$.\footnote{By the universal property, conjugation of the universal lift by any element of $\widehat{G}(R^{\square}_{\br})$ induces a morphism $R_{\br}^\square \to R_{\br}^\square$.}
\end{defn}
Schlessinger's criterion gives one way to cut down the amount of work needed to check that a sub-functor is a deformation condition:
\begin{lemma}[Schlessinger]\label{representability}
Let $\Lift_{\br}^{\mc{P}}$ be a sub-functor of $\Lift_{\br}$, assumed to be closed under strict equivalence. Then $\Lift_{\br}^{\mc{P}}$ is a deformation condition if and only if for all morphisms $A \to C$, $B \to C$ in $\CO^f$ with $B \to C$ small, the natural map
\[
\Lift_{\br}^{\mc{P}}(A \times_C B) \to \Lift_{\br}^{\mc{P}}(A) \times_{\Lift_{\br}^{\mc{P}}(C)} \Lift_{\br}^{\mc{P}}(B)
\]
is surjective.
\end{lemma}
\proof
This follows immediately from \cite[Theorem 2.11]{schlessinger:functors}: since $\Lift_{\br}^{\mc{P}}$ is a sub-functor of $\Lift_{\br}$, all of the injectivity statements in Schlessinger's criteria follow from the corresponding statements for $\Lift_{\br}$.
\endproof  
Here is an important example, the general analogue of fixing the determinant in the case $G= \mr{GL}_N$:
\begin{eg}\label{fixedchar}
Recall from Equation (\ref{multiplier}) that $\mu \colon G \to S$ is the map onto the maximal abelian quotient. Fix a lift $\nu \colon \gal{F} \to S(\mc{O})$ of $\mu \circ \br$; for all $\mc{O}$-algebras $R$, we also write $\nu$ for the induced homomorphism $\gal{F} \to S(R)$. Then we consider the sub-functor $\Lift_{\br}^{\nu}$ of $\Lift_{\br}$ of lifts $\rho \colon \gal{F} \to G(R)$ such that $\mu \circ \rho= \nu$. This is a deformation condition in the sense of Definition \ref{defcond}.
\end{eg}
We now recall the usual description of the tangent space of the deformation functor (the proofs in what follows are standard, or can be imitated from the case of $G= \mr{GL}_N$ in \cite{clozel-harris-taylor}). Without the crutch of matrices, we will instead use the exponential map for $G$, described for instance in \cite[\S 3.5]{tilouine:defs}: for any small extension $A \to B$ in $\CO^f$, with kernel $I$, the exponential map is a bijection
\[
\exp \colon \fg \otimes_k I \xrightarrow{\sim} \ker \left(\widehat{G}(A) \to \widehat{G}(B) \right).
\]
There are canonical isomorphisms
\begin{equation}\label{tspace1}
\Hom_k \left(\mf{m}_{R^{\square}_{\br}}/(\mf{m}_{R^{\square}_{\br}}^2, \varpi), k \right) \cong \Hom_{\CO}(R^{\square}_{\br}, k[\epsilon]) \cong \Lift_{\br}(k[\epsilon]) \xleftarrow[\sim]{\tau} Z^1(\Gamma, \br(\fg)).
\end{equation}
The isomorphism $\tau$ associates to a cocycle $\phi \in Z^1(\Gamma, \br(\fg))$ the lift
\[
g \mapsto \rho(g)= \exp(\epsilon \phi(g))\br(g) \in G(k[\epsilon]).
\]
It also induces $\Def_{\br}(k[\epsilon]) \cong H^1(\Gamma, \br(\fg))$. We also have the usual variant with fixed similitude character as in Example \ref{fixedchar}: $\Lift_{\br}^{\nu}(k[\epsilon]) \cong Z^1(\Gamma, \br(\mf{g}_{\mu}))$. 

Now suppose $\Lift_{\br}^{\mc{P}}$ is a deformation condition, represented by $R_{\br}^{\square, \mc{P}} \xleftarrow{\sim} R^{\square}_{\br}/J^{\mc{P}}$. We associate a subspace $L^{\square, \mc{P}} \subset Z^1(\Gamma, \br(\fg))$, with image $L^{\mc{P}} \subset H^1(\Gamma, \br(\fg))$, with the property that $\Lift^{\mc{P}}_{\br}(k[\epsilon]) \xleftarrow[\tau]{\sim} L^{\square, \mc{P}}$, as follows: $L^{\square, \mc{P}}$ is by definition the annihilator of $J^{\mc{P}}/(J^{\mc{P}}\cap(\mf{m}_{R^{\square}_{\br}}^2, \varpi))$ under the pairing induced by Equation (\ref{tspace1}). Since $\Lift_{\br}^{\mc{P}}$ is closed under strict equivalence, there is an exact sequence
\[
0 \to H^0(\Gamma, \br(\fg)) \to \fg \to L^{\square, \mc{P}} \to L^{\mc{P}} \to 0,
\]
i.e. $L^{\square, \mc{P}}$ contains all coboundaries.

Beyond simply describing the set of deformations to $k[\epsilon]$, the subspaces $L^{\mc{P}}$ are useful in describing deformations across any small morphism $R \to R/I$ in $\CO^f$ (`small' means that $\mf{m}_R \cdot I=0$; in particular, $I$ is a $k$-vector space). Namely, if $\rho \in \Lift_{\br}^{\mc{P}}(R/I)$, then the set of lifts of $\rho$ to an element of $\Lift_{\br}^{\mc{P}}(R)$ is an $L^{\square, \mc{P}} \otimes_k I$-torsor. This torsor may of course be empty, and therein lies the whole difficulty of the subject. At least one knows that the obstruction to lifting $\rho$ to an element of $\Lift_{\br}(R)$ (not necessarily satisfying the deformation condition $\mc{P}$) is measured by a class in $H^2(\gal{F}, \br(\fg)) \otimes_k I$. 
\begin{defn}
We say a deformation condition $\Lift_{\br}^{\mc{P}}$ is \textit{liftable} if for all small surjections $R \to R/I$, 
\[
\Lift_{\br}^{\mc{P}}(R) \to \Lift_{\br}^{\mc{P}}(R/I)
\]
is surjective. Equivalently, $R_{\br}^{\square, \mc{P}}$ is isomorphic to a power series ring over $\mc{O}$ in $\dim_k L^{\square, \mc{P}}$ variables.
\end{defn}
Finally, we remark that it is often convenient to define a local deformation condition after replacing $\mc{O}$ by the ring of integers $\mc{O}'$ of some finite extension of $\mr{Frac}(\mc{O})$; there are various ways of handling this, the simplest being just to enlarge, from the outset, $\mc{O}$ (and $k$) to be as big as necessary, and then to work with this updated version of $\mc{C}^f_{\mc{O}}$ (the reader who does not wish to keep track of successive enlargements can once and for all take $\mc{O}= W(\overline{\mathbb{F}}_{\ell})$).
\subsection{The global theory with local conditions}\label{deflocglob}
Again we take $\Gamma$ to be a profinite group, but we also assume it is equipped with maps, indexed by $v$ in some set $\Sigma$, $\iota_v \colon \Gamma_v \to \Gamma$ from profinite groups $\Gamma_v$ (also satisfying the $\ell$-finiteness condition). We continue to fix a continuous homomorphism $\br \colon \Gamma \to G(k)$, but now the discussion from \S \ref{defbasics} applies both to $\br$ and to its `restrictions' $\br_v= \br \circ \iota_v$. For each $v \in \Sigma$ we give ourselves (`local') deformation conditions $\Lift_{\br_v}^{\mc{P}_v}$, and then consider $\Lift_{\br}^{\mc{P}}$, the functor of lifts $\rho \in \Lift_{\br}(R)$ such that $\rho|_{\Gamma_v} \in \Lift_{\br_v}^{\mc{P}_v}(R)$ for all $v \in \Sigma$; it is of course also representable, by a ring we denote $R_{\br}^{\square, \mc{P}}$. We now formulate the analogue of the tangent space and obstruction theories in this local-global setting. Let $L_v^{\square}$ denote the subspace of $Z^1(\Gamma_v, \br(\fg))$ corresponding to $\mc{P}_v$, with image $L_v$ in $H^1(\Gamma_v, \br(\fg))$. Define
\begin{equation}\label{globtangent}
L_v^{\square, i}= 
\begin{cases}
\text{$C^0(\Gamma_v, \br(\fg))$ if $i=0$;}\\
\text{$L_v^{\square}$ if $i=1$;} \\
\text{$0$ if $i \geq 2$.}
\end{cases}
\end{equation}
Then consider the map of co-chain complexes given by restriction to $\Gamma_v$ for all $v \in \Sigma$:
\[
C^{\bullet}(\Gamma, \br(\fg)) \xrightarrow{f} \bigoplus_{v \in \Sigma} C^{\bullet}(\Gamma_v, \br(\fg))/L_v^{\square, \bullet}.
\]
Define a complex $C^i_{\mc{P}}(\Gamma, \br(\fg))= \mr{Cone}_f^{i-1}$ (for some choice of cone). We denote cocycles and cohomology of the complex $C^\bullet_{\mc{P}}$ by $Z^\bullet_{\mc{P}}$ and $H^\bullet_{\mc{P}}$, so there is a long exact sequence in cohomology
\begin{equation}\label{les}
0 \to H^1_{\mc{P}}(\Gamma, \br(\fg)) \to H^1(\Gamma, \br(\fg)) \to \bigoplus_{v \in \Sigma} H^1(\Gamma_v, \br(\fg))/L_v 
\to H^2_{\mc{P}}(\Gamma, \br(\fg)) \to \ldots,
\end{equation}
and it is immediate that $\Lift_{\br}^{\mc{P}}(k[\epsilon]) \cong Z^1_{\mc{P}}(\Gamma, \br(\fg)$. If the invariants $\fg^{\br(\Gamma)}$ equal the center $\mf{z}(\fg)$, then $\Def_{\br}$ and the corresponding $\Def_{\br}^{\mc{P}}$ are representable (by Schlessinger's criteria--see the proof of Proposition \ref{defLprop}), by objects $R_{\br}$ and $R_{\br}^{\mc{P}}$ of $\CO$, and then
\[
\Def_{\br}^{\mc{P}}(k[\epsilon]) \cong H^1_{\mc{P}}(\Gamma, \br(\fg)).
\]
The essential result is the following; the argument goes back to Mazur's original article \cite{mazur:deforming}.
\begin{prop}\label{globalobstruction}
Assume that $\fg^{\br(\Gamma)}= \mf{z}(\fg)$, and let $\Lift_{\br_v}^{\mc{P}_v}$ be a collection of \textit{liftable} local deformation conditions. Then the universal deformation ring $R_{\br}^{\mc{P}}$ is isomorphic to a quotient of a power series ring over $\mc{O}$ in $\dim_k H^1_{\mc{P}}(\Gamma, \br(\fg))$ variables by an ideal that can be generated by at most $\dim_k H^2_{\mc{P}}(\Gamma, \br(\fg))$ elements. In particular, if $H^2_{\mc{P}}(\Gamma, \br(\fg))=0$, then $R_{\br}^{\mc{P}}$ is isomorphic to a power series ring over $\mc{O}$ in $\dim_k H^1_{\mc{P}}(\Gamma, \br(\fg))$ variables.
\end{prop}
\begin{rmk}
Liftability of the local deformation conditions is used to define classes in $H^2_{\mc{P}}(\Gamma, \br(\fg)) \otimes_k I$ measuring the obstructions to surjectivity of $\Def_{\br}^{\mc{P}}(R) \to \Def_{\br}^{\mc{P}}(R/I)$, for any small surjection $R \to R/I$.
\end{rmk}
\begin{rmk}\label{simchar}
The whole discussion of this subsection continues to hold if we fix a similitude character $\nu \colon \Gamma \to S(\mc{O})$ as in Example \ref{fixedchar}, simply replacing $\fg$ by $\fg_{\mu}$ in all the group cohomology calculations; in order to have a transparent relationship between the two problems--with fixed similitude character and without--we use our ($\ell$ very good) assumption that $\mu \colon G \to S$ splits at the level of Lie algebras: $\fg \cong \fg_{\mu} \oplus \mf{s}$, $G$-equivariantly. A liftable local deformation condition $\Lift_{\br}^{\mc{P}_v}$ with tangent space $L_v$ then induces a liftable local deformation condition, now also requiring fixed multiplier character, $\Lift_{\br}^{\mc{P}_v, \nu}$ with tangent space $L_v \cap H^1(\Gamma_v, \brgm)$. If we are fixing similitude characters, then we will use $L_v$ to refer to this intersection, not to the larger tangent space.
\end{rmk}
\subsection{Deformations of Galois representations}
We specialize now to the setting of interest. The results of this subsection will be explained with proof, in a very slightly different (but perhaps less familiar) context, in \S \ref{defL}; see Proposition \ref{defLprop}. Fix a number field $F$ and a continuous representation $\br \colon \gal{F} \to G(k)$ that is unramified outside a finite set of places $\Sigma$, which we will assume contains all places above $\ell$ and $\infty$. Then for $\Gamma$ we will take $\gal{F, \Sigma}$, the Galois group of the maximal extension $F_{\Sigma}/F$ inside $\overline{F}$ that is unramified outside of $\Sigma$. For each $v \in \Sigma$, we consider the groups $\Gamma_v= \gal{F_v}$ with their maps $\gal{F_v} \to \gal{F, \Sigma}$ (enshrined in \S \ref{notation}). Note that $\gal{F, \Sigma}$ and $\gal{F_v}$ satisfy the $\ell$-finiteness hypothesis, so the discussion of sections \ref{defbasics}-\ref{deflocglob} applies. Fix a lift $\nu \colon \gal{F, \Sigma} \to S(\mc{O})$ of $\mu \circ \br$, and for the remainder of this section require all local and global lifting functors to have fixed multiplier $\nu$. For each $v \in \Sigma$, fix a liftable deformation condition $\mc{P}_v$, yielding a global deformation condition $\mc{P}= \{\mc{P}_v\}_{v \in \Sigma}$, so that we can consider the global functors $\Lift_{\br}^{\mc{P}}$ and $\Def_{\br}^{\mc{P}}$. The former is always representable, and the latter is as long as we assume that the centralizer of $\br$ in $\fg$ is equal to $\mf{z}(\fg)$, which we do from now on.

The obstruction theory of Proposition \ref{globalobstruction} is related to a more tractable problem in Galois cohomology via the following essential result, which follows from combining the long exact sequence of Equation (\ref{les}) with the Poitou-Tate long exact sequence. To state it, let $L_v^{\perp} \subset H^1(\gal{F_v}, \br(\fg_{\mu})(1))$ denote the orthogonal complement of $L_v$ under the local duality pairing (recall we have fixed an identification $\fg_{\mu}^* \cong \fg_{\mu}$), and consider the dual Selmer group 
\[
H^1_{\mc{P}^{\perp}}(\gal{F, \Sigma}, \br(\fg_{\mu})(1))= \ker \left(H^1(\gal{F, \Sigma}, \br(\fg_{\mu})(1)) \to \bigoplus_{v \in \Sigma} H^1(\gal{F_v}, \br(\fg_{\mu})(1))/L_v^{\perp} \right).
\]
Again, see Proposition \ref{defLprop} for proofs of (slight variants of) what follows.
\begin{prop}\label{h2isdualselmer}
$\dim H^2_{\mc{P}}(\gal{F, \Sigma}, \br(\fg_{\mu}))= \dim H^1_{\mc{P}^\perp}(\gal{F, \Sigma}, \br(\fg_{\mu})(1))$.
\end{prop}
We also recall Wiles's formula relating the size of a Selmer group to that of the corresponding dual Selmer group, a consequence of global duality and the global Euler characteristic formula.
\begin{prop}\label{statewilesformula}
Retain the above hypotheses and notation. Then
\begin{align*}
&\dim H^1_{\mc{P}}(\gal{F, \Sigma}, \br(\fg_{\mu}))- \dim H^1_{\mc{P}^{\perp}}(\gal{F, \Sigma}, \br(\fg_{\mu})(1))= \\
&h^0(\gal{F}, \br(\fg_{\mu}))-h^0(\gal{F}, \br(\fg_{\mu})(1))+ \sum_{v \in \Sigma} \left( \dim L_v-h^0(\gal{F_v}, \br(\fg_{\mu})) \right).
\end{align*}
\end{prop}
Of course, we have assumed $h^0(\gal{F}, \brgm)=0$, but the formula holds as stated without this assumption. Propositions \ref{globalobstruction}, \ref{h2isdualselmer}, and \ref{statewilesformula} yield the global cohomological foundation of Ramakrishna's method:
\begin{cor}\label{globdefring}
Let $\br \colon \gal{F} \to G(k)$ be a continuous homomorphism with infinitesimal centralizer equal to $\mf{z}(\fg)$, and let $\Sigma$ be a finite set of places of $F$ containing all primes at which $\br$ is ramified, all archimedean places, and all places above $\ell$. Fix a multiplier character $\nu$ lifting $\mu \circ \br$. For all $v \in \Sigma$, let $\mc{P}_v$ be a liftable local deformation condition with corresponding tangent space $L_v$. If $H^1_{\mc{P}^{\perp}}(\gal{F, \Sigma}, \br(\fg_{\mu})(1))=0$, then $R_{\br}^{\mc{P}}$ is isomorphic to a power series ring over $\mc{O}$ in
\[
\dim H^1_{\mc{P}}(\gal{F, \Sigma}, \br(\fg_{\mu}))= h^0(\gal{F, \Sigma}, \br(\fg_{\mu}))-h^0(\gal{F, \Sigma}, \br(\fg_{\mu})(1))+ \sum_{v \in \Sigma} (\dim L_v - h^0(\gal{F_v}, \br(\fg_{\mu}))
\]
variables.
\end{cor}
\begin{rmk}\label{selmernote}
Recall that for any finite $\gal{F, \Sigma}$-module $M$ and collection of subspaces $L_v \subset H^1(\gal{F_v}, M)$ for $v \in \Sigma$, the Selmer group $H^1_{\{L_v\}}(\gal{F, \Sigma}, M)$ can also be regarded as the subset $H^1_{\{L_v\}}(\gal{F}, M)$ of $H^1(\gal{F}, M)$ consisting of all classes whose restriction to $\gal{F_v}$ lies in $L_v$ for all $v \in \Sigma$, and whose restriction to $\gal{F_v}$ is unramified for all $v \not \in \Sigma$. In particular, for any finite place $w \not \in \Sigma$, if we take $L_w$ to be the unramified local condition, then the canonical inflation map
\[
H^1_{\{L_v\}_{v \in \Sigma}}(\gal{F, \Sigma}, M) \xrightarrow{\sim} H^1_{\{L_v\}_{v \in \Sigma \cup w}}(\gal{F, \Sigma \cup w}, M)
\]
is an isomorphism.
\end{rmk}
\section{Some liftable local deformation conditions}\label{localsection}
In this section we study some local deformation conditions that are particularly useful for the global application to exceptional monodromy groups.
\subsection{Ordinary deformations}\label{ordsection}
Fix a Borel subgroup $B \subset G$, let $N$ denote its unipotent radical, and let $T$ be the torus $B/N$. Let $F$ be a finite extension of $\Ql$, with absolute Galois group $\gal{F}= \Gal(\overline{F}/F)$. In this section we will study the functor of ordinary lifts of a residual representation
\[
\br \colon \gal{F} \to G(k)
\]
that satisfies $\br(\gal{F}) \subset B(k)$. This problem is discussed in \cite[\S 2]{tilouine:defs}, and our arguments will have a quite similar flavor, but will yield precise analogues for a general group $G$ of the results proven in \cite[\S 2.4.2]{clozel-harris-taylor} for the case $G= \mr{GL}_n$. Following Tilouine, we introduce (but do not yet impose) the following hypotheses:
\begin{align*}
\text{(REG)} \qquad &H^0(\gal{F}, \br(\fg/\fb))=0. \\
\text{(REG*)} \qquad &H^0(\gal{F}, \br(\fg/\fb)(1))=0.
\end{align*}
We now define the deformation functor of interest. We can push-forward $\br$ to a homomorphism
\begin{equation}\label{diagchar}
\br_T \colon \gal{F} \to T(k),
\end{equation}
and we once and for all fix a lift 
\[
\chi_T \colon I_F \to T(\mc{O})
\]
of $\br_T|_{I_F}$. For any $\mc{O}$-algebra $R$, we also denote by $\chi_T$ the associated homomorphism to $T(R)$. 
\begin{defn}\label{ordcond} Let
\[
\Lift_{\br}^{\chi_T} \colon \mc{C}^f_{\mc{O}} \to \Sets
\] 
be the subfunctor of $\Lift_{\br}$ whose set of $R$-points $\Lift_{\br}^{\chi_T}(R)$ is given by all lifts $\rho \colon \gal{F} \to G(R)$ of $\br$ such that
\begin{itemize}
\item there exists $g \in \widehat{G}(R)$ such that ${}^g \rho(\gal{F}) \subset B(R)$; and 
\item the restriction to inertia of the push-forward
\[
({}^g \rho)_T|_{I_F} \colon I_F \to B(R) \to T(R)
\]
is equal to $\chi_T$. 
\end{itemize}
We let $R^{\square, \chi_T}_{\br}$ denote the universal lifting ring representing $\Lift^{\chi_T}_{\br}$.
\end{defn}
\begin{lemma}\label{ordcondlem}
Assume $\br$ satisfies (REG). Then $\Lift_{\br}^{\chi_T}$ is well-defined (i.e. the second condition in Definition \ref{ordcond} does not depend on the choice of $g$), and it defines a local deformation condition in the sense of Definition \ref{defcond}.

Moreover, if $\rho \colon \gal{F} \to G(\mc{O})$ is a continuous lift of $\br$ such that $\rho_n= \rho \pmod {\varpi^n}$ belongs to $\Lift_{\br}^{\chi_T}(\mc{O}/\varpi^n)$ for all $n \geq 1$, then $\rho$ is $\widehat{G}(\mc{O})$-conjugate to a homomorphism $\gal{F} \to B(\mc{O})$ whose push-forward to $T(\mc{O})$ equals $\chi_T$ on $I_F$.
\end{lemma}
\proof
Both claims follow from the elementary fact (see \cite[Proposition 6.2]{tilouine:defs}), whose proof requires the assumption (REG), that if for some $g \in G(R)$, both $\rho(\gal{F}) \subset B(R)$ and ${}^g \rho(\gal{F}) \subset B(R)$, then $g \in \widehat{B}(R)$ (see Lemma \ref{ramcond} below for a similar argument). Consequently, the push-forwards $\rho_T$ and $({}^g \rho)_T$ agree, so $\Lift_{\br}^{\chi_T}$ is well-defined, and that it defines a local deformation condition follows, using Lemma \ref{representability}, as in \cite[Proposition 6.2]{tilouine:defs} (or again, see Lemma \ref{ramcond} below).

The second claim follows from the same consequence of (REG). Namely, for all $n \geq 1$, let $g_n \in \widehat{G}(\mc{O}/\varpi^n)$ conjugate $\rho_n$ into $B$, so that, by the previous paragraph, there is some $b_n \in \widehat{B}(\mc{O}/\varpi^n)$ such that
\[
g_{n+1} \pmod {\varpi^n}= b_n g_n.
\]
Lifting $b_n$ to an element $\tilde{b}_n \in \widehat{B}(\mc{O}/\varpi^{n+1})$ and replacing $g_{n+1}$ by $\tilde{b}_n^{-1} g_{n+1}$, we may assume $g_{n+1}$ lifts $g_n$. After inducting on $n$, the resulting $(g_n)_{n \geq 1}$ define the required element of $\widehat{G}(\mc{O})$.
\endproof
We can now describe the tangent space $L_{\br}^{\chi_T}$ of the deformation condition associated to $\Lift_{\br}^{\chi_T}$:
\begin{lemma}\label{ordtanspace}
Assume $\br$ satisfies (REG), so that $\Lift_{\br}^{\chi_T}$ defines a local deformation condition. Then its tangent space is
\[
\ker \left( H^1(\gal{F}, \br(\fb)) \to H^1(I_F, \br(\fb/\fn)) \right).
\]
\end{lemma}
\proof
Note that by the assumption (REG), $H^1(\gal{F}, \br(\fb))$ is a subspace of $H^1(\gal{F}, \br(\fg))$, so the claim of the lemma is meaningful. Let $\phi \in Z^1(\gal{F}, \br(\fg))$ represent a class in the tangent space, so that there exists $X \in \fg$ such that
\[
\exp(\epsilon X) \exp(\epsilon \phi(g)) \br(g) \exp(-\epsilon X) \in B(k[\epsilon]).
\]
Replacing $\phi$ by the cohomologous cocycle $g \mapsto \phi(g)+ X- \ad(\br(g)) X$, we may therefore assume that $\phi$ belongs to $Z^1(\gal{F}, \br(\fb))$. Pushing forward to $T(k[\epsilon])$, we see that $\phi(g)_T|_{I_F}=0$, and the claim follows.
\endproof
We now come to the main result of this section. To achieve a particularly simple description of $R^{\square, \chi_T}_{\br}$, we require additional hypotheses on $\br$; but note that this result implies the result of \cite[\S 2.4.2]{clozel-harris-taylor} in the case $G= \mr{GL}_n$.
\begin{prop}\label{ordliftring}
Assume that $\br$ satisfies (REG) and (REG*), and that $\zeta_{\ell} \notin F$. Then $\Lift^{\chi_T}_{\br}$ is liftable, and the dimension of $L_{\br}^{\chi_T}$ is $\dim_k (\fn)[F: \Ql]+ \dim_k H^0(\gal{F}, \br(\fg))$. That is, $R^{\square, \chi_T}_{\br}$ is a power series ring over $\mc{O}$ in $\dim_k \fg+ \dim_k (\fn)[F: \Ql]$ variables.
\end{prop}
\proof
We must show that $\Lift^{\chi_T}_{\br}$ is formally smooth, and that the tangent space $L_{\br}^{\chi_T}$ has the claimed dimension. Both claims rely on the observation that $H^2(\gal{F}, \br(\fn))=0$, which follows from the assumption (REG*). To see this, note that the Killing form, which is $\gal{F}$-equivariant, induces a non-degenerate pairing $\fn \times \fg/\fb \to k$, hence a $\gal{F}$-equivariant identification $\fn^* \cong \fg/\fb$. Local duality then implies
\[
H^2(\gal{F}, \br(\fn)) \cong H^0\left(\gal{F}, \br(\fg/\fb)(1)\right)^*=0.
\]
We first use this to compute $\dim L^{\chi_T}_{\br}$, the key point being that the map
\begin{equation}\label{tanspace}
f \colon H^1(\gal{F}, \br(\fb)) \to H^1(I_F, \br(\fb/\fn))^{\gal{F}/I_F}
\end{equation}
used in describing $L^{\chi_T}_{\br}$ is surjective. Indeed, it is the composite of maps
\[
H^1(\gal{F}, \br(\fb)) \xrightarrow{\alpha} H^1(\gal{F}, \br(\fb/\fn)) \xrightarrow{\beta} H^1(I_F, \br(\fb/\fn))^{\gal{F}/I_F},
\]
where $\alpha$ is surjective since $H^2(\gal{F}, \br(\fn))=0$; and where $\beta$ is surjective because $\gal{F}/I_F$ has cohomological dimension 1. Now (REG) and (REG*) imply, respectively, that $H^0(\gal{F}, \br(\fb))= H^0(\gal{F}, \br(\fg))$ and $H^2(\gal{F}, \br(\fb))= H^2(\gal{F}, \br(\fb/\fn))$, and then the local Euler characteristic formula and local duality imply that
\begin{align}\label{eqtn1}
h^1(\gal{F}, \br(\fb))-h^0(\gal{F}, \br(\fg))&= [F:\Ql] \dim \fb+ h^2(\gal{F},\br(\fb/\fn))  \\
&= [F:\Ql] \dim \fb+ h^0(\gal{F}, \br(\fb/\fn)(1))= [F:\Ql] \dim \fb;\nonumber
\end{align}
in the final equality we use the hypothesis $\zeta_{\ell} \notin F$ (note that $\br(\fb/\fn)$ is a trivial $\gal{F}$-module). Moreover, by local class field theory
\begin{equation}\label{eqtn2}
h^1(I_F, \br(\fb/\fn))^{\gal{F}/I_F}= h^1(\gal{F}, \br(\fb/\fn))- h^1(\gal{F}/I_F, \br(\fb/\fn))=\dim (\fb/\fn)[F:\Ql],
\end{equation}
and we can combine Equations (\ref{eqtn1}) and (\ref{eqtn2}) to conclude that
\[
\dim \ker(f)- h^0(\gal{F}, \br(\fg))= \dim(\fn)[F:\Ql]. 
\]

Now we show that $\Lift_{\br}^{\chi_T}$ is liftable. Let $R \to R/I$ be a small surjection in $\mc{C}^f_{\mc{O}}$, and let $\rho \in \Lift_{\br}^{\chi_T}(R/I)$. We must lift $\rho$ to an object of $\Lift_{\br}^{\chi_T}(R)$, and since $\widehat{G}$ is smooth, it suffices to do this in the case where $\rho$ factors through $B(R/I)$. The obstruction to lifting $\rho$ to a homomorphism $\rho_R \colon \gal{F} \to B(R)$ lies in $H^2(\gal{F}, \br(\fb)) \otimes_k I$, which is zero as we have already seen $h^2(\gal{F}, \br(\fn))=0$ and $h^0(\gal{F}, \br(\fb/\fn)(1))=0$. We can therefore find such a $\rho_R$, and it remains only to arrange its push-forward to $T$ to equal $\chi_T$. But the space of $B(R)$-valued lifts of $\rho$ is an $H^1(\gal{F}, \br(\fb)) \otimes_k I$-torsor, and the claim follows immediately from the surjectivity of the map $f$ in Equation (\ref{tanspace}).
\endproof
Let us quickly explain how the analogous results with fixed multiplier character immediately follow from Proposition \ref{ordliftring}. Fix a homomorphism $\mu \colon G \to S$ to an $\mc{O}$-torus $S$, and let $\mf{g}_{\mu}$ be the Lie algebra of $G_{\mu}= \ker \mu$. Fix a lift $\nu \colon \gal{F} \to S(\mc{O})$ of $\mu \circ \br$. The restriction of $\mu$ to $B$ factors through $B/N=T$, so $\mu \circ \br_T= \bar{\nu}$, and it makes sense to require that $\chi_T$ satisfy $\mu \circ \chi_T|_{I_F}= \nu|_{I_F}$. Then we can define the functor $\Lift_{\br}^{\nu, \chi_T}$ of lifts of $\br$ that are both of type $\chi_T$ and of type $\nu$. Since $\ell$ is very good for $G$, we conclude:
\begin{cor}
The functor $\Lift_{\br}^{\nu, \chi_T}$ defines a liftable deformation condition with tangent space $L_{\br}^{\nu, \chi_T}$ of dimension
\[
\dim_k L_{\br}^{\nu, \chi_T}= h^0(\gal{F}, \br(\mf{g}_{\mu}))+ [F:\Ql]\dim_k(\fn).
\]
\end{cor}
\begin{rmk}
This extra $[F:\Ql]\dim_k \mf{n}$ in the dimension of the local condition at $\ell$ is exactly what is needed in Ramakrishna's method (or originally in the Taylor-Wiles method) to offset the local archimedean invariants for `odd' representations. See Equation (\ref{wilesformula}) in the proof of Proposition \ref{formalliftingtheorem}.
\end{rmk}
\begin{rmk}\label{ordnotell}
There are analogous results in the case $p \neq \ell$, which follow from the same arguments. Namely, if we now take $F$ to be a finite extension of $\Q_p$, and let $\br \colon \gal{F} \to B(k)$ be a homomorphism satisfying (REG) and (REG*), then still assuming $\zeta_{\ell} \notin F$ (which, note, affects the calculation of $h^1(\gal{F}, \br(\fb/\fn))$) we deduce that the functor of type $\chi_T$ lifts gives rise to a liftable deformation condition whose associated tangent space has dimension $h^0(\gal{F}, \br(\fg))$, as is needed for Ramakrishna's global Galois cohomology argument. In \S \ref{steinberg}, we will consider one well-behaved example, deformations of Steinberg type, in which the condition (REG*) fails. In general, however, if (REG*) fails, or if $\zeta_{\ell}$ belongs to $F$, then the type-$\chi_T$ deformation ring can be singular, making it ill-suited for Ramakrishna's method.     
\end{rmk}
Finally, in our global applications we will want our characteristic zero lifts to be de Rham; this is not the case for all of the lifting functors $\Lift_{\br}^{\chi_T}$ considered in this section. We can ensure it as follows:
\begin{lemma}\label{ordderham}
Let $\chi_T \colon I_{F} \to T(\mc{O})$ be a lift of $\br_T$ such that for all $\alpha \in \Delta$, $\alpha \circ \chi_T= \kappa^{r_{\alpha}}$ for some positive integer $r_{\alpha}$. Let $(\rho^{\square, \chi_T}, R_{\br}^{\square, \chi_T})$ be the universal object for $\Lift_{\br}^{\chi_T}$. Then for all $\mc{O}$-algebra homomorphism $f \colon R_{\br}^{\square, \chi_T} \to \Qlb$, the push-forward $f(\rho^{\square, \chi_T})$ is de Rham.
\end{lemma}
\proof
First note that by the second part of Lemma \ref{ordcondlem}, the $\mc{O}$-points of $R_{\br}^{\square, \chi_T}$ admit a moduli description analogous to that of Definition \ref{ordcond}. While the de Rham condition is not stable under extensions, it is stable under direct sums and under extensions of the form
\begin{equation}\label{drext}
0 \to \kappa^r \to E \to \kappa^s \to 0 
\end{equation}
where $r > s$. For a lift $\rho \in \Lift_{\br}^{\chi_T}(\mc{O})$, which we may assume valued in $B(\mc{O})$, we apply this observation to the $\gal{F}$-stable filtration of $\mf{b}_{\Qlb}$ by root height, i.e. the decreasing filtration with $\Fil^i \mf{b}_{\Qlb}$ equal to the direct sum of all (positive) root spaces of height at least $i$. That $\mf{b}_{\Qlb}$ is de Rham follows by induction, since for any positive root $\gamma$, $\Ad(B)(\mf{g}_{\gamma})$ has non-zero $\mf{g}_{\gamma'}$-component only when $\gamma' -\gamma$ is a nonnegative linear combination of simple roots. The product map $B \to \mr{GL}(\mf{b}) \times T$ is faithful, so we can conclude $\rho$ itself is de Rham.
\endproof
\subsection{Ramakrishna deformations}\label{ramdefsection}
This section studies the local deformation condition used at the auxiliary primes of ramification in Ramakrishna's global argument. Let $F$ be a finite extension of $\Q_p$ for $p \neq \ell$, and let $\br \colon \gal{F} \to G(k)$ be an unramified homomorphism such that $\br(\fr_F)$ is a regular semi-simple element. 
Let $T$ be the connected component of the centralizer of $\br(\fr_F)$; this is a maximal $k$-torus of $G$, but we can lift it to an $\mc{O}$-torus (uniquely up to isomorphism), which we also denote $T$, and then we can lift the embedding over $k$ to an embedding $T \into G$ over $\mc{O}$ (see \cite[Corollary B.3.5]{conrad:luminy}); moreover, the latter lift is unique up to $\widehat{G}(\mc{O})$-conjugation. $T$ splits over an \'{e}tale extension $\mc{O}'/\mc{O}$ (corresponding to a finite extension of residue fields $k'/k$), 
and we for the rest of the section enlarge $\mc{O}$ to $\mc{O'}$ (and $k$ to $k'$). We therefore assume that $T$ is split, and we will invoke freely the resulting theory of roots and root subgroups (see \cite[\S 5.1]{conrad:luminy}).
\begin{defn}\label{deframtype}
An unramified residual representation $\br \colon \gal{F} \to G(k)$ is defined to be of Ramakrishna type if 
\begin{itemize}
\item $\br(\fr_F)$ is a regular semi-simple element; and
\item letting $T$ as above denote the connected component of the centralizer of $\br(\fr_F)$, there exists a root $\alpha \in \Phi(G, T)$ such that
\[
\alpha(\br(\fr_F))= \kbar(\fr_F)= q_F^{-1}.
\]
\end{itemize}
By the regularity assumption, the order $q_F$ of the residue field of $F$ is not congruent to  $1 \pmod \ell$. Let $H_{\alpha}= T \cdot U_{\alpha}$ be the subgroup of $G$ generated by $T$ and the root subgroup $U_{\alpha}$ corresponding to $\alpha$. We now define the lifts of $\br$ of Ramakrishna type: let $\Lift_{\br}^{\Ram}(R)$ be the subfunctor of $\Lift_{\br}$ consisting of all $\rho \in \Lift_{\br}(R)$ such that $\rho$ is $\widehat{G}(R)$-conjugate to a homomorphism $\gal{F} \xrightarrow{\rho'} H_{\alpha}(R)$, with the resulting composite 
\[
\gal{F} \xrightarrow{\rho'} H_{\alpha}(R) \xrightarrow{\Ad} \mr{GL}_R(\fg_{\alpha} \otimes R)=R^\times
\]
equal to $\kappa$.
\end{defn}
\begin{lemma}\label{ramcond}
For $\br$ of Ramakrishna type, $\Lift_{\br}^{\Ram}$ is well-defined and yields a liftable deformation condition.
\end{lemma}
\proof
We check the Mayer-Vietoris property of Lemma \ref{representability}. The argument is quite similar to that of Lemma \ref{ordcondlem}, and we give the details here since they were omitted above. Let $A \to C$ and $B \to C$ be morphisms in $\CO^f$, and assume $B \to C$ is small. Suppose we are given 
\[
\rho_A \times \rho_B \in \Lift_{\br}^{\Ram}(A) \times_{\Lift_{\br}^{\Ram}(C)} \Lift_{\br}^{\Ram}(B).
\]
By assumption, there exist $g_A \in \widehat{G}(A)$ and $g_B \in \widehat{G}(B)$ such that ${}^{g_A} \rho_A$ factors through $H_{\alpha}(A)$, and ${}^{g_B} \rho_B$ factors through $H_{\alpha}(B)$. We denote the push-forwards of $g_A$ and $g_B$ to $C$ by $g_{A, C}$ and $g_{B, C}$. Since the push-forwards of $\rho_A$ and $\rho_B$ to $G(C)$ are equal--denote this element of $\Lift_{\br}(C)$ by $\rho_C$--both ${}^{g_{B, C}} \rho_C$ and ${}^{(g_{A, C} g_{B, C}^{-1})g_{B, C}} \rho_C$ factor through $H_{\alpha}(C)$. We are thus led to consider, for any $\rho \colon \gal{F} \to H_{\alpha}(R)$ lifting $\br$, the set
\[
\mc{U}(\rho, R)= \{ g \in \widehat{G}(R): {}^g \rho(\gal{F}) \subset H_{\alpha}(R) \}.
\]
We claim $\mc{U}(\rho, R)= \widehat{H}_{\alpha}(R)$; from this it follows that the element $g_{A,C} g_{B, C}^{-1}$ can be lifted to $h \in \mc{U}({}^{g_B} \rho_B, B)= \widehat{H}_{\alpha}(B)$, so that $g_A$ and $h g_B$ have the same image in $C$; and this in turn easily implies that $\rho_A \times \rho_B$ is an element of $\Lift_{\br}^{\Ram}(A \times_C B)$.

To prove the claim, we argue by induction on the length of $R$, so we let $R \to R/I$ be a small morphism and assume the claim over $R/I$. Fix a $g \in \mc{U}(\rho, R)$, so the fiber of $\mc{U}(\rho, R) \to \mc{U}(\rho_{R/I}, R/I)$ containing $g$ consists of all elements of the form $\exp(Y) g$, for $Y \in \fg \otimes_k I$, such that ${}^{\exp(Y)g} \rho(\gal{F})\subset H_{\alpha}(R)$. But this implies that $\exp(Y- {}^{\br(\sigma)} Y) \cdot {}^g \rho(\sigma) \in H_{\alpha}(R)$ for all $\sigma \in \gal{F}$, hence that $\exp(Y- {}^{\br(\sigma)} Y) \in H_{\alpha}(k)$ for all $\sigma$. But now the regularity hypothesis implies that $Y$ belongs to $\mf{t} \oplus \fg_{\alpha} = \mr{Lie}(H_{\alpha})$, and the claim follows.

Having established that $\Lift_{\br}^{\Ram}$ is a deformation condition, we now check that it is liftable. Since $\widehat{G}$ is formally smooth, it suffices to show that we can lift an element $\rho_{R/I}$ of $\Lift_{\br}^{\Ram}(R/I)$ that factors through $H_{\alpha}(R/I)$ to $\Lift_{\br}^{\Ram}(R)$. Every element of $H_{\alpha}(R)$ can be written uniquely as a product of elements of $U_{\alpha}(R)$ and $T(R)$ (in this degree of generality, see \cite[Theorem 4.1.4]{conrad:luminy}), and writing
\[
\rho_{R/I}(g)= u_{\alpha}(x_g) \rtimes t_g \in U_{\alpha}(R/I) \rtimes T(R/I),
\]
we see that $g \mapsto t_g$ is a homomorphism, and $g \mapsto x_g$ is a cocycle in $Z^1(\gal{F}, R/I(1))$ (since $\alpha \circ t_g= \kappa(g)$, by assumption). Note also that $t_g$ is necessarily unramified, since $\br$ is unramified, $p \neq \ell$, and $q_F \not \equiv 1 \pmod \ell$. It is then easy to easy that $g \mapsto t_g$ lifts to a homomorphism (necessarily unramified) $\gal{F} \to T(R)$ whose composition with $\alpha$ is $\kappa$; and then to lift $g \mapsto x_g$ it suffices to see that $H^1(\gal{F}, R(1)) \to H^1(\gal{F}, R/I(1))$ is surjective. This follows from injectivity of $H^2(\gal{F}, I(1)) \to H^2(\gal{F}, R(1))$, which in turn (by local duality) follows from surjectivity of 
\[
\Hom(R, \Q/\Z) \to \Hom(I, \Q/\Z).
\]
\endproof
Next we describe the tangent space $L^{\Ram}_{\br}$; for later use, it is also convenient to describe the annihilator $L^{\Ram, \perp}_{\br}$ inside $H^1(\gal{F}, \br(\fg)(1))$. Consider the sub-torus $T_{\alpha}= \ker(\alpha)^0$ of $T$, and denote by $\mf{t}_{\alpha}$ its Lie algebra. There is a canonical decomposition $\mf{t}_{\alpha} \oplus \mf{l}_{\alpha} \xrightarrow{\sim} \mf{t}$ with $\mf{l}_{\alpha}$ the one-dimensional torus generated by the coroot $\alpha^\vee$.
\begin{lemma}\label{ramtangent}
Let $W= \mf{t}_{\alpha} \oplus \fg_{\alpha}$, and assume $\br$ is of Ramakrishna type. Then: 
\begin{enumerate}
\item The tangent space of $\Lift_{\br}^{\Ram}$ is (the preimage in $Z^1(\gal{F}, \br(\fg))$ of) 
\[
L_{\br}^{\Ram}= \im \left(H^1(\gal{F}, W) \to H^1(\gal{F}, \br(\fg)\right).
\]
\item $\dim L_{\br}^{\Ram}= h^0(\gal{F}, \br(\fg))$.
\item The orthogonal complement $L_{\br}^{\Ram, \perp} \subset H^1(\gal{F}, \br(\fg)(1))$ is equal to
\[
\im \left(H^1(\gal{F}, \br(W^\perp)(1)) \to H^1(\gal{F}, \br(\fg)(1)) \right),
\]
where $W^{\perp}$ denotes the annihilator of $W$ under the given $G$-invariant duality on $\fg$.
\item All cocycles in $L_{\br}^{\Ram, \square} \subset Z^1(\gal{F}, \br(\fg))$ have $\mf{l}_{\alpha}$-component, under the canonical decomposition $\fg= \bigoplus_{\gamma} \fg_{\gamma} \oplus \mf{t}_{\alpha} \oplus \mf{l}_{\alpha}$, equal to zero. All cocycles in $L_{\br}^{\Ram, \perp, \square} \subset Z^1(\gal{F}, \br(\fg)(1))$ have $\mf{g}_{-\alpha}$ component equal to zero.
\end{enumerate}
\end{lemma}
\proof
Let $\rho \in \Lift_{\br}^{\Ram}(k[\epsilon])$, with associated 1-cocycle $\phi \in Z^1(\gal{F}, \br(\fg))$. Then there exists $X \in \mf{g}$ such that
\[
\exp(\epsilon X) \exp(\epsilon \phi(g)) \br(g) \exp(-\epsilon X) \in H_{\alpha}(k[\epsilon])
\]
with image $\kbar(g)$ under $\alpha$. Modifying $\phi$ by $X- {}^{\br(g)} X$ we may then assume that $\phi \in Z^1(\gal{F}, \br(W))$, as desired.

To compute $\dim L_{\br}^{\Ram}$, consider the long exact sequence in $\gal{F}$-cohomology associated to the sequence of $\gal{F}$-modules
\[
0 \to \br(W) \to \br(\fg) \to \br(\fg/W) \to 0.
\]
Putting the definitions together with local duality, we find
\[
\dim L_{\br}^{\Ram}- h^0(\gal{F}, \br(\fg))= h^0(\gal{F}, \br(W^*)(1))- h^0(\gal{F}, \br(\fg/W))=0,
\]
since by assumption each term on the right-hand side is one-dimensional.

For part (3), note that $L_{\br}^{\Ram, \perp}$ clearly contains the image of $H^1(\gal{F}, \br(W^\perp)(1))$, so it suffices to check the two spaces have the same dimension. This follows as in part (2), now by considering the piece of the long exact sequence in cohomology beginning
\[
0 \to \im \left(H^1(\gal{F}, \br(W^\perp)(1))\right) \to H^1(\gal{F}, \br(\fg)(1)) \to \ldots.
\]
(Use the identifications $W^\perp \cong (\fg/W)^*$ and $\fg/W^{\perp} \cong W^*$ to reduce to the calculation in part (2).)

Part (4) follows immediately from parts (2) and (3) and an easy calculation: note that all root spaces except $\fg_{-\alpha}$ pair trivially with $W$.
\endproof
\begin{rmk}
If we fix the similitude character, the above goes through \textit{mutatis mutandis}.
\end{rmk}
\subsection{Steinberg deformations}\label{steinberg}
Let $F$ be a finite extension of $\Q_p$ for $p \neq \ell$. In this section we study deformations that generically correspond to Steinberg representations on the automorphic side. In the global application of \S \ref{-1notinW}, this deformation condition will be used at a prime $p \neq \ell$ to ensure the algebraic monodromy group of a characteristic zero lift is `big enough.'

For simplicity (this is all that will be needed in the application), assume $G$ is an adjoint group. We may then assume it is simple, with Coxeter number $h$. We assume $\kbar \colon \gal{F} \to \fl^\times$ has order greater than $h-1$. Let $B$ be a Borel subgroup of $G$ over $\mc{O}$, and let $N$ denote the unipotent radical of $B$, so $B/N= T$ is an $\mc{O}$-torus. 
We may as usual replace $\mc{O}$ by a finite \'{e}tale extension before defining our lifting functors, and so we may and do assume that $T$ is split, and that there is a section $B= N \rtimes T$ (see \cite[Proposition 5.2.3]{conrad:luminy}). The split torus gives us a root system $\Phi(G, T)$, and the Borel $B \supset T$ gives us a system of positive roots with simple roots $\Delta$.

\begin{defn}
Let $\br \colon \gal{F} \to B(k)$ be a representation factoring through $B$, and assume that for all $\alpha \in \Delta$ the composite
\[
\gal{F} \xrightarrow{\br} B(k) \to T(k) \xrightarrow{\alpha} k^\times
\]
is equal to $\kbar$. In this case we say $\br$ is of Steinberg type. We define the lifts of $\br$ of Steinberg type to be the sub-functor $\Lift^{\St}_{\br}(R)$ of $\rho \in \Lift_{\br}(R)$ such that $\rho$ is $\widehat{G}(R)$-conjugate to a homomorphism $\rho' \colon \gal{F} \to B(R)$, and for all $\alpha \in \Delta$ the composite 
\[
\gal{F} \xrightarrow{\rho'} B(R) \to T(R) \xrightarrow{\alpha} R^\times
\]
equals $\kappa$.
\end{defn}
\begin{lemma}
For $\br$ of Steinberg type, $\Lift^{\St}_{\br}$ is a well-defined deformation condition.
\end{lemma}
\proof
Our hypothesis that $1, \kbar, \kbar^2, \ldots, \kbar^{h-1}$ are all distinct implies that $H^0(\gal{F}, \br(\fg/\fb))=0$ (note that $h-1$ is the height of the highest root). Then the lemma follows exactly as in Lemma \ref{ordcondlem}.
\endproof
We now compute the tangent space of $\Lift_{\br}^{\St}$:
\begin{lemma}
The tangent space $L_{\br}^{\St}$ is equal to
\[
\ker \left( H^1(\gal{F}, \br(\fb)) \to H^1(\gal{F}, \br(\fb/\fn)) \right)
\]
and has dimension $h^0(\gal{F}, \br(\fg))$.
\end{lemma}
\proof
Note first that the assertion of the lemma makes sense, since $H^1(\gal{F}, \br(\fb)) \to H^1(\gal{F}, \br(\fg))$ is injective; so, too, is $H^1(\gal{F}, \br(\fb/\fn)) \to H^1(\gal{F}, \br(\fg/\fn))$. The description of $L_{\br}^{\St}$ follows easily, as in Lemma \ref{ordtanspace}. Note that $L_{\br}^{\St}$ also equals
\[
\ker \left( H^1(\gal{F}, \br(\fg)) \to H^1(\gal{F}, \br(\fg/\fn)) \right).
\]
The long exact sequence in $\gal{F}$-cohomology associated to $0 \to \br(\fn) \to \br(\fg) \to \br(\fg/\fn) \to 0$ then implies
\begin{align*}
\dim(L^{\St}_{\br})-h^0(\gal{F}, \br(\fg))= h^1(\gal{F}, \br(\fn))-h^0(\gal{F}, \br(\fn))-h^0(\gal{F}, \br(\fg/\fn)) \\
= h^2(\gal{F}, \br(\fn))-h^0(\gal{F}, \br(\fg/\fn))= h^0(\gal{F}, \br(\fg/\fb)(1))-h^0(\gal{F}, \br(\fg/\fn)). 
\end{align*}
In the last expression, each negative simple root contributes one dimension to the $\br(\fg/\fb)(1)$ term, and the image of $\mf{t}$ gives the $\br(\fg/\fn)$ term; the net contribution is then $\rk(G)-\rk(G)=0$, so $\dim(L_{\br}^{\St})= h^0(\gal{F}, \br(\fg))$.  
\endproof
Although the argument to this point closely resembles that of \S \ref{ordsection}, the failure in the Steinberg case of the condition $h^0(\gal{F}, \br(\fg/\fb)(1))=0$ means we have to work somewhat harder to establish liftability:
\begin{lemma}
$\Lift_{\br}^{\St}$ is liftable.
\end{lemma}
\proof
Let $R \to R/I$ be a small surjection, and let $\rho \in \Lift_{\br}^{\St}(R/I)$ be a lift of Steinberg type. We may assume that $\rho$ factors $\rho \colon \gal{F} \to B(R/I)$. Now, the Lie algebra $\fb$ admits a $B$-stable filtration by root height: let $F^0 \fb= \fb$, and for $r>0$ let 
\[
F^r \fb= \bigoplus_{\mr{ht}(\alpha) \geq r} \fg_{\alpha},
\]
so we have
\[
\fb= F^0 \fb \supset F^1 \fb \supset \cdots \supset F^{h-1} \fb \supset F^{h} \fb = 0.
\]
For each $i = 1, \ldots h$, let $N_{\geq i}$ be the closed subgroup of $N$, also a normal subgroup of $B$, whose Lie algebra is $F^i \fb$. We will construct a lift of $\rho$ by inductively constructing lifts to each $B/N_{\geq i}$. For $i=1$, the lift is forced on us by the definition of the Steinberg condition: we take the unique character $\gal{F} \to B/N_{\geq 1}(R)=T(R)$ whose composition with each simple root equals $\kappa$. The case of lifting to $B/N_{\geq 2}$ is rather special, and we postpone it. Assume then that $i  \geq 2$, and that by induction we are given the following commutative diagram, where our task is to fill in the dotted arrow:
\[
\xymatrix{
 && B/N_{\geq i+1}(R) \ar[d] \ar[r] & B/N_{\geq i}(R) \ar[d] \\
\gal{F} \ar@{-->}[urr] \ar[urrr]^-{\tilde{\rho}_{\geq i}} \ar[rr]_-{\rho_{\geq i+1}} & & B/N_{\geq i+1}(R/I) \ar[r] & B/N_{\geq i}(R/I).
}
\]
The obstruction to lifting $\rho_{\geq i+1}$ is, by standard obstruction theory and existence of the lift $\tilde{\rho}_{\geq i}$, an element of
\[
\ker \left(H^2(\gal{F}, \br(\fb/ F^{i+1} \fb))\otimes_k I \to H^2(\gal{F}, \br(\fb/ F^i \fb)) \otimes_k I \right).
\]
But observe that
\[
h^2(\gal{F}, \br(F^i \fb/ F^{i+1} \fb)) = h^2(\gal{F}, \kbar^i)^{\oplus \dim F^i \fb/F^{i+1} \fb}= h^0(\gal{F}, \kbar^{1-i})^{\oplus \dim F^i \fb/ F^{i+1} \fb}=0,
\]
since we are considering the cases $i= 2, \ldots, h-1$, and we have assumed the order of $\kbar$ is greater than $h-1$.

Thus it remains only to consider the case when $i=1$ in the above diagram. Via our fixed splitting $B= N \rtimes T$, we can write
\[
\rho_{\geq 2} \colon \gal{F} \to N/N_{\geq 2}(R/I) \rtimes T(R/I)
\]
in the form 
\[
\rho_{\geq 2}(g)= \exp(\phi(g)) \rtimes \rho_{\geq 1}(g)
\] 
for some function $\phi \colon \gal{F} \to F^1 \fb/F^2 \fb$. We can clearly lift $\rho_{\geq 1}$ to a homomorphism $\tilde{\rho}_{\geq 1} \colon \gal{F} \to T(R)$ satisfying $\alpha \circ \tilde{\rho}_{\geq 1}= \kappa$ for all $\alpha \in \Delta$, so having done this we need only address lifting $\phi$. But note that $\phi$ is simply an element of
\[
Z^1(\gal{F}, \rho_{\geq 1}(F^1 \fb/ F^2 \fb))= Z^1(\gal{F}, \oplus_{\alpha \in \Delta} R/I(1)),
\]
so to construct a lift $\tilde{\rho}_{\geq 2}$ (with the fixed push-forward $\tilde{\rho}_{\geq 1}$ to $T$), we need only lift this co-cycle to an element of $Z^1(\gal{F}, \oplus_{\alpha \in \Delta} R(1))$; this is done exactly as at the end of Lemma \ref{ramcond}, and the proof is complete.
\endproof
\subsection{Minimal prime to $\ell$ deformations}\label{minimal}
In the application, we will require one more especially simple local condition. Continue to assume $F$ is a finite extension of $\Q_p$ with $p \neq \ell$. Now suppose that $\br \colon \gal{F} \to G(k)$ satisfies $\ell \nmid \br(I_F)$. Let $\Lift_{\br}^{\mc{P}}$ be the deformation condition consisting of all $\rho \in \Lift_{\br}(R)$ such that $\rho|_{I_F}$ factors through the fixed field of $\br|_{I_F}$.
\begin{lemma}
Under the above hypotheses, $\Lift_{\br}^{\mc{P}}$ is a liftable deformation condition whose tangent space has dimension $h^0(\gal{F}, \br(\fg))$.
\end{lemma}
\proof
This follows immediately from the Hochschild-Serre spectral sequence and the following standard facts: 
\begin{itemize}
\item $H^i(\br(I_F), \br(\fg))=0$
for all $i>0$;
\item $H^2(\gal{F}/I_F, M)=0$ for all finite $\gal{F}/I_F$-modules $M$ (take $M= \br(\fg)^{I_F}$); and
\item $h^1(\gal{F}/I_F, M^{I_F})=h^0(\gal{F}, M)$ for all finite $\gal{F}$-modules $M$.
\end{itemize}
\endproof
\subsection{The archimedean condition}\label{archcond}
Our global deformation problems will not explicitly impose any condition at the archimedean places, but the archimedean deformations will implicitly be dictated by properties of the residual representation. Basic to Ramakrishna's method, and to the original form of the Taylor-Wiles method, is the requirement that the residual representation be suitably `odd.' We will now explain this oddness condition. The reader might wish to glance ahead to Equation (\ref{wilesformula}) in \S \ref{formalglobal} (with reference to assumptions (\ref{dualinvs}) and (\ref{localconditions}), also in \S \ref{formalglobal}). For an appropriately chosen collection of local deformation conditions (we restrict for simplicity here to semi-simple $\fg$, so there is no `multiplier character'), this equation gives an equality\footnote{The logic of the present section only depends on knowing the right-hand-side of this equation; we mention the rest only for motivation.}
\begin{equation}\label{balance}
\dim H^1_{\mc{P}}(\gal{F, \Sigma}, \br(\fg))- \dim H^1_{\mc{P}^{\perp}}(\gal{F, \Sigma}, \br(\fg)(1))= [F:\Q] \dim \mf{n}- \sum_{v \vert \infty}h^0(\gal{F_v}, \br(\fg)).
\end{equation}
We recall that an involution $\tau$ of $G$ is called a split Cartan involution if $\dim \fg^{\tau}= \dim (\fn)$.
\begin{lemma}
In order for the right-hand-side of Equation (\ref{balance}) to be nonnegative, $F$ must be totally real, and $\Ad(\br(c_v))$ must be a split Cartan involution for all complex conjugations $c_v$. 
\end{lemma}
\proof
Observe that
\begin{align*}
&\sum_{v \vert \infty} h^0(\gal{F_v}, \br(\fg))= \sum_{\text{$v$ complex}} \dim \fg + \sum_{\text{$v$ real}} \dim \fg^{c_v} \geq \\ 
& \sum_{\text{$v$ complex}} \dim \fg + \sum_{\text{$v$ real}} \dim \mf{n}= [F:\Q] \dim \mf{n}+ \#\{\text{$v$ complex}\}\cdot \dim \mf{t},
\end{align*}
with equality only when every $\Ad(\br(c_v))$ is a split Cartan involution of $\mf{g}$ (see \cite[Proposition 2.2]{yun:exceptional} for this result, due to Cartan, about involutions of a reductive group). The value of Equation (\ref{balance}) being nonnegative thus forces $F$ to be totally real and all $\br(c_v)$ to induce split Cartan involutions of $\fg$.
\endproof
We then ask: what (connected reductive) groups $G$ contain an order two element $c$ such that $\Ad(c)$ is a split Cartan involution?
\begin{lemma}\label{splitcartan}
If $-1$ belongs to the Weyl group of $G$, and the co-character $\rho^\vee$ (the half-sum of the positive co-roots) of $G^{\mr{ad}}$ lifts to a co-character of $G$, then $G(k)$ contains an element $c$ of order 2 such that $\Ad(c)$ is a split Cartan involution of $\mf{g}^{\mr{der}}$. If $-1$ does not belong to the Weyl group of $G$, then $G(k)$ contains no such element $c$.
\end{lemma}
\proof
Fix a split maximal torus $T$ of $G$ and a choice of positive system of roots with respect to $T$. If $-1 \in W_G$, then $\rho^\vee(-1) \in G^{\mr{ad}}(k)$ is a split Cartan involution by \cite[Lemma 2.3]{yun:exceptional}. If $\rho^\vee$ lifts to $G$, then of course the same is true of $\rho^\vee(-1) \in G(k)$. If we do not assume $-1 \in W_G$, but just consider any split Cartan involution $\tau$, then decomposing $\fg= \fg^+ \oplus \fg^-$ into $\pm 1$-eigenspaces for $\tau$, one checks that the maximal abelian semi-simple subalgebra $\mf{s}$ of $\mf{g}^-$ must be a Cartan sub-algebra (see the proof of \cite[Proposition 2.2]{yun:exceptional}). Thus $\tau$ acts as $-1$ on $\mf{s}$, and if $\tau$ were an inner automorphism of $G$, we would necessarily have $-1 \in W_G$.
\endproof
\begin{rmk}
\begin{itemize}
\item Note also that while $\rho^\vee$ does not always lift to $G$ (eg, $G= \mr{SL}_2$), we can always enlarge the center of $G$ to make the lift possible (eg, $G= \mr{GL}_2$). See Equation (\ref{enlargeG}) in \S \ref{principalsl2.1}.
\item Our later arguments do not logically require this description of the case $-1 \not \in W_G$, but we provide it for motivation. When $-1$ does not belong to $W_G$, we will have to work with a suitable non-connected extension of $G$, in order for the split Cartan involution to be an inner (in the larger group) automorphism: see \S \ref{outerodd}.
\end{itemize}
\end{rmk}
\section{The global theory: axiomatizing Ramakrishna's method for annihilating the dual Selmer group}\label{formalglobal}
Let $\br \colon \gal{F, \Sigma} \to G(k)$ be a continuous homomorphism such that the infinitesimal centralizer of $\br$ is $\mf{z}(\fg)$, so that the deformation functor is representable.
We now begin to explain the global Galois cohomological argument, due to Ramakrishna \cite{ramakrishna02} for $G= \mr{GL}_2$ and $F= \Q$, that under favorable circumstances allows us to find a geometric characteristic zero lift of $\br$. To be precise, in this section we will explain an `axiomatized' version of the method; then in \S \ref{bigimage} and \S \ref{principalSL2} we explain precise conditions on $\br$ that allow this axiomatized method to run successfully. Recalling that $S$ is the maximal torus quotient of $G$, We fix once and for all a de Rham `similitude character' $\nu$ lifting $\mu \circ \br$: 
\[
\xymatrix{
& S(\mc{O}) \ar[d] \\
\gal{F, \Sigma} \ar[r]^{\mu \circ \br} \ar[ur]^{\nu} & S(k).
}
\]
We will henceforth only consider, both locally and globally, lifts of $\br$ with fixed similitude character $\nu$ (see Example \ref{fixedchar}); thus $\fg_{\mu}= \fg^{\mr{der}}$ is the Galois module appearing as coefficients in all of our cohomology groups measuring deformations of $\br$. We remark, however, that the reader will lose little simply by assuming $G$ is adjoint.

Let $K= F(\br(\fg_{\mu}), \mu_{\ell})$. For the rest of this section we assume the following:
\begin{enumerate}
\item\label{dualinvs} $h^0(\gal{F}, \brgm)= h^0(\gal{F}, \brgm(1))=0$. We note here for later use that since $\ell$ is very good for $G$, our centralizer condition on $\br$ is equivalent to the condition $h^0(\gal{F}, \brgm)=0$.
\item\label{localconditions} There is a global deformation condition $\mc{P}= \{\mc{P}_v\}_{v \in \Sigma}$ consisting of liftable local deformation conditions for each place $v \in \Sigma$ (taking fixed multiplier character, both locally and globally); the dimensions of their tangent spaces are
\[
\dim L_v=
\begin{cases}
\text{$0$ if $v \vert \infty$;}\\
\text{$h^0(\gal{F_v}, \br(\fg_{\mu}))$ if $v \nmid \ell \cdot \infty$;}\\
\text{$h^0(\gal{F_v}, \br(\fg_{\mu})) + [F_v :\Q_{\ell}] \dim (\mf{n})$ if $v \vert \ell$.}
\end{cases}
\]
\item\label{infiniteplaces} $F$ is totally real, and for all $v \vert \infty$,
\[
h^0(\gal{F_v}, \br(\fg_{\mu})= \dim (\mf{n})
\]
\item\label{cohvanish} $H^1(\Gal(K/F), \br(\fg_{\mu}))=0$ and $H^1(\Gal(K/F), \brgm(1))=0$.
\item\label{lindisjoint} Assume item (\ref{cohvanish}) holds. For any pair of non-zero Selmer classes $\phi \in H^1_{\mc{P}^\perp}(\gal{F, \Sigma}, \br(\fg_{\mu})(1))$ and $\psi \in H^1_{\mc{P}}(\gal{F, \Sigma}, \br(\fg_{\mu}))$, we can of course restrict $\phi$ and $\psi$ to $\gal{K}$, where they become homomorphisms (rather than twisted homomorphisms), which are non-zero by item (\ref{cohvanish}). Letting $K_{\phi}/K$ and $K_{\psi}/K$ be their respective fixed fields, we assume that $K_{\phi}$ and $K_{\psi}$ are linearly disjoint over $K$.
\item\label{avoid} Consider any $\phi$ and $\psi$ as in the hypothesis of item (\ref{lindisjoint}) (we do not require the conclusion to hold). Then there is an element $\sigma \in \gal{F}$ such that $\br(\sigma)$ is a regular semi-simple element of $G$, the connected component of whose centralizer we denote $T$, and such that there exists a root $\alpha \in \Phi(G, T)$ satisfying
\begin{enumerate}
\item $\kbar(\sigma)= \alpha \circ \br(\sigma)$; 
\item $k[\psi(\gal{K})]$ has an element with non-zero $\mf{l}_{\alpha}$ component;\footnote{Recall from Lemma \ref{ramtangent} that $\mf{l}_{\alpha}$ is the span of the $\alpha$-coroot vector.} and
\item $k[\phi(\gal{K})]$ has an element with non-zero $\mf{g}_{-\alpha}$ component.
\end{enumerate}
\end{enumerate}
\begin{rmk}
We note that these conditions continue to hold if we replace $k$ by a finite extension, hence if we replace $\mc{O}$ by the ring of integers in any finite extension of $\Frac(\mc{O})$. We will use this flexibility freely in the applications of Proposition \ref{formalliftingtheorem}.
\end{rmk}
\begin{prop}\label{formalliftingtheorem}
Under assumptions (\ref{dualinvs})-(\ref{avoid}) above, there exists a finite set of primes $Q$ disjoint from $\Sigma$, and a lift 
\[
\xymatrix{
& G(\mc{O}) \ar[d] \\
\gal{F, \Sigma \cup Q} \ar[r]_{\br} \ar@{-->}[ur]^{\rho} & G(k) 
}
\]
such that $\rho$ is type $\mc{P}_v$ at all $v \in \Sigma$ and of Ramakrishna type at all $v \in Q$. 
\end{prop}
\proof
The strategy is to allow additional ramification at auxiliary primes (those in $Q$) to define a global deformation problem whose corresponding dual Selmer group vanishes; then we will conclude from Corollary \ref{globdefring}. We are already done unless there is a nonzero element $\phi \in H^1_{\mc{P}^\perp}(\gal{F, \Sigma}, \br(\fg_{\mu})(1))$. Recall Wiles's formula:
\begin{align}
&\dim H^1_{\mc{P}}(\gal{F, \Sigma}, \br(\fg_{\mu}))- \dim H^1_{\mc{P}^{\perp}}(\gal{F, \Sigma}, \br(\fg_{\mu})(1))= \label{wilesformula}\\
&h^0(\gal{F}, \br(\fg_{\mu}))-h^0(\gal{F}, \br(\fg_{\mu})(1))+ \sum_v (\dim L_v-h^0(\gal{F_v}, \br(\fg_{\mu}))= \nonumber \\
&\sum_{v \vert \ell} [F_v: \Q_{\ell}] \dim \mf{n}- \sum_{v \vert \infty}h^0(\gal{F_v}, \br(\fg_{\mu}))=0, \nonumber
\end{align}
the equalities of the final line following from assumptions (\ref{dualinvs}), (\ref{localconditions}), and (\ref{infiniteplaces}). In particular, having such a non-zero $\phi$ forces the existence of a non-zero $\psi \in H^1_{\mc{P}}(\gal{F}, \br(\fg_{\mu}))$.

We will see that the hypotheses of the proposition allow us to achieve the following: 
\begin{lemma}\label{findingw}
There exist infinitely many primes $w \notin \Sigma$ such that $\br|_{\gal{F_w}}$ is of Ramakrishna type and
\begin{align}
&\phi|_{\gal{F_w}} \notin L_w^{\Ram, \perp}; \label{phihyp}\\
&\psi|_{\gal{F_w}} \notin L_w^{\mr{unr}} \cap L_w^{\Ram}; \label{psihyp}
\end{align}
here $L_w^{\mr{unr}}$ denotes the tangent space of the unramified local condition (which is what is implicitly taken we study deformations of type $\mc{P}$), and $L_w^{\Ram}$ denotes (as in Lemma \ref{ramtangent}) the tangent space of the local condition of Ramakrishna type. 
\end{lemma}
We admit the existence of such a $w$ for the time being, and show how to conclude the argument of Proposition \ref{formalliftingtheorem}. Let $L_w= L_w^{\mr{unr}} \cap L_w^{\Ram}$, so that $L_w^{\perp}= L_w^{\unr, \perp}+ L_w^{\Ram, \perp}$. There are evident inclusions
\begin{align}
H^1_{\mc{P}^\perp \cup L_w^{\Ram, \perp}}(\gal{F, \Sigma \cup w}, \br(\fg_{\mu})(1)) \to &H^1_{\mc{P}^\perp \cup L_w^{\perp}}(\gal{F, \Sigma \cup w}, \br(\fg_{\mu})(1)), \label{raminf}\\
H^1_{\mc{P}^\perp}(\gal{F, \Sigma}, \br(\fg_{\mu})(1)) \to &H^1_{\mc{P}^\perp \cup L_w^{\perp}}(\gal{F, \Sigma \cup w}, \br(\fg_{\mu})(1)), \label{unrinf}
\end{align}
and we claim the second of these is an isomorphism. A double invocation of Wiles's formula (alleviating the notation with the self-explanatory shorthand) gives
\[
h^1_{\mc{P}^\perp \cup L_w^{\perp}}- h^1_{\mc{P}^\perp}= h^1_{\mc{P} \cup L_w}-h^1_{\mc{P}}-\dim L_w+ h^0(\gal{F_w}, \br(\fg_{\mu})),
\]
and the right-hand side of this equality is zero: indeed, this follows by exactness of the sequence
\[
0 \to H^1_{\mc{P} \cup L_w}(\gal{F, \Sigma \cup w}, \brgm) \to H^1_{\mc{P}}(\gal{F, \Sigma}, \brgm) \to L_w^{\mr{unr}}/L_w \to 0,
\]
where for surjectivity we use the assumption that $\psi|_{\gal{F_w}} \notin L_w^{\unr} \cap L_w^{\Ram}$, and the fact that
\[
\dim \left( L_w^{\unr}/(L_w^{\unr} \cap L_w^{\Ram}) \right)= 1.
\]
Note that the intersection here is (see Lemma \ref{ramtangent}) 
\[
\im \left(H^1(\gal{F_w}, \br(\mf{t}_{\alpha})) \to H^1(\gal{F_w}, \br(\fg)) \right).
\]
So, the map of Equation (\ref{unrinf}) is an isomorphism, and combined with Equation (\ref{raminf}) we get an exact sequence
\[
0 \to H^1_{\mc{P}^\perp \cup L_w^{\Ram, \perp}}(\gal{F, \Sigma \cup w}, \br(\fg_{\mu})(1)) \to H^1_{\mc{P}^\perp}(\gal{F, \Sigma}, \br(\fg_{\mu})(1)) \to H^1(\gal{F_w}, \br(\fg_{\mu})(1))/L_w^{\Ram, \perp}.
\]
By assumption (Equation (\ref{phihyp})), $\phi \in H^1_{\mc{P}^\perp}(\gal{F, \Sigma}, \br(\fg_{\mu})(1))$ does not restrict to an element of $L_w^{\Ram, \perp}$, so $H^1_{\mc{P}^\perp \cup L_w^{\Ram, \perp}}(\gal{F, \Sigma \cup w}, \br(\fg_{\mu})(1))$ has strictly smaller order than $H^1_{\mc{P}^\perp}(\gal{F, \Sigma}, \br(\fg_{\mu})(1))$, and by induction the proof of Proposition \ref{formalliftingtheorem}, modulo Lemma \ref{findingw}, is complete.
\endproof
It remains to prove Lemma \ref{findingw}:
\proof[Proof of Lemma \ref{findingw}]
For our fixed $\phi$ and $\psi$, we obtain a $\sigma \in \gal{F}$, a maximal torus $T$, and an $\alpha \in \Phi(G, T)$ as in assumption (\ref{avoid}). Choose any lift $\tilde{\sigma}$ of $\sigma$ to $\gal{F}$. By assumptions (\ref{lindisjoint}) and (\ref{avoid}), we can find a $\tau \in \Gal(K_{\phi} K_{\psi}/K)$ such that for any lift $\tilde{\tau}$ of $\tau$ to $\gal{K}$,
\begin{align}
&\text{$\phi(\tilde{\tau} \tilde{\sigma})$ has non-zero $\mf{g}_{-\alpha}$-component; and} \label{phidesid}\\
&\text{$\psi(\tilde{\tau} \tilde{\sigma})$ has non-zero $\mf{l}_{\alpha}$-component.} \label{psidesid}
\end{align}
To be precise, $\br(\tilde{\tau})$ acts trivially on $\br(\fg_{\mu})(1)$, so 
\[
\phi(\tilde{\tau} \tilde{\sigma})= \phi(\tilde{\tau})+ {}^{\br(\tilde{\tau})} \phi(\tilde{\sigma})= \phi(\tilde{\tau})+ \phi(\tilde{\sigma}),
\]
and likewise $\psi(\tilde{\tau} \tilde{\sigma})= \psi(\tilde{\tau})+\psi(\tilde{\sigma})$. Whatever $\phi(\tilde{\sigma})$ and $\psi(\tilde{\sigma})$ may be, we can, by hypotheses (\ref{lindisjoint}) and (\ref{avoid}), then find $\tau \in \Gal(K_{\phi}K_{\psi}/K)$ satisfying the conditions in equations (\ref{phidesid}) and (\ref{psidesid}).

Finally, by the \v{C}ebotarev density theorem, applied to the Galois extension $F(\br)K_{\phi}K_{\psi}/F$,\footnote{$K_{\phi}$ and $K_{\psi}$ are Galois over $F$ because $\phi$ and $\psi$ are cocycles for $\gal{F}$.} we can find a positive density set of primes $w$ of $F$ at which $\br$ is unramified and such that $\fr_w= \tilde{\tau} \tilde{\sigma}$ in $\Gal(F(\br)K_{\phi}K_{\psi}/F)$. Note that by construction $\br(\tilde{\tau})$ belongs to $Z_G(k)$ (it acts trivially in the adjoint representation) and $\kbar(\tilde{\tau})=1$. Thus $\tilde{\tau} \tilde{\sigma}$ satisfies the same hypothesis (\ref{avoid}) that $\sigma$ was assumed to satisfy, and $\br|_{\gal{F_w}}$ is therefore of Ramakrishna type. The lemma now follows from the explicit description (Lemma \ref{ramtangent}) of cocycles in $L_w^{\Ram, \square}$ and $L_w^{\Ram, \perp, \square}$.
\endproof

\section{The case $\br(\gal{F}) \supset G_{\mu}(\mathbb{F}_{\ell})$}\label{bigimage}
In this section, under the most generous assumptions on the image of $\br$, we show that all the hypotheses of Proposition \ref{formalliftingtheorem} can be satisfied (except possibly the existence of suitable local deformation conditions at places in $\Sigma$). For simplicity, we will assume the derived group of $G$ is almost-simple. There is no essential difficulty in passing to the general semi-simple case, but this way we save ourselves a little bookkeeping. As always, let $h$ denote the Coxeter number of $G$. For any extension $k'$ of $\mathbb{F}_{\ell}$ we introduce the notation
\[
\overline{G}^{\mr{sc}}_{\mu}(k')= \im \left( G^{\mr{sc}}_{\mu}(k') \to G_{\mu}(k') \right),
\]
where $G^{\mr{sc}}_{\mu} \to G_{\mu}$ dnotes the simply connected cover of $G_{\mu}$ (i.e., the simply-connected Chevalley group of the appropriate type). We then assume throughout this subsection that for some sub-extension $k \supset k' \supset \mathbb{F}_{\ell}$, 
\begin{equation}\label{bigimageassumption}
\overline{G}^{\mr{sc}}_\mu(k') \subseteq \br(\gal{F}) \subseteq Z_G(k) \cdot G(k').
\end{equation}
Note that $G_{\mu}(k')$ is then normal in $\br(\gal{F})$.
\begin{eg}
The template for the assumption (\ref{bigimageassumption}) is the following theorem of Ribet (building on ideas of Serre and Swinnerton-Dyer) about the images of Galois representations associated to holomorphic modular forms:
\begin{thm}[Theorem 3.1 of \cite{ribet:modlimage2}]\label{ribetimage}
If $f$ is a new eigenform in $S_k(\Gamma_1(N))$, with field of coefficients $E_f$ (with ring of integers $\mc{O}_f$), then for almost all $\lambda$, there is a subfield $k'_{\lambda} \subset \mc{O}_{f, \lambda}/\lambda= k_{\lambda}$ such that the associated mod $\lambda$ Galois representation $\bar{r}_{f, \lambda} \colon \gal{\Q} \to \mr{GL}_2(k_\lambda)$ has, after suitable conjugation, image containing $\mr{SL}_2(k'_{\lambda})$ as a normal subgroup.
\end{thm}
An elementary (Hilbert Theorem 90) argument shows that the normalizer of $\mr{SL}_2(k'_{\lambda})$ in $\mr{GL}_2(k_{\lambda})$ is $k_{\lambda}^\times \cdot \mr{GL}_2(k'_{\lambda})$, so the condition (\ref{bigimageassumption}) is indeed a natural one to impose.
\end{eg}

More generally, if $\br(\gal{F})$ contains $G_{\mu}(k')$ as a normal subgroup, then the obstruction to having $\br(\gal{F})$ contained in $Z_G(k) \cdot G(k')$ is an element of $H^1(\Gal(k/k'), Z_G(k))$. This group vanishes if $Z_G$ is a split torus, but not in general. We will not require $Z_G$ to be a torus, but the reader should keep in mind that the plausibility of hypothesis (\ref{bigimageassumption}) may depend on such an assumption. Note that every $G$ can be enlarged to a group $\tG$ with $Z_{\tG}$ a split torus and $G_{\mu}= \tG^{\mr{der}}$. 

For convenience, we recall the assumptions on $\ell$:
\begin{assumption}\label{whatsell}
As always (see \S \ref{review}), we assume $\ell>2$ is a `very good prime' for $G$. Additionally, we require $\ell$ to be greater than the maximum value of $\langle \alpha^\vee, \theta \rangle$, where $\theta$ denotes the highest root of $G$ and $\alpha^\vee$ ranges over simple coroots of $G$; but this condition is satisfied for all $\ell \geq 3$ (checking each simple type), so it in fact is no further constraint. We deduce that the adjoint representation $\mf{g}_{\mu}$ is not only irreducible as an algebraic representation (as ensured by $\ell$ being very good), but also irreducible as $k[G_{\mu}(\mathbb{F}_{\ell})]$-module (\cite[Theorem 43, pg 230]{steinberg:chevalley}). Note also the following:
\begin{itemize}
\item The index of $\overline{G}^{\mr{sc}}_{\mu}(k')$
in $G_{\mu}(k)$ is prime to $\ell$.
\item $G_{\mu}(k')$ is equal to its commutator subgroup whenever $\ell >3$ or whenever $\ell=3$ and $G$ is not of type $\mr{A}_1$, by \cite[Lemma $32^\prime$]{steinberg:chevalley}. We omit the (easy) modifications needed for our arguments in type $\mr{A}_1$, since this case has of course already been treated by Ramakrishna.
\end{itemize}
\end{assumption}
First note that condition (\ref{dualinvs}) of \S \ref{formalglobal}, that $h^0(\gal{F}, \brgm)= h^0(\gal{F}, \brgm(1))=0$, is satisfied, since $\brgm$ and $\brgm(1)$ are irreducible $\gal{F}$-representations under the assumptions on the image of $\br$ and on $\ell$. In particular, $\Def_{\br}$ is representable. We now state the main result of this section, which, note, incorporates strictly stronger assumptions on $\ell$ than those of Assumption \ref{whatsell}.
\begin{thm}\label{maximaltheorem}
Let $F$ be a totally real field with $[F(\mu_{\ell}):F]= \ell -1$, and let $\br \colon \Gamma_{F, \Sigma} \to G(k)$ be a continuous representation satisfying the following conditions:
\begin{enumerate}
\item\label{maximage} There is a subfield $k' \subset k$ such that
\[
\overline{G}^{\mr{sc}}_{\mu}(k') \subset \br(\gal{F}) \subset Z_G(k) \cdot G(k').
\]
\item\label{powers} $\ell -1$ is greater than the maximum of $8\cdot \# Z_{G^{\mr{sc}}_{\mu}}$ and
\[
\begin{cases}
\text{$(h-1)\# Z_{G^{\mr{sc}}_{\mu}}$ if $\# Z_{G^{\mr{sc}}_{\mu}}$ is even; or} \\
\text{$(2h-2) \# Z_{G^{\mr{sc}}_{\mu}}$ if $\# Z_{G^{\mr{sc}}_{\mu}}$ is odd.}
\end{cases}
\]
\item $\br$ is odd, i.e. for all complex conjugations $c_v$, $\Ad(\br(c_v))$ is a split Cartan involution of $G$.
\item For all places $v \in \Sigma$ not dividing $\ell \cdot \infty$, $\br|_{\gal{F_v}}$ satisfies a liftable local deformation condition $\mc{P}_v$ with tangent space of dimension $h^0(\gal{F_v}, \br(\fg_{\mu}))$ (eg, the conditions of \S \ref{steinberg} or \S \ref{minimal}).
\item For all places $v \vert \ell$, $\br|_{\gal{F_v}}$ is ordinary in the sense of \S \ref{ordsection}, satisfying the conditions (REG) and (REG*).
\end{enumerate}
Then there exists a finite set of primes $Q$ disjoint from $\Sigma$, and a lift 
\[
\xymatrix{
& G(\mc{O}) \ar[d] \\
\gal{F, \Sigma \cup Q} \ar[r]_{\br} \ar@{-->}[ur]^{\rho} & G(k) 
}
\]
such that $\rho$ is type $\mc{P}_v$ at all $v \in \Sigma$ (taking $\mc{P}_v$ to be an appropriate ordinary condition at $v \vert \ell$) and of Ramakrishna type at all $v \in Q$. In particular $\br$ admits a characteristic zero lift that is geometric in the sense of Fontaine-Mazur.
\end{thm}
\begin{rmk}
\begin{itemize}
\item Although not strictly necessary for the statement of the theorem, it is essential to note that the oddness hypothesis on $\br$ will never be met unless $-1$ belongs to the Weyl group of $G$.
\item We have made no attempt to optimize the hypotheses on $\ell$; sharper results can be extracted by examining the proof below, and still sharper results can be obtained by minor variations on the argument and case-by-case analysis. For example, taking $G= \mr{GSp}_{2n}$, a version of this theorem was established in \cite{stp:undergrad} with somewhat tighter bounds on $\ell$.
\end{itemize}
\end{rmk}
\proof
The proof of this theorem will occupy the rest of this section. Conditions (\ref{dualinvs})-(\ref{infiniteplaces}) of \S \ref{formalglobal} follow immediately from the hypotheses of the theorem.
 
We proceed to the condition (\ref{cohvanish}) of \S \ref{formalglobal}, namely that
\[
H^1(\Gal(K/F), \br(\fg_{\mu}))=0.
\]
Recall that $K= F(\br(\fg_{\mu}), \mu_{\ell})$.
\begin{lemma}\label{cps}
Under assumptions (\ref{maximage}) and (\ref{powers}) of Theorem \ref{maximaltheorem}, $H^1(\Gal(K/F), \br(\fg_{\mu}))=0$.

If we further assume that $F(\zeta_{\ell})$ is not contained in $F(\br(\fg_{\mu}))$, then $H^1(\Gal(K/F), \br(\fg_{\mu})(1))=0$ as well.
\end{lemma}
\proof
Repeated inflation-restriction arguments, using that $[Z_G(k) G(k'): G_{\mu}(k')]$ and $[F(\zeta_{\ell}):F]$ are coprime to $\ell$, reduce the desired vanishing to the assertion that $H^1(G_{\mu}(k'), \mf{g}_{\mu}(k'))=0$. Irreducibility of $\mf{g}_{\mu}$ lets us apply \cite[Corollary 2.9]{cline-parshall-scott}, and by \cite[Proposition 3.3]{cline-parshall-scott} our assumptions on $\ell$ (in particular, $\ell >9$) imply there are no non-trivial `Galois equivalences' between roots; then the output of \cite[Corollary 2.9]{cline-parshall-scott} is precisely that $H^1(G^{\mr{sc}}_{\mu}(k'), \mf{g}_{\mu}(k'))=0$. Another inflation-restriction argument (recalling Assumption \ref{whatsell}) implies that $H^1(G_{\mu}(k'), \mf{g}_{\mu}(k'))=0$ as well.

For the second point, let $H= \gal{F(\br(\fg_{\mu}))}$ and $H'=\gal{F(\br(\fg_{\mu})(1))}$. Note that $F(\br(\fg_{\mu})(1)) \supseteq F(\br(\fg_{\mu}), \mu_{\ell})$, and we claim that equality in fact holds. $H'$ acts on $\br(\fg_{\mu})$ as scalar multiplication by the character $\kbar^{-1}$, but no element of $G$ can act by a non-trivial scalar in the adjoint representation. Thus, $H'$ acts trivially both on $\br(\mf{g}_{\mu})$ and on $\mu_{\ell}$, so $F(\br(\fg_{\mu})(1))= F(\br(\fg_{\mu}), \mu_{\ell})$.
To conclude, since $H H'/H' \cong H/(H \cap H')$ has order prime to $\ell$ we get an (inflation) isomorphism
\[
H^1(\gal{F}/HH', \br(\fg_{\mu})(1)^H) \xrightarrow{\sim} H^1(\gal{F}/H', \br(\fg_{\mu})(1))= H^1(\Gal(K/F), \brgm(1)).
\]
These groups are clearly zero once $\kbar|_H \neq 1$ (i.e. $F(\zeta_{\ell})$ is not contained in $F(\brgm)$), since $H$ acts on $\br(\fg_{\mu})(1)$ via $\kbar$. 
\endproof
In light of Lemma \ref{cps}, we would like to understand the intersection $F(\zeta_{\ell}) \cap F(\br(\fg_{\mu}))$; this will also prove important in satisfying condition (\ref{avoid}) of \S \ref{formalglobal}. Note that there is a sandwich
\[
\overline{G}^{\mr{sc}}_{\mu}(k')/(\overline{G}^{\mr{sc}}_{\mu}(k') \cap Z_G(k')) \subseteq \mathbb{P} \left( \br(\gal{F})\right) \subseteq G(k')/Z_G(k'),
\]
where $\mathbb{P}\left(\br(\gal{F})\right)$ denotes the `projective image' of $\br$.
The maximal abelian quotient of $\mathbb{P}\left( \br(\gal{F}) \right) \cong \Gal(F(\br(\fg_{\mu}))/F)$ has order dividing $\# G(k')/[\overline{G}^{\mr{sc}}_{\mu}(k'),\overline{G}^{\mr{sc}}_{\mu}(k')] Z_G(k')$, hence order dividing (recall Assumption \ref{whatsell})
\[
\# G(k')/\left( Z_G(k') \cdot \overline{G}^{\mr{sc}}_{\mu}(k') \right).
\]  
This latter group in turn has order dividing 
\[
\# H^1(\gal{k'}, Z_{G^{\mr{sc}}_{\mu}}) \vert \# Z_{G^{\mr{sc}}_{\mu}}.
\]
Let $L$ be the abelian extension $F(\mu_{\ell}) \cap F(\brgm)$ of $F$, so that $\br(\gal{L})$ contains the commutators 
$[\br(\gal{F}), \br(\gal{F})] \supseteq \overline{G}^{\mr{sc}}_{\mu}(k')$. By the preceding calculation, $[L:F]$ divides $\# Z_{G^{\mr{sc}}_{\mu}}$, so as long as $[F(\mu_{\ell}):F]$ exceeds $\# Z_{G^{\mr{sc}}_{\mu}}$, $F(\mu_{\ell})$ cannot be contained in $F(\brgm)$; condition (\ref{cohvanish}) of \S \ref{formalglobal} therefore follows from the hypotheses of Theorem \ref{maximaltheorem}.

We now treat condition (\ref{avoid}) of \S \ref{formalglobal}. We claim that there is a regular semisimple element $x \in \br(\gal{L}) \cap G_{\mu}(\mathbb{F}_{\ell})$, contained in a torus $T$ of $G$, and a simple root $\alpha$ of $T$, such that $\alpha(x)$ lies in $\Gal(F(\mu_{\ell})/L)$ when we identify $\kbar \colon \Gal(F(\mu_{\ell})/F) \xrightarrow{\sim} \mathbb{F}_{\ell}^\times$. Once we have found such an $x$, we will be able to satisfy the condition (\ref{avoid}) of \S \ref{formalglobal}. The sharpest version of the following lemma would require a case-by-case analysis in the Dynkin classificaton. Since we have no particular need to optimize the set of allowable $\ell$, we content ourselves with a crude bound that works for any group; moreover, the argument in the present form will be reused in \S \ref{principalSL2}.
\begin{lemma}\label{findsselt}
Retain the hypotheses of Theorem \ref{maximaltheorem}. Then we can satisfy condition (\ref{avoid}) of \S\ref{formalglobal}.
\end{lemma}
\proof
By assumption, $\Gal(F(\mu_{\ell})/F)$ is cyclic of order $\ell-1$, and we have just seen that the subgroup $\Gal(F(\mu_{\ell})/L)$ contains all $\# Z_{G^{\mr{sc}}_{\mu}}$-powers in $\mathbb{F}_{\ell}^\times$. For any generator $g$ of $\mathbb{F}_{\ell}^\times$, let $t$ be either $g^{\frac{\# Z_{G^{\mr{sc}}_{\mu}}}{2}}$ or $g^{\# Z_{G^{\mr{sc}}_{\mu}}}$ depending on whether $\# Z_{G^{\mr{sc}}_{\mu}}$ is even or odd. The element $x= 2 \rho^\vee(t)$ belongs to $\overline{G}^{\mr{sc}}_{\mu}(\mathbb{F}_{\ell})$, hence to $\br(\gal{L})$; it is regular since for all positive roots $\alpha$, $\alpha(x)= t^{2 \mr{ht}(\alpha)}$, the maximum height of a root is $h-1$, and none of $t^2, t^4, \ldots, t^{2h-2}$ equals 1, by our assumption on $\ell$. Finally, note that for any simple root $\alpha$, $\alpha(x)= t^2$ belongs to $\Gal(F(\mu_{\ell})/L)$. Thus, for some $\tau \in \gal{L}$, $\br(\tau)=x$, and we can find $\sigma \in \Gal(K/L)$ such that 
\begin{align*}
\sigma &\mapsto \tau|_{F(\brgm)} \in \Gal(F(\brgm)/L), \\
\sigma &\mapsto \alpha(x)=t^2 \in \Gal(F(\mu_{\ell})/L).
\end{align*}
We will now take any lift $\tilde{\sigma}$ of $\sigma$ to $\gal{L}$, and any simple root $\alpha$, to satisfy the first part of condition (\ref{avoid}) of \S \ref{formalglobal}. For the rest, observe that since $K_{\phi}$ and $K_{\psi}$ are both Galois over $F$ itself, $k[\phi(\gal{K})]$ and $k[\psi(\gal{K})]$ are in fact $k[\gal{F}]$-submodules of $\br(\fg_{\mu})$; by irreducibility, they both must equal the whole $\fg_{\mu}$. For our chosen $\tilde{\sigma}$, and any choice of simple root $\alpha$, all of condition (\ref{avoid}) is then satisfied. 
\endproof
Finally, we handle condition (\ref{lindisjoint}); again, our goal has not been to find sharp bounds on $\ell$, but to find a simple, uniform argument.
\begin{lemma}\label{killbigimage}
Retain the hypotheses of Theorem \ref{maximaltheorem}. Then we can satisfy condition (\ref{lindisjoint}) of \S \ref{formalglobal}. 
\end{lemma}
\proof

We must show that $K_{\phi}$ and $K_{\psi}$ are linearly disjoint over $K$. Let 
\[
\Lambda \colon \Rep_k(\gal{F}) \to \Rep_{\fl}(\gal{F})
\]
denote the forgetful functor. If $V$ is any irreducible $k[\gal{F}]$-module, and if $W$ is an irreducible $\fl[\gal{F}]$-submodule of $\Lambda(V)$, then for some integer $r$, $\Lambda(V) \cong W^{\oplus r}$ (by adjunction), and then necessarily $\Lambda(V(1)) \cong W(1)^{\oplus r}$. Taking $V= \brgm$, if the common $\fl[\gal{F}]$-subquotient $\Gal(K_{\phi} \cap K_{\psi}/K)$ of $\Lambda(\brgm)$ and $\Lambda(\brgm(1))$ were non-trivial, then we could conclude that $\Lambda(\brgm)$ and $\Lambda(\brgm(1))$ were isomorphic.
It therefore suffices to produce an element of $\gal{F}$ that acts with different eigenvalues in these two representations (over $\mathbb{F}_{\ell}$). Let $\alpha$ be a root of $G$ with respect to a split maximal torus over $\fl$, and consider an element $x= \alpha^\vee(c)$ for some $c \in \mathbb{F}_{\ell}^\times$ such that $c^2 \in \Gal(F(\mu_{\ell})/L)$; as in the proof of Lemma \ref{findsselt}, we can find $\sigma \in \gal{L}$ such that $\br(\sigma)\in x \cdot Z_G(k)$ and $\kbar(\sigma)=c^2$. For all roots $\beta$ of $G$, $|\langle \alpha^\vee, \beta \rangle| \leq 3$, so for $\sigma$ to act with different eigenvalues on $\Lambda(\brgm)$ and $\Lambda(\brgm(1))$, it suffices for the order of $c$ to be greater than 8. We can arrange this all with $c= g^{\# Z_{G^{\mr{sc}}_{\mu}}}$ for a generator $g$ of $\mathbb{F}_{\ell}^\times$, as long as $\ell>1+8 \cdot\# Z_{G^{\mr{sc}}_{\mu}}$.
\endproof
The hypotheses of Theorem \ref{maximaltheorem} therefore enable us to satisfy conditions (\ref{dualinvs})-(\ref{avoid}) of \S \ref{formalglobal}, so we can invoke Proposition \ref{formalliftingtheorem} to complete the proof of Theorem \ref{maximaltheorem}.
\endproof
\section{The principal $\mr{SL}_2$}\label{principalSL2}
\subsection{Another lifting theorem}\label{principalsl2.1}
The aim of this section is to show that the axiomatized argument of \S \ref{formalglobal} continues to apply when $\br$ is of the form
\[
\gal{F} \xrightarrow{\bar{r}} \mr{GL}_2(k) \xrightarrow{\varphi} G(k) 
\]
where $\varphi$ is (an extension to $\mr{GL}_2$ of) the principal homomorphism $\varphi \colon \mr{SL}_2 \to G$. For background on the principal $\mr{SL}_2$, we can do no better than refer to the lucid expositions of \cite{gross:principalsl2} and \cite{serre:principalsl2}, but here we will recall what is necessary to fix our notation. We begin by recalling the situation in characteristic zero (\cite[\S 2.3]{serre:principalsl2}). Fix a Borel $B$ containing a (split) maximal torus $T$ of $G$, with corresponding base $\Delta$ of the root system, and also fix a pinning 
\[
\{u_{\alpha} \colon \mathbf{G}_a \xrightarrow{\sim} U_{\alpha}\}_{\alpha \in \Delta};
\]
here $U_{\alpha}$ is of course the root subgroup in $B$ corresponding to $\alpha$. Setting $X_{\alpha}= d u_{\alpha}(1)$ for all $\alpha \in \Delta$, we obtain a principal (regular) nilpotent element
\[
X= \sum_{\alpha \in \Delta} X_{\alpha}.
\]
$X$ can be extended to an $\mf{sl}_2$-triple $(X, H, Y)$ inside $\fg$ as follows. For each $\alpha$, let $H_{\alpha}$ be the coroot vector corresponding to $\alpha$ (i.e. $H_{\alpha}= d \alpha^\vee (1)$), and define $Y_{\alpha} \in \fg_{-\alpha}$ (uniquely) by requiring $(X_{\alpha}, H_{\alpha}, Y_{\alpha})$ to be an $\mf{sl}_2$-triple. Then define
\begin{align}
&H= \sum_{\alpha >0} H_{\alpha}= \sum_{\alpha \in \Delta} c_{\alpha} H_{\alpha} \\
&Y= \sum_{\alpha \in \Delta} c_{\alpha} Y_{\alpha};
\end{align}
here the $c_{\alpha}$ are integers determined by the first equation. The resulting homomorphism $\mf{sl}_2 \to \mf{g}$ then uniquely lifts to a homomorphism
\[
\varphi \colon \mr{SL}_2 \to G,
\]
called the principal $\mr{SL}_2$. By construction, $\varphi \begin{pmatrix} t & \\ & t^{-1} \end{pmatrix}= 2 \rho^\vee(t)$, so for any $G$ we can extend $\varphi$ to a homomorphism
\begin{equation}\label{enlargeG}
\varphi \colon \mr{GL}_2= \frac{\mr{SL}_2 \times \mathbf{G}_m}{\langle (-1, -1 )\rangle} \xrightarrow{\varphi \times \id} \frac{G \times \mathbf{G}_m}{\langle (2 \rho^\vee(-1), -1 ) \rangle}.
\end{equation}
Denote by $G_1$ the enlarged target group,\footnote{We remark that this enlargement of $G$ is frequently technically convenient: for instance, it is the Tannakian group appearing when one studies the `geometric' Satake correspondence over finite fields or number fields.} so we have a principal homomorphism $\varphi \colon \mr{GL}_2 \to G_1$; of course, the derived group of $G_1$ is still $G_{\mu}$. Note that the cocharacter $\rho^\vee$ of $G^{\mr{ad}}$ lifts to $G_1$,\footnote{Namely, for the torus $T_1= (T \times \mathbf{G}_m)/\langle 2 \rho^\vee(-1), -1 \rangle$, 
\[
(\rho^\vee, \frac{\beta}{2}) \in X_{\bullet}(T_1) \otimes_{\Z} \Q,
\]
where $\beta$ denotes a generator of $X_{\bullet}(\mathbf{G}_m)$, in fact defines an element of $X_{\bullet}(T_1)$.
}
and if $G$ is any group for which the co-character $\rho^\vee$ of $G^{\mr{ad}}$ lifts to a co-character $\tilde{\rho}^\vee$ of $G$,
\footnote{This lift is of course not unique; for any choice, $2 \tilde{\rho}^\vee$ differs from $2\rho^\vee$, the usual co-character of $G^{\mr{der}}$, by an element of $X_{\bullet}(Z_G)$.} then we can extend the principal $\mr{SL}_2$ to a principal $\mr{GL}_2$ by setting $\varphi \begin{pmatrix} z & \\ & 1 \end{pmatrix}= \tilde{\rho}^\vee(z)$.

The upshot is that for the remainder of this section, we will assume that the co-character $\rho^\vee$ of $G^{\mr{ad}}$ lifts to $G$; we have seen that any $G$ can be embedded in a group with this property, without changing the derived group. We fix such a lift of $\rho^\vee$ as well as the corresponding principal $\mr{GL}_2$
\[
\varphi \colon \mr{GL}_2 \to G.
\]

A crucial piece of structure theory for us will be the decomposition, due to Kostant (see \cite[Proposition 5.2]{gross:principalsl2}), of $\fg_{\mu}$ as $\mf{sl}_2$-module. Let $P$ be the centralizer of $X$. The action of $H$ on $\fg_{\mu}$ preserves $P$, so there is a grading by $H$-eigenvalues 
\[
P= \oplus_{m >0} P_{2m};
\]
that the eigenvalues are even integers follows from standard $\mf{sl}_2$-theory. More precisely, Kostant showed:
\begin{prop}[\cite{kostant:principalsl2}]\label{kostant}
The centralizer $P$ is an abelian subalgebra of $\mf{g}_{\mu}$ of dimension equal to the rank of $\fg_{\mu}$. Letting $\mr{GL}_2$ act on $\mf{g}_{\mu}$ via $\varphi$, there is an isomorphism of $\mr{GL}_2$-representations
\[
\mf{g}_{\mu} \cong \bigoplus_{m >0} \Sym^{2m}(k^2)\otimes {\det}^{-m} \otimes P_{2m},
\]
and $P_{2m}$ is non-zero if and only if $m$ is an exponent of $G$. In particular, the maximal such $m$ is $h-1$.
\end{prop}
Now, the whole theory of the principal $\mr{SL}_2$ works $\ell$-integrally for $\ell >>0$:
\begin{lemma}[\S 2.4 of \cite{serre:principalsl2}]\label{integralsl2}
If $\ell \geq h$, the homomorphism $\varphi \colon \mr{SL}_2 \to G$ is defined over the localization $\Z_{(\ell)}$.
\end{lemma}
\begin{lemma}\label{kostantmodl}
If $\ell \geq 2h-1$, then the decomposition of Proposition \ref{kostant} continues to hold over $\Z_{(\ell)}$.
\end{lemma}
\proof
Let $S$ denote the (diagonal) torus of $\mr{SL}_2$, and write $X^\bullet(S)= \Z \chi$. Then the highest weight of $S$ acting on $\mf{g}_{\mu}$ is $2h-2$ ($h-1$ is the height of the highest root), so \cite[\S 2.2, Proposition 2]{serre:semisimpletensor} implies that $\mf{g}_{\mu}$ is semi-simple for $\ell \geq 2h-1$; \cite[\S 2.2 Remarque]{serre:semisimpletensor} moreover implies the decomposition of $\mf{g}_{\mu}$ is characteristic $\ell$ is just the reduction of the usual decomposition in characteristic zero.
\endproof
We now come to the main result of this section; note that, again, the bounds on $\ell$ in the following result can be somewhat sharpened, but I don't believe in a way that would justify the added complexity. The argument applies uniformly to all $G$ whose Weyl group contains $-1$, until the final step, the group-theoretic Lemma \ref{magmacalc}, which we have checked only for the exceptional groups. 
\begin{thm}\label{principalsl2lift}
Let $G$ be a connected reductive group to which $\rho^\vee$ lifts and whose adjoint form is simple of type $\mr{G}_2$, $\mr{F}_4$, $\mr{E}_7$, or $\mr{E}_8$ (in particular, $-1 \in W_G$), and let $\ell$ be a rational prime greater than $4h-1$. In type $\mr{E}_8$, also exclude $\ell= 229, 269, 367$. Let $F$ be a totally real field for which $[F(\zeta_{\ell}):F]= \ell -1$, and let $\bar{r} \colon \gal{F} \to \mr{GL}_2(k)$ be a continuous representation. Let $\br$ be the composite
\[
\xymatrix{
\gal{F} \ar[r]^{\bar{r}} \ar@/^2 pc/[rr]^{\br} & \mr{GL}_2(k) \ar[r]^{\varphi} & G(k).
} 
\]
Assume that $\bar{r}$ satisfies:
\begin{enumerate}
\item For some subfield $k' \subset k$, 
\[
\mr{SL}_2(k') \subset \bar{r}(\gal{F}) \subset k^\times \cdot \mr{GL}_2(k');
\]
\item $\bar{r}$ is odd;
\item for each $v \vert \ell$, $\bar{r}|_{\gal{F_v}}$ is ordinary, satisfying
\[
\bar{r}|_{I_{F_v}} \sim \begin{pmatrix}
\chi_{1, v} & * \\
0 & \chi_{2, v} \\
\end{pmatrix},
\]
where $\left( \chi_{1,v}/\chi_{2,v}\right)|_{I_{F_v}}= {\kbar}^{r_v}$ for an integer $r_v \geq 2$ such that $\ell>r_v(h-1)+1$;
\item for all primes $v \nmid \ell$ at which $\bar{r}$ is ramified, we can find liftable local deformation conditions $\mc{P}_v$ for $\br$ such that $\dim L_v= h^0(\gal{F_v}, \brgm)$.
\end{enumerate}
Then there exists a lift
\[
\xymatrix{
& G(\mc{O}) \ar[d] \\
\gal{F} \ar@{-->}[ur]^{\rho} \ar[r]^{\br} & G(k)
}
\]
such that $\rho$ is of type $\mc{P}_v$ at all primes $v$ at which $\br|_{\gal{F_v}}$ is ramified. 
\end{thm}
\proof
We proceed one by one through the conditions (\ref{dualinvs})-(\ref{avoid}) of \S \ref{formalglobal}.

Condition (\ref{dualinvs}) demands that $H^0(\gal{F}, \brgm)$ and $H^0(\gal{F}, \brgm(1))$ are both zero; the former vanishing moreover implies that $\Def_{\br}$ is representable. By Lemma \ref{kostantmodl}, we have to check that, for all $m$ such that $P_{2m} \neq 0$, 
\[
H^0(\gal{F}, (\Sym^{2m} \otimes {\det}^{-m})(\bar{r}))=0 
\]
and 
\[
H^0(\gal{F}, (\Sym^{2m} \otimes {\det}^{-m})(\bar{r})(1))=0.
\]
This is obvious, since for $\ell \geq 2h-1$, all such $\Sym^{2m}(k^2)$ are irreducible $k[\mr{SL}_2(k')]$-modules.

Next we check that $\br$ is odd, i.e. that condition (\ref{infiniteplaces}) of \S \ref{formalglobal} is satisfied. Since we have assumed $\bar{r}$ is odd, $\br(c_v)$ is conjugate to $\varphi \begin{pmatrix} -1 & \\ & 1 \end{pmatrix}$ for all $v \vert \infty$. 
We have already noted (Lemma \ref{splitcartan}) that $\rho^\vee(-1)$ is a split Cartan involution of $\fg_{\mu}$, since $-1 \in W_G$, so $\br$ is odd.

Next we treat item (\ref{cohvanish}) of \S \ref{formalglobal}, the cohomological vanishing result. To begin, we record the following elementary lemma in group cohomology, which follows (just as in Lemma \ref{cps}) from \cite[2.9, 3.3]{cline-parshall-scott}:
\begin{lemma}\label{cpssl2}
Assume $\#k' \not \in \{3,5,9\}$. Fix an integer $r < \ell$. Then $H^1(\mr{SL}_2(k'), \Sym^r(k^2))=0$.
\end{lemma}
Again applying Kostant's result (Lemma \ref{kostantmodl}), we deduce from Lemma \ref{cpssl2} that
\[
H^1(\mr{SL}_2(k'), \varphi(\fg_{\mu}))= \bigoplus_m H^1(\mr{SL}_2(k'), \Sym^{2m}(k^2))^{\oplus \dim P_{2m}}=0.
\]
By repeated inflation-restriction as in Lemma \ref{cps}, it follows easily that (recall $K= F(\brgm, \mu_{\ell})$)
\[
H^1(\Gal(K/F), \brgm)=0.
\]
Also by the argument of Lemma \ref{cps}, we deduce that $H^1(\Gal(K/F), \brgm(1))=0$, as long as $\bar{\kappa}|_{H} \neq 1$: here as before we let $H= \gal{F(\brgm)}= \gal{F(\bar{r}(\mf{sl}_2))}$. But now the discussion after Lemma \ref{cps} (using the group $\mr{GL}_2$ instead of $G$) applies, and we see that to ensure $\kbar|_H \neq 1$ it suffices to take $[F(\zeta_{\ell}):F]>2$.

We will next construct an element $\sigma \in \gal{F}$ such that $\br(\sigma)$ is regular semi-simple and satisfies $\alpha(\br(\sigma))= \bar{\kappa}(\sigma)$ for any simple root $\alpha$ (of the unique maximal torus containing $\br(\sigma)$), arguing just as in Lemma \ref{findsselt}. First note that $F(\brgm)= F(\ad^0 \bar{r})$, and that the assumption $\bar{r}(\gal{F}) \supset \mr{SL}_2(k')$ implies $[F(\ad^0 \bar{r}) \cap F(\zeta_{\ell}):F] \leq 2$. For all $a \in \mathbb{F}_{\ell}^\times$, $\varphi \begin{pmatrix} a & \\ & a^{-1} \end{pmatrix}= 2 \rho^\vee(a) \in \br(\gal{F})$, and this element is regular if and only if for all positive roots $\alpha$, $a^{\langle 2 \rho^\vee, \alpha \rangle} \neq 1$. Now, $\langle 2 \rho^\vee, \alpha \rangle$ is equal to twice the height of the root $\alpha$; its maximum value is $2h-2$, so consider any element $a \in \fl^\times$ with order greater than $2h-2$. Set $L= F(\mu_{\ell}) \cap F(\ad^0 \bar{r})$, so that $\bar{r}(\gal{L})$ contains the commutators of $\bar{r}(\gal{F})$, hence contains $\mr{SL}_2(\mathbb{F}_{\ell})$. Let $\sigma \in \gal{L}$ be an element mapping to $x= \begin{pmatrix} a & \\ & a^{-1} \end{pmatrix}$, so that $\br(\sigma)$ is regular semi-simple. For any simple root $\alpha$, $\alpha(\br(\sigma))=a^2$; to ensure that we can arrange this to equal $\bar{\kappa}(\sigma)$, we need
\[
a^2 \in \Gal(F(\zeta_{\ell})/L) \subset \mathbb{F}_{\ell}^\times.
\]
But $\Gal(F(\zeta_{\ell}/L)$ has index at most two in $\Gal(F(\zeta_{\ell})/F) \xrightarrow[\kbar]{\sim} \mathbb{F}_{\ell}^\times$, hence contains all squares.

We will use the $\sigma$ just constructed to satisfy condition (\ref{avoid}) of \S \ref{formalglobal}. Recall that we are given non-zero classes $\phi \in H^1_{\mc{P}^\perp}(\gal{F, \Sigma}, \brgm(1))$ and $\psi \in H^1_{\mc{P}}(\gal{F, \Sigma}, \brgm)$, and that we have a maximal torus $T$ containing $\br(\sigma)$. We will find a simple root $\alpha$ such that $k[\psi(\gal{K})]$ has non-zero $\mf{l}_{\alpha}$ component, and $k[\phi(\gal{K})]$ has nonzero $\mf{g}_{-\alpha}$ component. This suffices to establish condition (\ref{avoid}), since the regular semi-simple element $\br(\sigma)$ constructed above satisfies $\alpha(\br(\sigma))= \kbar(\sigma)$ for any simple root $\alpha$. We again appeal to Kostant's result; since $k[\psi(\gal{K})]$ (respectively, $k[\phi(\gal{K})]$) is a $k[\gal{F}]$-submodule of $\brgm$ (respectively, $\brgm(1)$), we know that each contains one of the summands $\Sym^{2m}(k^2)$ under the action of our principal $\mf{sl}_2$ on $\fg_{\mu}$. We will use this fact to construct the desired $\alpha$. An analogue of the following lemma surely holds for any simple type, but I do not have a satisfactory general argument. It is easy to check for any particular group by explicit computation, and we do this here for the exceptional types (the case of $\mr{E}_6$ is recorded here for later use: see Theorem \ref{principalsl2liftout}).
\begin{lemma}\label{magmacalc}
Let $\mf{g}$ be a simple Lie algebra of type $\mr{G}_2$, $\mr{F}_4$, $\mr{E}_7$, or $\mr{E}_8$ equipped with a principal homomorphism $\varphi \colon \mf{sl}_2 \to \mf{g}$. Let $m$ be an exponent of $\mf{g}$, so that there is a unique summand $\Sym^{2m}(k^2)$ of $\mf{g}_{\mu}$ when regarded as $\mf{sl}_2$-representation. This summand $\Sym^{2m}(k^2)$ contains a unique line inside $\oplus_{\alpha \in \Delta} \mf{g}_{-\alpha}$, and the projection of this line to each (negative) simple root space $\mf{g}_{-\alpha}$ is non-zero except in the following cases:
\begin{itemize}
\item type $\mr{G}_2$: $\ell \in \{2,3,5 \}$;
\item type $\mr{F}_4$: $\ell \in \{2,3,5,7,11 \}$;
\item type $\mr{E}_7$: $\ell \in \{2,3,5,7,11,13,17,19,31,37,53 \}$;
\item type $\mr{E}_8$: $\ell \in \{2,3,5,7,11,13,17,19,23,29,61,67,71,97,103,109,229,269,397 \}$;
\end{itemize}

Now asume $\mf{g}$ is of type $\mr{E}_6$. If $m \neq 4, 8$, then the line inside $\oplus_{\alpha \in \Delta} \mf{g}_{-\alpha}$ of the summand $\Sym^{2m}(k^2)$ has non-zero projection to each $\mf{g}_{-\alpha}$ except when $\ell \in \{2,3,5,7,11\}$. For $m \in \{4, 8\}$, there is a simple root $\alpha_1$ such that this projection is non-zero in $\mf{g}_{-\alpha_1}$, and for all exponents $n$ of $\mf{g}$, the projection of $\mf{t} \cap \Sym^{2n}(k^2) \subset \mf{g}$ to $\mf{l}_{\alpha_1}$ is non-zero, again except when $\ell \in \{2,3,5,7,11\}$ (in fact, any simple root not fixed by the outer automorphism of $\mr{E}_6$ works here).  
\end{lemma}
\proof 
We describe the algorithm for checking this; it is especially easy to implement using any computational software that has built-in Chevalley bases for the simple Lie algebras, eg GAP or Magma. I carried the calculation out in Magma, so I have noted in what follows how Magma normalizes some indexing, etc.
\begin{enumerate}
\item Let $l$ be the rank of $G$. Construct a simple Lie algebra of type $G$ with a Chevalley basis 
\[
\{x[i], y[j], h[k]\}_{1 \leq i, j \leq |\Phi^+|, 1 \leq k \leq l};
\]
for a suitable bijection $\Phi^+ \xrightarrow[f]{\sim} \{1, \ldots, |\Phi^+|\}$ making $x[f(\alpha)]$ a basis of $\mf{g}_{\alpha}$; arrange this bijection so that $f(\Delta)= \{1, \ldots, l\}$.\footnote{Magma automatically does this when it produces a Chevalley basis.} Conventions for Chevalley bases are not universal; Magma's bracket relations include
 \[
 \left[y[f(\alpha)], x[f(\alpha)]\right]= \alpha^\vee\] 
(not $-\alpha^\vee$) and 
\[
\left[x[f(\alpha)], h\right]= \alpha(h)x[f(\alpha)]
\] 
for all $h \in \mf{t}$, the Cartan sub-algebra spanned by the $h[k]$. The elements $h[k]$ are not the coroots but rather are the dual basis:
\[
[x[i], h[j]]= \delta_{ij} x[i]
\]
for all $1 \leq i, j \leq l$. Since $\ell$ is very good for $G$, these span the same $k$-subspace as the coroot vectors.
\item Then we define elements
\begin{itemize}
\item $X= \sum_{i=1}^l x[i]$ (a regular unipotent);
\item $H= \sum_{i=1}^{|\Phi^+|} \left[y[i], x[i]\right]= 2 \rho^\vee$, so that $[X, H]= 2X$;
\item $Y= \sum_{i=1}^l c[i] \cdot y[i]$ where the $c[i]$ are the appropriate structure constants to make $\{X, H, Y\}$ an $\mf{sl}_2$-triple in the following sense:\footnote{The $c[i]$ can be found in the tables of \cite{bourbaki:lie456}. Within Magma, they can be derived, for our Lie algebra $\mf{g}$, by computing $\mr{rd}:= \mbf{RootDatum}(\mf{g})$, then $A:=\mbf{SimpleCoroots}(\mr{rd})$ (an $l \times l$ matrix), then by forming a vector $c$ whose entries are twice the sum of the rows of $A^{-1}$; the $i^{th}$ entry of $c$ is the desired $c[i]$. But we double-check in each case the bracket properties for $X, Y, H$, so how the $c[i]$ are arrived at no longer matter.} $[X,H]=2X$, $[Y,H]=-2Y$, $[Y, X]=H$. In terms of a Chevalley basis $\langle e, f, h \rangle$ of $\mf{sl}_2$, satisfying $[e,h]=e$, $[f, h]=-f$, $[f, e]=2h$, we then get a principal $\mf{sl}_2$ (defined over $\Z[\frac{1}{2}]$)
\begin{align*}
\mf{sl}_2 \to \mf{g} \\
h \mapsto \frac{H}{2}\\
x \mapsto X\\
y \mapsto Y.
\end{align*}
\end{itemize}
\item Compute the centralizer $P$ of the element $X$ in $\mf{g}$, i.e. find a basis in terms of the Chevalley basis.
\item $H$ preserves $P$, acting semi-simply, so replace the above basis of $P$ with a basis of $H$-eigenvectors.\footnote{In Magma, when you call $P:= \textbf{Centraliser}(\mf{g}, X)$ and then $\textbf{ExtendBasis}(P, \mf{g})$, it hands you a basis of $\mf{g}$, in terms of the Chevalley basis, whose first $l$ entries are a basis of $P$, and in fact already an $H$-eigenbasis.} Let us call this eigen-basis $p[1], \ldots, p[l]$, ordered so that when we write $m_1 \leq m_2 \leq \cdots \leq m_l$ for the exponents of $\mf{g}$, 
\[
\left[p[i], H\right]= 2m_i p[i]
\]
for $i=1, \ldots, l$ (we remark that $p[1]=X$).
\item We then simply compute the Lie brackets, for each $i = 1, \ldots, l$,
\[
\ad(Y)^{m_i+1}(p[i]) \in \bigoplus_{\alpha \in \Delta} \mf{g}_{-\alpha};
\]
in each case we get some $\Z$-linear combination of $y[1], \ldots, y[l]$, and we record (as the exceptional cases) all primes $\ell$ dividing one of these coefficients.
\end{enumerate}
The above argument in type $\mr{E}_6$ fails for all primes $\ell$ (i.e. in characteristic zero, one of the integers recorded in the previous sentence is zero), but we can replace it with a simple variant. The exponents of $\mr{E}_6$ are $\{1,4,5,7,8,11\}$, and the above argument works except for the exponents $m_2=4$ and $m_5=8$; that is, for $i \in \{1,3,4,6\}$, $\ad(Y)^{m_i +1}(p[i])$ has non-zero projection to every $\mf{g}_{-\alpha}$, except in characteristics $\ell \in \{2,3,5,7,11\}$. For $m_2=4$ and $m_5=8$, $\ad(Y)^{m_i+1}(p[i])$ has non-zero projection to the $\mf{g}_{-\alpha_i}$ root space for $i \in \{1,3,5,6\}$,\footnote{Intrinsically, $\alpha_2$ and $\alpha_4$ are those fixed under the outer automorphism of $\mr{E}_6$. Our labeling convention is $f(\alpha_i)=i$, where $x[1], \ldots, x[6]$ is the ordered set of simple roots produced by Magma.} at least if $\ell \not \in \{2,3,5,7,11\}$. We then show that for any exponent $m_j$, the simple root $\alpha_1$ is non-vanishing on the line 
\[
\mf{t} \cap \Sym^{2m_j}(k^2) \subset \mf{g}.
\]
That is, the $\mf{l}_{\alpha_1}$ component of this summand is non-zero, while $\Sym^{2m_i}(k^2)$ has non-zero $\mf{g}_{-\alpha_1}$ component for $i=2,5$. To do this check, we compute the $h[1]$-component of
\[
\ad(Y)^{m_j}(p[j]),
\]
and again record the rational primes dividing the output.
\endproof
\proof[Conclusion of the proof of Theorem \ref{principalsl2lift}] We can therefore satisfy condition (\ref{avoid}) with any simple root $\alpha$ for which $\alpha \left( \mf{t} \cap k[\psi(\gal{K})] \right) \neq 0$ (there is always at least one such $\alpha$, since this intersection is non-zero; in fact, further calculation shows that any simple $\alpha$ will work). 

To complete the proof of the theorem, we address condition (\ref{lindisjoint}) of \S \ref{formalglobal}; that is, we must show that for non-zero Selmer classes $\phi$ and $\psi$, the fixed fields $K_{\phi}$ and $K_{\psi}$ are linearly disjoint over $K$. As before (Lemma \ref{killbigimage}), we must rule out the possibility that $\brgm$ and $\brgm(1)$ have a common $\fl[\gal{F}]$-subquotient. By Lemma \ref{kostantmodl} (and in particular semi-simplicity), it suffices to show that, for all pairs of integers $m$ and $n$ such that $P_{2m}$ and $P_{2n}$ are non-zero, $(\Sym^{2m} \otimes \det^{-m} )(\bar{r})$ and $(\Sym^{2n} \otimes \det^{-n} )(\bar{r})(1)$ have no common $\fl[\gal{F}]$-submodule. Recall that since $\ell \geq 2h-1$ these are irreducible $k[\gal{F}]$-modules. For any two irreducible $k[\gal{F}]$-modulues $V_1$ and $V_2$ such that $\Lambda(V_1)$ and $\Lambda(V_2)$ have a common sub-quotient, there is some $\iota \in \Gal(k/\fl)$ such that $V_1 \otimes_{k, \iota} k$ is isomorphic to $V_2$.\footnote{Note that $\Lambda$ is right-adjoint to $\otimes_{\fl} k$, and $k \otimes_{\fl} k$ is isomorphic to $\prod_{\iota} k$ by the map $x \otimes y \mapsto (\iota(x)y)_{\iota}$.} Applying this observation to $V_1= (\Sym^{2m} \otimes \det^{-m} )(\bar{r})$ and $V_2= (\Sym^{2n} \otimes \det^{-n} )(\bar{r})(1)$, we are clearly done unless $m=n$. Then, by an earlier argument in the present proof, for any $a \in \fl^\times$ we can find a $\sigma \in \gal{F}$ such that $\bar{r}(\sigma)= \begin{pmatrix} a & \\ & a^{-1} \end{pmatrix}$ and $\kbar(\sigma)=a^2$. It therefore suffices to ensure that the (multi-)sets of eigenvalues $\{a^{2i}\}_{i \in [-m, m]}$ and $\{a^{2i+2}\}_{i \in [-m, m]}$ are not the same ($a \in \fl^\times$ implies $\iota$ has no effect on the eigenvalues in question). Since $m$ is at most $h-1$, taking $a$ to be a generator of $\fl^\times$ and $\ell$ to be greater than $4h-1$ guarantees this distinctness.
\endproof
We now record some further properties of the lifts produced by Theorem \ref{principalsl2lift}. These will be used in our application to the construction of geometric Galois representations with exceptional monodromy groups.
\begin{lemma}\label{reductivemono}
Let $F$ be a number field, let $G$ be a Chevalley group, and let $\br= \varphi \circ \bar{r}$ for some representation $\bar{r} \colon \gal{F} \to \mr{GL}_2(k)$ with image as in Theorem \ref{principalsl2lift}. Suppose there exists a lift $\rho$ of $\br$ to $G(\mc{O})$ (as in Theorem \ref{principalsl2lift}, for instance). Then the algebraic monodromy group $G_{\rho} \subset G_{E}$ is reductive (recall $E= \Frac{\mc{O}}$).
\end{lemma}
\proof
We may assume $G$ is of adjoint type. Let $R_u \subset G_{\rho}^0$ be the unipotent radical of $G_{\rho}$. We aim to show $R_u$ is trivial. The sequence of $E[\gal{F}]$-stable Lie algebras $\Lie(R_u) \subset \Lie(G_{\rho}) \subset \mf{g}_{ E}$ can be intersected with $\mf{g}_{\mc{O}}$ to give a sequence of $\mc{O}[\gal{F}]$-modules
\[
\Lie(R_u) \cap \mf{g}_{\mc{O}} \subset \Lie(G_{\rho}) \cap \mf{g}_{\mc{O}} \subset \mf{g}_{\mc{O}}
\]
where each of the inclusions has torsion-free cokernel. Reducing to $k$, we get a corresponding sequence, still inclusions, of $k[\gal{F}]$-modules
\[
(\Lie(R_u) \cap \mf{g}_{\mc{O}} )\otimes_{\mc{O}} k \subset (\Lie(G_{\rho}) \cap \mf{g}_{\mc{O}}) \otimes_{\mc{O}} k \subset \mf{g}.
\]
Recall that $\mf{g}$ decomposes into irreducible constituents, as $\mf{sl}_2$-module and hence as $k[\gal{F}]$-module, as $\oplus_{m > 0} \Sym^{2m}(k^2) \otimes P_{2m}$. In particular, any of these summands, and hence any $k[\gal{F}]$-submodule, contains a non-zero semi-simple element of $\mf{g}$, coming from the weight zero component, i.e. the centralizer of the regular element $H$. We conclude that $(\Lie(R_u) \cap \mf{g}_{\mc{O}}) \otimes_{\mc{O}} k$, and therefore $R_u$ itself, must be trivial, since a nilpotent Lie algebra contains no non-zero semi-simple elements.
\endproof
\begin{lemma}\label{avoidsl2}
Let $F$ be a number field, let $G$ be an exceptional group, and suppose that $\rho \colon \gal{F} \to G(\mc{O})$ is a continuous representation such that
\begin{enumerate}
\item $G_{\rho}$ is reductive;
\item $G_{\rho}$ contains a regular unipotent element of $G$;
\item for some $v \vert \ell$, $\rho|_{\gal{F_v}} \colon \gal{F_v} \to B(\mc{O})$ is ordinary, factoring through a Borel subgroup $B$ of $G$, and for some (any\footnote{The set $\{r_{\alpha}\}_{\alpha \in \Delta(B, T)}$ does not depend on the choice of $T$.}) maximal torus $T \subset B$, and all simple roots $\alpha \in \Delta(B, T)$, the composites $\alpha \circ \rho|_{I_{F_v}}$ have the form $\kappa^{r_{\alpha}}$ for distinct integers $r_{\alpha}$.
\end{enumerate}
Then $G_{\rho} \supseteq G^{\mr{der}}$.
\end{lemma}
\proof
We may assume $G$ is adjoint, so our task is to show $G_{\rho}=G$. The essential input is a result of Dynkin (see \cite[Theorem A]{saxl-seitz:dynkin}) establishing that any $G_{\rho}$ satisfying the first two conditions of the lemma must be either a principal $\mr{PGL}_2$ or an embedded $\mr{F}_4= G_{\rho} \into \mr{E}_6=G$ containing a principal $\mr{PGL}_2$. Assume first that $G$ is not of type $\mr{E}_6$. Then it suffices for the integers $r_{\alpha}$, $\alpha \in \Delta(B,T)$, not all to be equal to conclude $G_{\rho} \not \subset \mr{PGL}_2$, hence $G_{\rho}= G$.

Now take $G$ to be of type $\mr{E}_6$. The argument is essentially the same as above, but we make it more precise. By \cite[Theorem 1]{seitz:maximalinexceptional}, all embeddings $\mr{F}_4 \into \mr{E}_6$ are related by $\Aut(\mr{E}_6)$. One such $\mr{F}_4$ can be constructed by fixing a pinning, denoting by $\tau \in \Aut(\mr{E}_6)$ the image of the non-trivial element of $\mr{Out}(\mr{E}_6)$ under the associated section $\mr{Out}(\mr{E}_6) \to \Aut(\mr{E}_6)$, and setting $\mr{F}_4= (\mr{E}_6)^{\tau=1}$. It follows that all embeddings $\mr{F}_4 \into \mr{E}_6$ are $\mr{E}_6$-conjugate to this one, and all have the form $(\mr{E}_6)^{\tau'=1}$ for some involution $\tau'$ of $\mr{E}_6$ inducing its non-trivial outer automorphism. 

Let us suppose, then, that our $\rho \colon \gal{F} \to \mr{E}_6(\mc{O})$ (which from now on we regard as $\Frac(\mc{O})$, or even $\Qlb$-valued) factors through $(\mr{E}_6)^{\tau=1} \cong \mr{F}_4$ for some such involution $\tau$. In particular, $\rho|_{\gal{F_v}}$ factors through $B \cap (\mr{E}_6)^{\tau=1}$. For all Borel subgroups $B'$ of $G$, the character groups $\Hom(B', \mbf{G}_m)$ are \textit{canonically} identified: this assertion combines the conjugacy of Borel subgroups with the fact that they are their own normalizers. For emphasis, denote by $\mc{B}^*$ this canonically-defined group. In particular, pre-composition with $\tau$ defines an automorphism of $\mc{B}^*$. Regarding $\Delta(B, T)$ as a subset of $\mc{B}^*$, $\tau$-invariance of $\rho$ immediately implies that $r_{\alpha}= r_{\tau(\alpha)}$ for all  $\alpha \in \Delta(B, T) \subset \mc{B}^*$. This contradicts the fact that $\tau$ is a non-trivial outer automorphism, since we have assumed the $r_{\alpha}$ are all distinct.

\endproof
\section{Exceptional monodromy groups: the case $-1 \not \in W_G$}\label{-1notinW}
Having shown that the monodromy groups $G_{\rho}$ produced by Theorem \ref{principalsl2lift} are reductive, we would like to arrange that they contain the principal $\mr{SL}_2$ of $G$, i.e. that they contain a regular unipotent element of $G$. I do not know if this is automatic for all $\rho$ produced by Theorem \ref{principalsl2lift}, but the following beautiful result of Serre recommends caution:
\begin{eg}[Theorem 1 of \cite{serre:principalsl2}]\label{serreeg}
Let $G$ be a Chevalley group of adjoint type, with, as always, Coxeter number $h$, and let $\ell$ be a prime number. Let $K$ be any algebraically closed field (eg, $K= \Qlb$). 
\begin{enumerate}
\item If $\ell= h+1$, then there exists an embedding $\mr{PGL}_2(\mathbb{F}_{\ell}) \into G(K)$, except when $h=2$ and $\mr{char}(K)=2$.
\item If $\ell = 2h+1$, then there additionally exists an embedding $\mr{PSL}_2(\mathbb{F}_{\ell}) \into G(K)$.
\end{enumerate}
There is a pleasant parallel with the techniques of the present paper: the idea of Serre's argument is to begin with $K$ of characteristic $\ell$ and the principal $\mr{PGL}_2(\mathbb{F}_{\ell}) \into G(K)$ (for any $\ell \geq h$), and then to try to deform this homomorphism to characteristic zero. This follows from a simple and satisfying group cohomology calculation for $\ell = h+1$, but is rather trickier for $\ell= 2h+1$. \qed
\end{eg}
Our first task, then, is to circumvent examples of this sort.
\subsection{Non-splitness of residual Galois representations of Steinberg type}
In this subsection we explain how to choose our $\bar{r}$ and our deformation problem for $\br= \varphi \circ \bar{r}$ to guarantee that $G_{\rho}$ contains a principal $\mr{SL}_2$. The basic idea is to choose $\bar{r}$ having an auxiliary prime $v \not \mid \ell$ of ramification `of Steinberg type,' and then to use the local condition of \S \ref{steinberg}. The subtlety is that we will need (in order to get non-trivial unipotent elements in the image) to guarantee that the local lift $\rho|_{\gal{F_v}} \colon \gal{F_v} \to G(\mc{O})$ is \textit{ramified}; this is very difficult to force through purely Galois-theoretic argument unless the original $\br|_{\gal{F_v}}$ was itself already ramified.\footnote{There are some elaborate Galois-theoretic means under hypotheses that are in practice too difficult to arrange; if the characteristic zero lifts are known to be \textit{automorphic}, then suitable cases of the Ramanujan conjecture imply this non-splitness--this observation should suggest the difficulty of the problem. Current potential automorphy techniques would only help with the case $G= \mr{G}_2$.} 

We will therefore consider $\bar{r}$ arising as the residual Galois representations associated to cuspidal automorphic representations $\pi$ (corresponding to classical holomorphic modular forms) of $\mr{GL}_2(\mathbb{A}_{\Q})$ such that $\pi_p$ is the twist of a Steinberg representation for some $p$. Given such a $\pi$, there is a number field $E$ and a strongly compatible system of $\ell$-adic representations
\[
r_{\pi, \lambda} \colon \gal{\Q} \to \mr{GL}_2(E_{\lambda})
\]
as $\lambda$ ranges over the finite places of $E$. We will show that for almost all $\lambda$, $\bar{r}_{\pi, \lambda}|_{\gal{\Q_p}}$ is (reducible but) indecomposable. The argument we give will surely also apply in the Hilbert modular case,\footnote{In fact, it will even apply to appropriate unitary groups of any rank, using \cite{blggt:potaut}.} but I unfortunately don't know a reference for the necessary modular-lifting theorem in that context; for our application, the case of $\Q$ suffices. After writing this paper, I learned from Khare that this has been proven by Weston in \cite{weston:unobstructed}. Weston attributes the argument to Ribet, perhaps unsurprisingly, since the technique is a variant of the `level-lowering' trick that reduces Fermat's Last Theorem to the Shimura-Taniyama conjecture.
\begin{prop}[Proposition 5.3 of \cite{weston:unobstructed}]\label{levellower}
Let $\pi$ be as above a cuspidal automorphic representation corresponding to a holomorphic eigenform of weight at least 2. Assume that for some prime $p$ of $\Q$, $\pi_p$ is isomorphic to the twist of a Steinberg representation of $\mr{GL}_2(\Q_p)$. Then for almost all $\lambda$, the local Galois representation $\bar{r}_{\pi, \lambda}|_{\gal{\Q_p}}$ has the form
\[
\bar{r}_{\pi, \lambda}|_{\gal{\Q_p}} \sim \begin{pmatrix} \chi \kbar & * \\ 0 & \chi \end{pmatrix}
\] 
where the extension $* \in H^1(\gal{\Q_p}, \mathbb{F}_{\lambda}(\kbar))$ is non-zero.
\end{prop}
We give two proofs. The first is quite heavy-handed, but shows this proposition is immediate from known results. The second will apply to Hilbert modular forms as well, but I haven't found the necessary reference in the literature to complete it, and a full proof would take us unnecessarily far afield.
\proof[First proof]
Since any $\pi$ for which some $\pi_p$ is Steinberg is non-CM, Ribet's result (see Theorem \ref{ribetimage} above) shows that for almost all $\lambda$, the image $\bar{r}_{\pi, \lambda}(\gal{\Q})$ contains $\mr{SL}_2(\mathbb{F}_{\ell})$, and in particular $\bar{r}_{\pi, \lambda}$ is irreducible. Suppose there is an infinite set $\Lambda$ of $\lambda$ for which $\bar{r}_{\pi, \lambda}|_{\gal{\Q_p}}$ is split. Then for all $\lambda \in \Lambda$, the refined level-aspect of Serre's conjecture (see \cite[Theorem 1.12]{edixhoven:serre} and the references given there) yields a cuspidal representation $\pi(\lambda)$ congruent to $\pi$ modulo $\lambda$, of level strictly less than that of $\pi$, and of weight at most that of $\pi$. The number of cuspidal automorphic representations of bounded weight and level is finite, so infinitely many of these $\pi(\lambda)$ must in fact be the same representation $\pi'$. But then $\pi \equiv \pi' \pmod \lambda$ for infinitely many $\lambda$, which implies $\pi= \pi'$ (by strong multiplicity one). Contradiction.
\endproof
\proof[Sketch of second proof]
Let us now suppose more generally that $\pi$ is a cuspidal Hilbert modular representation over a totally real field $F$, and that the archimedean weights of $\pi$ all have the same parity. By assumption, $\pi_v$ is (twist of) Steinberg for some place $v$ of $F$. Again Ribet's result shows that for almost all $\lambda$, the image $\bar{r}_{\pi, \lambda}(\gal{F})$ contains $\mr{SL}_2(\mathbb{F}_{\ell})$. For such $\lambda$, it follows that the index of a Borel subgroup in $\bar{r}_{\pi, \lambda}(\gal{F})$ is at least $\ell+1$, which is greater than $[F(\zeta_{\ell}):F]$. Thus, for almost all $\lambda$, $\bar{r}_{\pi, \lambda}|_{\gal{F(\zeta_{\ell})}}$ is irreducible. Moreover, for almost all $\lambda$ (say $\lambda \vert \ell$), the local restrictions $\bar{r}_{\pi, \lambda}|_{\gal{F_w}}$ for $w \vert \ell$ are torsion-crystalline in the Fontaine-Laffaille range.

Now consider a putative infinite set $\Lambda$ of $\lambda$ for which $\bar{r}_{\pi, \lambda}|_{\gal{F_v}}$ is split; throwing away a finite number of $\lambda$, we may assume the conclusions of the last paragraph hold for all $\lambda \in \Lambda$. For simplicity (we can always make a global twist), assume $\pi_v$ is an unramified twist of Steinberg. For all places $w \neq v$ of $F$ at which $\pi_w$ is ramified, let $\tau_w$ be the inertial type of $\pi_w$. Then define a global deformation ring $R_{F, \lambda}(\tau)$ classifying deformations of $\bar{r}_{\pi, \lambda}$ that are 
\begin{itemize}
\item fixed determinant equal to $\det (r_{\pi, \lambda})$;
\item crystalline of the same weight as $\pi$ at each $w \vert \ell$ of $F$;
\item type $\tau_w$ for all ramified primes $w \neq v$ for $\pi$;
\item unramified elsewhere.
\end{itemize} 
We claim that $R_{F, \lambda}(\tau)$ has a modular characteristic-zero point. This should follow as in \cite[\S 4.2-4.3]{kisin:cdm} (see also \cite[Theorem 3.14]{kisin:durham} from an appropriate $R= \mathbb{T}$ theorem over some well-chosen finite extension $F'$ of $F$. Kisin's variant of the Taylor-Wiles method reduces this to understanding the geometry of local deformation rings, \textit{modulo having an appropriate generalization of the Skinner-Wiles potential level-lowering arguments} (\cite{skinner-wiles:basechange}), which (see the comment in the proof of \cite[Theorem 3.14]{kisin:durham}) as far as I can tell have only been written down in weight 2. In our case, there are no local difficulties, since at $w \vert \ell$ we are in the Fontaine-Laffaille range, where the torsion-crystalline deformation ring is a power series ring over $\mc{O}_{\lambda}$. In any case, given such an $R= \mathbb{T}$ theorem, the methods of Khare-Wintenberger (see \cite[\S 4.2]{kisin:cdm}) show $R_{F, \lambda}(\tau)$ is a finite $\mc{O}_{\lambda}$-algebra of dimension at least 1, hence that it has a characteristic zero point. Then Kisin's modular lifting method would give a cuspidal $\pi(\lambda)$ with weight and level bounded independently of $\lambda$, and congruent to $\pi$ modulo $\lambda$, such that $\pi(\lambda)_v$ is unramified. The argument of the first proof now applies, and we deduce a contradiction. 
\endproof
\begin{rmk}
For a classical case in which the result of Proposition \ref{levellower} is more explicit, consider the compatible system $\{ r_{E, \ell}\}_{\ell}$ associated to an elliptic curve over $F$ with multiplicative reduction at some place $v$. Then the `Tate parameter' $q_{E_v} \in F_v^\times $ associated to $E_{F_v}$ (namely, $\ell$-divisibility of $v(q_{E_v})$) tells exactly when the extensions $*$ are non-split.
\end{rmk}
\subsection{Application to exceptional monodromy groups}
We can now complete the proof of the main theorem of this section:
\begin{thm}\label{notE6theorem}
Let $G$ be an adjoint Chevalley group of type $G_2$, $F_4$, $E_7$, or $E_8$. Then for a density one set of rational primes $\ell$, we can find $\ell$-adic represenations
\[
\rho_{\ell} \colon \gal{\Q} \to G(\Qlb)
\] 
whose image is Zariski-dense in $G$.
\end{thm}
\proof
Let $f$ be a (non-CM) weight 3 cuspidal eigenform that is a newform of level $\Gamma_0(p) \cap \Gamma_1(q)$ for some primes $p$ and $q$; the nebentypus of $f$ is a character $\chi \colon (\Z/pq\Z)^\times \onto (\Z/q\Z)^\times \to \CC^\times$. Such $f$ exist: for instance, there are such $f$ in $S_3(15)$ (\cite[15.3.1a, 15.3.3a]{lmfdb}). There is a number field $E$ such that for all finite places $\lambda$ of $E$ we have $\lambda$-adic representations $r_{f, \lambda} \colon \gal{\Q} \to \mr{GL}_2(E_{\lambda})$ and semi-simple residual representations 
\[
\bar{r}_{f, \lambda} \colon \gal{\Q} \to \mr{GL}_2(\mathbb{F}_{\lambda}),
\]
where we let $\mathbb{F}_{\lambda}$ denote the residue field of $E$ at $\lambda$. A well-known argument using the Weil bound and the \v{C}ebotarev density theorem (see \cite[Lemma 3.2]{gee:sato-tatewt3}) implies that for a density one set of $\ell$, there is some $\lambda \vert \ell$ such that $r_{f, \lambda}|_{\gal{\Ql}}$ is ordinary. Moreover Proposition \ref{levellower} implies that for all but finitely many $\lambda$, the restriction
\[
\bar{r}_{f, \lambda}|_{\gal{\Q_p}} \colon \gal{\Q_p} \to \mr{GL}_2(\mathbb{F}_{\lambda})
\]
is reducible but indecomposable. Finally, $\bar{r}_{f, \lambda}(I_{\Q_q})$ has order prime to $\ell$ for almost all $\ell$. We then define the composite $\br_{\lambda}$ by
\[
\xymatrix{
\gal{\Q} \ar[r]_-{\bar{r}_{f, \lambda}} \ar@/^2 pc/[rr]^{\br_{\lambda}} & \mr{GL}_2(\mathbb{F}_{\lambda}) \ar[r]_{\varphi} & G(\mathbb{F}_\lambda).
} 
\]
Throwing out a further finite set of primes $\ell$ (those less than $4h$, as in Theorem \ref{principalsl2lift}; note that the integer ``$r_v$'' in the statement of that result is in our case 2, by \cite[Theorem 2.1.4]{wiles:ordinary}), we therefore have a density one set of $\ell$, and for each such $\ell$ a $\lambda \vert \ell$, such that $\br_{\lambda}$ satisfies the hypotheses of Theorem \ref{principalsl2lift}: at the places $p$, $q$, and $\ell$ of ramification we take the following local deformation conditions:
\begin{itemize}
\item at $p$ we use the Steinberg deformation condition of \S \ref{steinberg} and let $B(p)$ denote the Borel of $G$ containing $\br_{\lambda}(\gal{\Q_p})$;
\item at $q$ we take the minimal deformation condition of \S \ref{minimal};
\item at $\ell$ we take ordinary deformations as in \S \ref{ordsection}; to be precise, the character $\bar{\rho}_T$ of Equation (\ref{diagchar}) by construction satisfies $\alpha \circ \bar{\rho}_T |_{I_{\Ql}}= \kbar$ for all $\alpha \in \Delta$, and we choose a lift $\chi_T \colon I_{\Ql} \to T(\mc{O}_{E_{\lambda}})$ satisfying $\alpha \circ \chi_T= \kappa^{r_{\alpha}}$ for positive integers $r_{\alpha} \equiv 1 \pmod {\ell-1}$ with not all $r_{\alpha}$ equal to one another.
\end{itemize}
In defining these local conditions, if necessary we enlarge the coefficient field $\mathbb{F}_{\lambda}$ and work with $\mc{C}_{\mc{O}}^f$ for an appropriate extension $\mc{O}$ of $W(\mathbb{F}_{\lambda})$. Applying Theorem \ref{principalsl2lift}, we obtain a $\rho_{\lambda} \colon \gal{\Q} \to G(\mc{O})$ that is de Rham and unramified outside a finite set of primes containing $p$, $q$, and $\ell$. By Lemma \ref{reductivemono}, the algebraic monodromy group $G_{\rho_{\lambda}}$ is reductive. By Lemma \ref{avoidsl2}, we will be done as long as we can show that $G_{\rho_{\lambda}}$ contains a regular unipotent element of $G$.

To see this, we consider the restriction $\rho_{\lambda}|_{\gal{\Q_p}}$. Since $\bar{r}_{f, \lambda}|_{\gal{\Q_p}}$ is indecomposable, $\bar{r}_{f, \lambda}(I_{\Q_p})$ contains a regular unipotent element of $\mr{GL}_2(\mathbb{F}_{\lambda})$. From the construction (\cite[\S 2.3-2.4]{serre:principalsl2}) of the principal homomorphism $\varphi$, and the elementary regularity criterion \cite[Proposition 5.1.3]{carter:finitelie} for a unipotent element, it follows that $\br_{\lambda}(I_{\Q_p})$ contains a regular unipotent element of $G(\mathbb{F}_{\lambda})$. Of course $\kappa$ is unramified at $p$, so by definition of the Steinberg deformation condition $\alpha \circ \rho_{\lambda}(I_{\Q_p})=1$ for all $\alpha \in \Delta$, and therefore $\rho_{\lambda}|_{I_{\Q_p}}$ is valued in the unipotent radical of $B(p)$. For any element $g \in I_{\Q_p}$ such that $\br_{\lambda}(g)$ is regular unipotent, we deduce that $\rho_{\lambda}(g)$ is regular unipotent in $G(\Frac(\mc{O}))$.\footnote{Note that if we did not restrict to adjoint $G$, the above argument would show that $\rho_{\lambda}(g)$ is a product $z \cdot u$, where $z \in Z_G$ and $u$ is regular unipotent in $G$; but of course then by Jordan decomposition $u \in G_{\rho}$ as well.}
\endproof
\begin{rmk}
\begin{enumerate}
\item The density-one set of the theorem is `explicit' in the sense that, given an $\ell$, we can check whether the argument of the theorem applies to it: we take our chosen modular form and compute its $\ell^{th}$ Fourier coefficient.
\item It is not hard to see that there are infinitely many $f$ for which this argument applies. One might ask whether they can somehow be played off of one another to deduce a version of the theorem for almost all $\ell$ rather than a density one set.
\item We have restricted to the case of adjoint $G$ only for simplicity; as in \S \ref{principalsl2.1}, it suffices to work with a $G$ having a `principal $\mr{GL}_2$,' i.e. for which $\rho^\vee$ makes sense as a co-character.
\end{enumerate}
\end{rmk}
\section{Deformation theory for L-groups}\label{defL}
In this and the following section we carry out the technical modifications needed to treat the case of $\mr{E}_6$ in Theorem \ref{mainintro}. We begin here with a brief discussion of Galois deformation theory for non-connected L-groups. For simplicity we have not sought optimal generality in this discussion, but what we do is more than enough for our purposes. The reader should note that \cite[\S 2]{clozel-harris-taylor} provides a template for this discussion; those authors work with the L-group of an outer form of $\mr{GL}_n \times \mr{GL}_1$.

\subsection{Group theory background}\label{splitLgroup} 
Let $\Psi= (X^{\bullet}, X_{\bullet}, \Phi, \Delta, \Phi^\vee, \Delta^\vee)$ be a based root datum. By a dual group for $\Psi$ we will mean a pinned split reductive group scheme over $\Z$ whose associated root datum is the dual root datum $\Psi^\vee$. To be precise, following \cite[Definition 5.1.1]{conrad:luminy}, this consists of
\begin{enumerate}
\item a reductive group $G^\vee$ over $\Z$;
\item a maximal torus $T^\vee$ of $G^\vee$ equipped with an isomorphism $\iota \colon T^\vee \xrightarrow{\sim} \mr{D}(X_{\bullet})$, where $\mr{D}(\bullet)$ denotes the `functor of characters' (\cite[Appendix B]{conrad:luminy}), and where $\iota$ satisfies the (co)root conditions of \cite[Definition 5.1.1]{conrad:luminy};
\item for each simple coroot $\alpha^\vee \in \Delta^\vee$, a choice $X_{\alpha^\vee}$ of basis of the free rank one $\Z$-module $\mf{g}^\vee_{\alpha^\vee}(\Z)$.
\end{enumerate}


By \cite[Theorem 7.1.9]{conrad:luminy} the outer automorphism group $\Out_{G^\vee}$ is identified with the constant group scheme of automorphisms $\Aut(\Psi^\vee)$ of our based root datum, and the pinning induces a splitting
\begin{equation}\label{splitting}
\Aut(\Psi)= \Aut(\Psi^\vee) \cong \Out_{G^\vee} \into \Aut_{G^\vee}.
\end{equation}
of the canonical projection $\Aut_{G^\vee} \to \Out_{G^\vee}$. In particular, we can define the semi-direct product group scheme $G^\vee \rtimes \Aut(\Psi)$ over $\Z$.

Now let $F$ be a field with separable closure $\overline{F}$ and as always absolute Galois group $\gal{F}= \Gal(\overline{F}/F)$. Let $G$ be a connected reductive group over $F$; over $\overline{F}$, we fix a Borel and maximal torus $T_{\overline{F}} \subset B_{\overline{F}} \subset G_{\overline{F}}$ and define $\Psi$ to be the associated based root datum. In the usual way (\cite[\S 1]{borel:L}) we obtain a canonical homomorphism $\mu_G \colon \gal{F} \to \Aut(\Psi)$, depending only on the class of inner forms to which $G$ belongs. Consider a dual group $G^\vee$ (the other data being implicit) for $\Psi$. We then define the (split form of the) L-group of $G$ by combining $\mu_G$ with the splitting (\ref{splitting}): ${}^L G= G^\vee \rtimes \gal{F}$. Thus we have defined a group scheme over $\Z$ whose base change to an algebraically closed field is the `usual' L-group of $G$.

We now recall (\cite[\S 2]{gross:principalsl2}) that the principal $\mr{SL}_2$ extends to a homomorphism
\[
\varphi \colon \mr{SL}_2 \times \gal{F} \to {}^L G.
\]
To be precise, as in Lemma \ref{integralsl2} let $\ell$ be a prime greater than or equal to $h^\vee$, the Coxeter number of $G^\vee$, so that there is (using the pinning) a principal homomorphism $\varphi \colon \mr{SL}_2 \to G^\vee$ defined over $\Z_{(\ell)}$. Since $\gal{F}$ permutes both the elements of $\Delta^\vee$ and the corresponding set of positive coroots, it is easily seen to preserve the $\mf{sl}_2$-triple $(X, H, Y)$ in $\mf{g}^\vee$ (see \S \ref{principalsl2.1}). Thus $\gal{F}$ preserves $\varphi$, which consequently extends to the desired $\varphi \colon \mr{SL}_2 \times \gal{F} \to {}^L G$.
\subsection{Deformation theory}
We now consider an L-group ${}^L G$ over $\Z$ as in \S \ref{splitLgroup}. For simplicity, from now on we assume $G$ is simply-connected, so the dual group $G^\vee$ is an adjoint Chevalley group; this will avoid the need to `fix the determinant' in what follows, and it suffices for our application. Moreover, $\Aut(\Psi)$ is now finite, and the homomorphism $\mu_G$ factors through a faithful homomorphism $\Gal(\widetilde{F}/F) \to \Aut(\Psi)$ for some finite extension $\tF/F$. We can and do replace ${}^L G$ with the finite form ${}^L G= G^\vee \rtimes \Gal(\tF/F)$ in all that follows. For simplicity, we will moreover assume that $G$ is chosen so that $\tF$ is totally imaginary.

Let $\Sigma$ denote a finite set of finite places of $F$; we will always assume that $\Sigma$ contains all places above $\ell$ and that all elements of $\Sigma$ split in $\tF/F$. Let $\tF_{\Sigma}$ be the maximal extension of $\tF$ unramified outside (places above) $\Sigma$. The extension $\tF_{\Sigma}/F$ is Galois, and we have an inclusion of Galois groups
\[
\gal{\tF, \Sigma} = \Gal(\tF_{\Sigma}/\tF) \subset \Gal(\tF_{\Sigma}/F)= \gal{\Sigma}.
\]
Here we have introduced the new notation $\gal{\Sigma}$ for the latter group; it is \textit{not} the same as the group $\gal{F, \Sigma}$, since $\tF/F$ may be ramified outside $\Sigma$. For each $v \in \Sigma$, choose a place $\tv$ of $\tF$ above $v$; we write $\widetilde{\Sigma}$ for the collection of all such $\tv$. We fix one member of the $\gal{\tF, \Sigma}$-conjugacy class of homomorphisms $\gal{\tF_{\tv}} \to \gal{\tF, \Sigma}$, and we thereby also obtain a homomorphism $\gal{F_v} \to \gal{\tF, \Sigma}$ (whose $\gal{\tF, \Sigma}$-conjugacy class depends on the choice of $\tv \vert v$). 

Let $\br \colon \gal{\Sigma} \to {}^L G(k)$ be an `L-homomorphism,' i.e. a continuous homomorphism such that the diagram
\[
\xymatrix{
\gal{\Sigma} \ar[rr]^{\br} \ar[rd] && {}^L G(k) \ar[ld] \\
& \Gal(\tF/F) &
}
\]
commutes. In particular, the restriction of $\br$ to $\gal{\tF, \Sigma}$ factors through $G^\vee(k)$. In parallel to Definition \ref{liftfunctor}, we now define the relevant functors of lifts and deformations of $\br$:
\begin{defn}
\begin{itemize}
\item Let $\Lift_{\br} \colon \mc{C}^f_{\mc{O}} \to \Sets$ be the pro-representable functor whose $R$-points is the set of all lifts of $\br$ to a continuous homomorphism (automatically an L-homomorphism) $\gal{\Sigma} \to {}^L G(R)$. We denote the representing object by $R_{\br}^{\square}$.
\item We say two lifts $\rho_1$ and $\rho_2$ are strictly equivalent if they are conjugate by an element of 
\[
\widehat{G^\vee}(R)= \ker \left( G^\vee(R) \to G^\vee(k) \right).
\]
\item Denote by $\Def_{\br} \colon \CO^f \to \Sets$ the functor of strict equivalence classes in $\Lift_{\br}$.
\item A collection of local deformation problems for $\br$ is, for each $v \in \Sigma$, a representable sub-functor $\Lift^{\mc{P}_v}_{\br|_{\gal{\tF_{\tv}}}}$ of $\Lift_{\br|_{\gal{\tF_{\tv}}}}$ that is closed under strict equivalence. Note that here $\br |_{\gal{\tF_{\tv}}}$ is a $G^\vee$-valued homomorphism, so the local deformation conditions in question are no different from those considered in \S \ref{review}-\S \ref{principalSL2}.
\item Writing $\mc{P}= \{\mc{P}_v\}_{v \in \Sigma}$ for a collection of local deformation conditions as in the previous item, we define the global functor $\Lift_{\br}^{\mc{P}}$ by taking the sub-functor of $\Lift_{\br}$ of lifts whose restrictions to each $\gal{\tF_{\tv}}$ lie in $\Lift^{\mc{P}_v}_{\br |_{\gal{\tF_{\tv}}}}$; the quotient of $\Lift_{\br}^{\mc{P}}$ by strict equivalence defines $\Def_{\br}^{\mc{P}}$.
\item Recall the subspace $L_{\tv}^{\square} \subset Z^1(\gal{\tF_{\tv}}, \br(\mf{g}^\vee))$ associated to $\mc{P}_v$, with image $L_{\tv}$ in $H^1(\gal{\tF_{\tv}}, \br(\fg^\vee))$. For $i \geq 0$, define $L_{\tv}^{\square, i}$ and the complex $C^i_{\mc{P}}(\gal{\Sigma}, \br(\fg^\vee))$ exactly as in \S \ref{deflocglob} (see Equation (\ref{globtangent}) and following).
\end{itemize}
\end{defn}
We impose the usual requirement that $\ell$ be very good for $G$ (\S \ref{review}); if $G$ has a simple factor of type $\mr{D}_4$, also exclude $\ell=3$. It is now important that $\ell$ not divide $[\tF:F]$, but in fact that follows from the other assumptions, since a simple Dynkin diagram has automorphism group with order divisible at most by the primes 2 and 3 (with 3 only occurring in type $\mr{D}_4$). The next proposition summarizes the basic facts about Galois deformation theory in this setting; this is essentially the same as \cite[2.3.3-2.3.5]{clozel-harris-taylor}. Since we did not give the corresponding proofs in \S \ref{review}, we linger over them here.
\begin{prop}\label{defLprop}
Let $\br \colon \gal{\Sigma} \to {}^L G(k)$ be a continuous L-homomorphism. For each $v \in \Sigma$, fix a local deformation condition $\mc{P}_v$. Assume that the centralizer of $\br|_{\gal{\tF, \Sigma}}$ in $\fg^\vee$ is trivial.
\begin{enumerate}
\item There is a canonical isomorphism $\Def_{\br}^{\mc{P}}(k[\epsilon]) \cong H^1_{\mc{P}}(\gal{\Sigma}, \br(\fg^\vee))$.
\item $\Def_{\br}^{\mc{P}}$ is (pro-)representable. We denote by $R_{\br}^{\mc{P}}$ the representing object.
\item The analogue of Wiles's formula holds:
\begin{align}\label{chi0}
&h^1_{\mc{P}}(\gal{\Sigma}, \brgd)-h^1_{\mc{P}^\perp}(\gal{\Sigma}, \brgd(1))  \\
&= h^0(\gal{\Sigma}, \brgd)- h^0(\gal{\Sigma}, \brgd(1))- \sum_{v \vert \infty} h^0(\gal{F_v}, \brgd)+ \sum_{v \in \Sigma} \left(\dim_k L_{\tv} - h^0(\gal{\tF_{\tv}}, \brgd) \right). \nonumber
\end{align}
\item Moreover assume the local conditions $\Lift^{\mc{P}_v}_{\br |_{\gal{\tF_{\tv}}}}$ are liftable. Then $R_{\br}^{\mc{P}}$ is isomorphic to a quotient of a power series ring over $\mc{O}$ in $\dim_k H^1_{\mc{P}}(\gal{\Sigma}, \br(\fg^\vee))$ variables by an ideal that can be generated by at most $\dim_k H^1_{\mc{P}^\perp}(\gal{\Sigma}, \br(\fg^\vee)(1))$ elements.
\end{enumerate}
\end{prop}
\begin{rmk}
Just as in Equation (\ref{wilesformula}) of \S \ref{formalglobal}, in the application our local calculations and assumptions on the image of $\br$ will imply that the right-hand-side of Equation (\ref{chi0}) is zero, hence that $h^1_{\mc{P}}(\gal{\Sigma}, \brgd)= h^1_{\mc{P}^\perp}(\gal{\Sigma}, \brgd(1))$, as needed for Ramakrishna's method.
\end{rmk}
\proof
The description of the tangent space $\Lift_{\br}(k[\epsilon])$ follows from the usual argument, noting that when we write $\br(g)= (\br_0(g), g) \in G^\vee(k) \times \Gal(\tF/F)$, $\br_0$ is a cocycle in $Z^1(\gal{\Sigma}, G^\vee(k))$, where $\gal{\Sigma}$ is regarded as acting on $G^\vee$ via the outer $\Gal(\tF/F)$-action. Item (1) follows easily.

Representability of $\Def_{\br}$ follows, using Schlessinger's criteria, from our assumption on the centralizer of $\br|_{\gal{\tF, \Sigma}}$. We check this somewhat more generally, when $\fg^\vee$ is not necessarily adjoint, but the invariants $(\fg^\vee)^{\br(\gal{\tF, \Sigma})}$ are assumed to equal the center $\mf{z}(\fg^\vee)$. Since $\widehat{G^\vee}$ is formally smooth, and $h^1(\gal{\Sigma}, \brgd)$ is finite, the only thing that really requires checking is the injectivity statement in Schlessinger's condition (H4), namely that for all small extensions $A \to B$ in $\CO^f$, the map
\[
\Def_{\br}(A \times_B A) \to \Def_{\br}(A) \times_{\Def_{\br}(B)} \Def_{\br}(A)
\]
is injective. Consider two elements $\rho$ and $r$ of $\Lift_{\br}(A \times_B A)$, and let $g$ and $g'$ be elements of $\widehat{G^\vee}(A)$ attesting to their equivalence in the right-hand-side. Pushing down to $B$, which we denote by a subscript $B$, the element $g_B^{-1} g'_B$ of $\widehat{G^\vee}(B)$ commutes with the image of $\rho_B \in \Lift_{\br}(B)$; in particular, it commutes with the image of $\rho_B|_{\gal{\tF, \Sigma}}$. We claim that under our hypotheses on $\br$, the centralizer $Z_{G^\vee(B)}(\rho_B|_{\gal{\tF, \Sigma}})$ equals $Z_{G^\vee}(B)$ for all lifts $\rho_B$ of $\br$. To see this, we argue by induction on the length of $B$: if $B \to B/I$ is a small extension with $I$ a 1-dimensional $k$-vector space, then the induction hypothesis (for $B/I$) and smoothness of $Z_{G^\vee}$ imply that any element $z$ of $Z_{G^\vee(B)}(\rho_B|_{\gal{\tF, \Sigma}})$ has the form $z_B \cdot \exp(X)$ for some $z_B \in Z_{G^\vee}(B)$ and
\[
X \in \left(\mf{g}^\vee \otimes_k I \right)^{\br(\gal{\tF, \Sigma})}= \mf{z}(\fg^\vee) \otimes_k I.
\]
Thus $z \in Z_{G^\vee}(B)$. Lifting $g_B^{-1}g'_B$ to some $\tilde{z} \in Z_{G^\vee}(A)$, and replacing $g$ by $g \tilde{z}^{-1}$ we may assume $g$ and $g'$ map to a common element of $B$, proving that $\rho$ and $r$ are equivalent as elements of $\Def_{\br}(A \times_B A)$.

We proceed to part (3). By construction of the complex $C^{\bullet}_{\mc{P}}(\gal{\Sigma}, \br(\fg^\vee))$, we have a long exact sequence
\begin{equation}\label{selmerles}
\xymatrix{
0 \ar[r] & H^1_{\mc{P}}(\gal{\Sigma}, \brgd) \ar[r] & H^1(\gal{\Sigma}, \brgd) \ar[r] & \bigoplus_{v \in \Sigma} H^1(\gal{\tF_{\tv}}, \brgd)/L_{\tv} \ar[r] &\\
\ar[r] & H^2_{\mc{P}}(\gal{\Sigma}, \brgd) \ar[r] & H^2(\gal{\Sigma}, \brgd) \ar[r] & \bigoplus_{v \in \Sigma} H^2(\gal{\tF_{\tv}}, \brgd) \ar[r] &\\
\ar[r] & H^3_{\mc{P}}(\gal{\Sigma}, \brgd) \ar[r] & 0, & &
}
\end{equation}
where the final term is zero because $H^3(\gal{\Sigma}, \brgd)= H^3(\gal{\tF, \Sigma}, \brgd)^{\Gal(\tF/F)}=0$ (recall that $\tF$ is totally imaginary; indeed, the same argument, using that $\ell$ is coprime to $[\tF:F]$, implies that $H^i_{\mc{P}}(\gal{\Sigma}, \brgd)$ vanishes for all $i>3$). The 9-term Poitou-Tate exact sequence of global duality also yields a long exact sequence
\begin{equation}\label{poitoutate}
\xymatrix{
0 \ar[r] & H^1_{\mc{P}}(\gal{\Sigma}, \brgd) \ar[r] & H^1(\gal{\Sigma}, \brgd) \ar[r] & \bigoplus_{v \in \Sigma} H^1(\gal{\tF_{\tv}}, \brgd)/L_{\tv} \ar[r] &\\
\ar[r] & H^1_{\mc{P}^\perp}(\gal{\Sigma}, \brgd(1))^\vee \ar[r] & H^2(\gal{\Sigma}, \brgd) \ar[r] & \bigoplus_{v \in \Sigma} H^2(\gal{\tF_{\tv}}, \brgd) \ar[r] &\\
\ar[r] & H^0(\gal{\Sigma}, \brgd(1))^\vee \ar[r] & 0. & &
}
\end{equation}
To be precise, global duality gives such a sequence with $\gal{\tF, \Sigma}$ in place of $\gal{\Sigma}$ (and all places of $\tF$ above $\Sigma$ contributing to the local terms); the above sequence is the result of taking $\Gal(\tF/F)$-invariants. Comparing this and the previous 7-term sequence, we conclude that
\begin{itemize}
\item $\dim_k H^3_{\mc{P}}(\gal{\Sigma}, \brgd)= \dim_k H^0(\gal{\Sigma}, \brgd(1))$; and
\item $\dim_k H^2_{\mc{P}}(\gal{\Sigma}, \brgd)= \dim_k H^1_{\mc{P}^\perp}(\gal{\Sigma}, \brgd(1))$.
\end{itemize} 
Now, the definition of the complex $C^\bullet_{\mc{P}}(\gal{\Sigma}, \brgd)$ implies its Euler-characteristic is equal to 
\[
\chi(\Gamma_{\Sigma}, \brgd)-\sum_{v \in \Sigma} \chi(\gal{\tF_{\tv}}, \brgd)+ \sum_{v \in \Sigma} \left( h^0(\gal{\tF_{\tv}}, \brgd)- \dim_k L_{\tv} \right).
\]
A minor variant of the global Euler-characteristic formula (demonstrated in \cite[Lemma 2.3.3]{clozel-harris-taylor}), combined with the local Euler-characteristic formula, yields the formula
\begin{align}\label{chi1}
\chi(C^\bullet_{\mc{P}}(\gal{\Sigma}, \brgd))&= \sum_{v \vert \infty} h^0(\gal{F_v}, \brgd) - [F:\Q] \dim_k \fg^\vee \\
&+ \sum_{v \vert \ell} [F_v:\Ql]\dim_k \fg^\vee+ \sum_{v \in \Sigma} \left( h^0(\gal{\tF_{\tv}}, \brgd)- \dim_k L_{\tv} \right) \nonumber \\
&= \sum_{v \vert \infty} h^0(\gal{F_v}, \brgd) + \sum_{v \in \Sigma} \left( h^0(\gal{\tF_{\tv}}, \brgd)- \dim_k L_{\tv} \right). \nonumber
\end{align}
On the other hand, we have just seen that
\begin{equation}\label{chi2}
\chi(C^\bullet_{\mc{P}}(\gal{\Sigma}, \brgd))= h^0(\gal{\Sigma}, \brgd)-h^1_{\mc{P}}(\gal{\Sigma}, \brgd)+h^1_{\mc{P}^\perp}(\gal{\Sigma}, \brgd(1))- h^0(\gal{\Sigma}, \brgd(1)).
\end{equation}
Combining Equations (\ref{chi1}) and (\ref{chi2}), we obtain the formula (\ref{chi0}) stated in the proposition.

Finally, under the assumption that each $\mc{P}_v$ is liftable, for any small extension $R \to R/I$, and any lift $\rho \in \Lift_{\br}^{\mc{P}}(R/I)$, we can construct an obstruction class $\mr{obs}_{\rho, R, \mc{P}} \in H^2_{\mc{P}}(\gal{\Sigma}, \brgd) \otimes_k I$, defined as in the statement of \cite[Lemma 2.2.11]{clozel-harris-taylor}.\footnote{In their notation, take $S= \Sigma$ and $T= \emptyset$.} This class vanishes if and only if $\rho$ admits a lift to $\Lift_{\br}^{\mc{P}}(R)$. Then a classic deformation theory argument (see \cite[Proposition 2]{mazur:deforming}) shows that $R_{\br}^{\mc{P}}$ is a quotient of a power series ring over $\mc{O}$ in $\dim_k H^1_{\mc{P}}(\gal{\Sigma}, \brgd)$ variables by an ideal generated by at most $\dim_k H^2_{\mc{P}}(\gal{\Sigma}, \brgd)$ elements; the last part of the proposition follows from the equality (established in the course proving (3)) $h^2_{\mc{P}}(\gal{\Sigma}, \brgd)= h^1_{\mc{P}^\perp}(\gal{\Sigma}, \brgd(1))$.
\endproof
\begin{rmk}
We note a basic compatibility in the construction of these Selmer groups. For any finite $\gal{\Sigma}$-module $M$, let $\{L_{\tv}\}_{v \in \Sigma}$ be a collection of sub-modules $L_{\tv} \subset H^1(\gal{\tF_{\tv}}, M)$. For any $w \not \in \Sigma$ split in $\tF/F$, with a specified place $\tilde{w}$ of $\tF$ above $w$, the inflation map 
\[
H^1_{\{L_{\tv}\}}(\gal{\Sigma}, M) \to H^1_{\{L_{\tv}\} \cup L_{\tilde{w}}^{\mr{unr}}}(\gal{\Sigma \cup w}, M)
\]
is an isomorphism: for surjectivity, note that a cocycle $\phi$ such that $\phi|_{I_{\tF_{\tilde{w}}}}=0$ also vanishes on $I_{\tF_{\tilde{w}'}}$ for all $\tilde{w}' \vert w$,\footnote{The groups $I_{\tF_{\tilde{w}}}$ and $I_{\tF_{\tilde{w}'}}$ are conjugate in $\gal{\Sigma \cup w}$, and the cocycle relation implies $\phi(ghg^{-1})=0$ whenever $\phi(h)=0$ and $h$ acts trivially on $M$.} so factors through $\gal{\Sigma}$.
\end{rmk}
\section{The case $-1 \not \in W_G$}\label{E6section}
A minor variant of the argument of \S\S \ref{formalglobal} and \ref{principalSL2} will treat the case of type $\mr{E}_6$ and consequently complete the proof of Theorem \ref{mainintro}. We carry out this variant in the present section, by establishing a version of Ramakrishna's lifting theorem for groups whose Weyl group does not contain $-1$.  As in \S \ref{defL}, we restrict for simplicity to the case in which $G^\vee$ is an adjoint group; thus, we are really concerned only with the groups $\mr{PGL}_n$, $\mr{PSO}_{4n+2}$, and the adjoint form of $\mr{E}_6$. Note that in these cases the automorphism group of the Dynkin diagram is always $\Z/2\Z$.
\subsection{Constructing ${}^L G$}\label{outerodd}
We retain the general setup of \S \ref{defL}. Thus $\Psi$ is a based root datum for which we construct a dual group  $G^\vee$ (in the sense of \S \ref{splitLgroup}). We fix a number field $F$, which from now on will be assumed totally real, and to a connected reductive group $G$ over $F$ with (absolute) root datum $\Psi$ we can associate a $\Z$-form of the L-group ${}^L G$. Our first aim is to choose this $G$, in the case $-1 \not \in W_G$, so that ${}^L G$ admits `odd' homomorphisms from $\gal{F}$, allowing us to work out a (minor) variant of Ramakrishna's method for ${}^L G$. Recall from \S \ref{archcond} that (now letting $\rho$ be a half sum of positive roots of $G$, hence a co-character of $G^\vee$) $\Ad(\rho(-1))$ is no longer a split Cartan involution of $\fg^\vee$, so composing an odd two-dimensional representation $\gal{F} \to \mr{PGL}_2(k)$ with the principal homomorphism $\mr{PGL}_2 \to G^\vee$ will no longer yield a homomorphism $\gal{F} \to G^\vee(k)$ satisfying item (\ref{infiniteplaces}) of \S \ref{formalglobal}. We now explain how to rectify this.

Let $\tF/F$ be a quadratic totally imaginary extension of the totally real field $F$. The choice of $\tF$ induces a canonical non-trivial element $\delta_{\tF/F}$ of $\Hom(\gal{F}, \Z/2\Z)=\Hom(\gal{F}, \Aut(\Psi))$ (note that $\Z/2\Z$ has no non-trivial automorphism, so we are justified in writing `$=$'). We choose $G$ to be any form over $F$ of the root datum $\Psi$ so that the associated homomorphism $\mu_G \colon \gal{F} \to \Aut(\Psi)$ is equal to $\delta_{\tF/F}$. To be precise, if $G_0$ is a pinned split form of $\Psi$ over $F$, giving a base-point in $H^1(\gal{F}, \Aut(G_{0, \overline{F}}))$, we can take $G$ to be any form whose cohomology class lifts $\delta_{\tF/F}$ under the homomorphism
\[
H^1(\gal{F}, \Aut(G_{0, \overline{F}})) \to H^1(\gal{F}, \Out(G_{0, \overline{F}})) = H^1(\gal{F}, \Aut(\Psi)).
\]
As in \S \ref{splitLgroup} we obtain the principal homomorphism
\[
\varphi \colon \mr{PGL}_2 \times \gal{F} \to {}^L G.
\]
The next lemma, which is essentially \cite[Proposition 7.2]{gross:principalsl2}, establishes the necessary `oddness' for the representations $\gal{\Sigma} \to {}^L G(k)$ that we will consider in Theorem \ref{E6theorem}.
\begin{lemma}\label{oddout}
For any infinite place $v \vert \infty$ of $F$, let $c_v \in \gal{F}$ denote a choice of complex conjugation, and let $\theta$ be the element
\[
\theta=  \varphi \left( \begin{pmatrix} -1 & 0 \\ 0 & 1 \end{pmatrix} \times c_v \right)
\]
of ${}^L G$. Then $\Ad(\theta)$ is a split Cartan involution of $\fg^\vee$, i.e. $\dim_k(\fg^\vee)^{\Ad(\theta)=1} = \dim_k \fn$.
\end{lemma}
\proof
For our choice of $G$, the homomorphism $\mu_G$ factors through $\delta_{\tF/F} \colon \Gal(\tF/F) \to \Aut(\Psi)$, with any complex conjugation $c_v$ mapping to the non-trivial element of $\Aut(\Psi)$ (the opposition involution). By \cite[Proposition 7.2]{gross:principalsl2}, which continues to hold under the assumption $\ell \geq 2h^\vee -1$ (see Lemma \ref{kostantmodl} and the proof of \cite[Proposition 7.2]{gross:principalsl2}), the trace of $\Ad(\theta)$ is $-\rk(\fg^\vee)$. Since $\theta^2 =1$, it follows easily that
\[
\dim_k (\fg^\vee)^{\Ad \theta =1}= \frac{\dim_k \fg^\vee - \rk \fg^\vee}{2}= \dim_k \mf{n}.
\]
\endproof
\subsection{The lifting theorems}\label{formalout}
We continue with an $\tF/F$ and ${}^L G$ as in \S \ref{outerodd}. Let $\Sigma$ be a finite set of finite places of $F$, containing all places above $\ell$, such that all members of $\Sigma$ split in $\tF/F$. For each $v \in \Sigma$, we fix extensions $\tv$ of $v$ to $\tF$ and fix embeddings $\gal{\tF_{\tv}} \into \gal{\tF, \Sigma}$. Assume we are given a continuous L-homomorphism 
\[
\br \colon \gal{\Sigma} \to {}^L G(k)
\]
(in particular, $\br(\gal{\tF, \Sigma}) \subset G^\vee(k)$) such that the centralizer in $\mf{g}^\vee$ of $\br|_{\gal{\tF, \Sigma}}$ is trivial, and that moreover satisfies the following properties:
\begin{enumerate}
\item\label{dualinvsout} $h^0(\gal{\Sigma}, \brgd)=h^0(\gal{\Sigma}, \brgd(1))=0$.
\item\label{locdefcondout} For all $v \in \Sigma$, there is a liftable local deformation condition $\mc{P}_v$ for $\br|_{\gal{\tF_{\tv}}}$ satisfying
\[
\dim L_{\tv}= \begin{cases}
\text{$h^0(\gal{\tF_{\tv}}, \brgd)$ if $v \nmid \ell$;}\\
\text{$h^0(\gal{\tF_{\tv}}, \brgd)+[F_v:\Ql]\dim(\fn)$ if $v \vert \ell$.}
\end{cases}
\]
\item\label{splitcartanout} For all $v \vert \infty$, $\Ad (\br(c_v))$ is a split Cartan involution of $\fg^\vee$, i.e. $h^0(\gal{F_v}, \brgd)= \dim(\fn)$.
\item\label{cohvanishout} Let $K= \tF(\brgd, \mu_{\ell})$. Note that $\tF_{\Sigma}$ contains $K$. Then $H^1(\Gal(K/F), \brgd)=0$ and $H^1(\Gal(K/F), \brgd(1))=0$.
\item\label{lindisjointout} For any pair of non-zero Selmer classes $\phi \in H^1_{\mc{P}^\perp}(\gal{\Sigma}, \brgd(1))$ and $\psi \in H^1_{\mc{P}}(\gal{\Sigma}, \brgd)$, the restrictions of $\phi$ and $\psi$ to $\Gal(\tF_{\Sigma}/K)$ are homomorphisms with fixed fields $K_{\phi}$ and $K_{\psi}$ that are disjoint over $K$. (From now on we denote $\Gal(\tF_{\Sigma}/K)$ by $\gal{K, \Sigma}$; this notation is consistent with the notation $\gal{\tF, \Sigma}$ since $K/\tF$ is ramified only at places above $\Sigma$.)
\item\label{avoidout} Consider any $\phi$ and $\psi$ as in the hypothesis of item (\ref{lindisjointout}) (we do not require the conclusion to hold). Then there is an element $\sigma \in \gal{\tF, \Sigma}$ such that $\br(\sigma)$ is a regular semi-simple element of $G^\vee$, the connected component of whose centralizer we denote $T^\vee$, and such that there exists a root $\alpha^\vee \in \Phi(G^\vee, T^\vee)$ satisfying
\begin{enumerate}
\item $\kbar(\sigma)= \alpha^\vee \circ \br(\sigma)$; 
\item $k[\psi(\gal{K, \Sigma})]$ has an element with non-zero $\mf{l}_{\alpha^\vee}$ component;\footnote{Here recall that we write $\mf{t}^\vee= \mf{l}_{\alpha^\vee} \oplus \mf{t}^\vee_{\alpha^\vee}$ for the decomposition of $\mf{t}^\vee$ into the span $\mf{l}_{\alpha^\vee}$ of the $\alpha$-coroot vector for $G^\vee$ and $\mf{t}^\vee_{\alpha^\vee}= \ker(\alpha)$.} and
\item $k[\phi(\gal{K, \Sigma})]$ has an element with non-zero $\mf{g}^\vee_{-\alpha^\vee}$ component.
\end{enumerate}
\end{enumerate}
\begin{prop}\label{formalliftouter}
Under assumptions (\ref{dualinvsout})-(\ref{avoidout}) above, there exists a finite set of primes $Q$ of $F$, disjoint from $\Sigma$ and split in $\tF/F$, and a lift
\[
\xymatrix{
& {}^L G(\mc{O}) \ar[d] \\
\gal{\Sigma \cup Q} \ar[r]_{\br} \ar@{-->}[ur]^{\rho} & {}^L G(k)
}
\]
such that $\rho$ is type $\mc{P}_v$ at all $\tv \in \widetilde{\Sigma}$ and of Ramakrishna type at each $\tv \in \widetilde{Q}$, where $\widetilde{Q}$ consists of, for each $v$ in $Q$, a specified extension $\tv$ of $v$ to $\tF$. 
\end{prop}
\proof
With the following modifications, the proof of Proposition \ref{formalliftingtheorem} applies verbatim:
\begin{itemize}
\item Replace $\gal{F, \Sigma}$ with $\gal{\Sigma}$ (likewise for $\Sigma \cup w$).
\item In Lemma \ref{findingw}, require $w \not \in \Sigma$ to be split in $\tF/F$; note that now we assume the existence of a $\sigma \in \Gal(\tF(\br)/\tF)$ satisfying the conclusion of item (\ref{avoidout}), and we find the desired split primes of $\tF/F$ by applying the \v{C}ebotarev density theorem to the Galois extension $\tF(\br)K_{\phi}K_{\psi}/\tF$ (recall that the primes of $\tF$ that are split over $F$ have density one in $\tF$).
\end{itemize}
\endproof
There is no difficulty now in deducing an analogue for ${}^L G$ of Theorem \ref{maximaltheorem}:
\begin{thm}\label{maximalimageout}
Let $F$ be a totally real field with $[F(\mu_{\ell}):F]= \ell -1$, and let $\tF/F$, $\Sigma$, and ${}^L G$ be as in \S \ref{outerodd}. Suppose $\br \colon \Gamma_{\Sigma} \to {}^L G(k)$ is a continuous representation satisfying the following conditions:
\begin{enumerate}
\item There is a subfield $k' \subset k$ such that $\br(\gal{\tF, \Sigma})$ contains $\im \left( G^{\vee, \mr{sc}}(k') \to G^\vee(k') \right)$.
\item $\ell -1$ is greater than the maximum of $8\cdot \# Z_{G^{\vee, \mr{sc}}}$ and
\[
\begin{cases}
\text{$(h-1)\# Z_{G^{\vee, \mr{sc}}}$ if $\# Z_{G^{\vee, \mr{sc}}}$ is even; or} \\
\text{$(2h-2) \# Z_{G^{\vee, \mr{sc}}}$ if $\# Z_{G^{\vee, \mr{sc}}}$ is odd.}
\end{cases}
\]
\item $\br$ is odd, i.e. for all complex conjugations $c_v$, $\Ad(\br(c_v))$ is a split Cartan involution of $\mf{g}^\vee$.
\item For all places $v \in \Sigma$ not dividing $\ell \cdot \infty$, $\br|_{\gal{\tF_{\tv}}}$ satisfies a liftable local deformation condition $\mc{P}_v$ with tangent space of dimension $h^0(\gal{F_v}, \br(\fg^\vee))$ (eg, the conditions of \S \ref{steinberg} or \S \ref{minimal}).
\item For all places $v \vert \ell$, $\br|_{\gal{\tF_{\tv}}}$ is ordinary in the sense of \S \ref{ordsection}, satisfying the conditions (REG) and (REG*).
\end{enumerate}
Then there exists a finite set of primes $Q$, disjoint from $\Sigma$ and split in $\tF/F$, and a lift 
\[
\xymatrix{
& {}^L G(\mc{O}) \ar[d] \\
\gal{\Sigma \cup Q} \ar[r]_{\br} \ar@{-->}[ur]^{\rho} & {}^L G(k) 
}
\]
such that $\rho$ is type $\mc{P}_v$ at all $\tv \in \widetilde{\Sigma}$ (taking $\mc{P}_v$ to be an appropriate ordinary condition at $v \vert \ell$) and of Ramakrishna type at all $\tv \in \widetilde{Q}$ (again, having fixed an extension $\tv$ of $v$ in $Q$ to $\tF$). In particular $\br$ admits a characteristic zero lift that is geometric in the sense of Fontaine-Mazur.
\end{thm}
\proof
Proceed as in the proof of Theorem \ref{maximaltheorem}, but now using Proposition \ref{formalliftouter}. We leave the details to the reader.
\endproof
Next we deduce an analogue of Theorem \ref{principalsl2lift}:
\begin{thm}\label{principalsl2liftout}
Assume now that $G$ is of type $\mr{E}_6$, and let $\ell$ be a rational prime greater than $4h^\vee-1=47$. Let $F$ be a totally real field for which $[F(\zeta_{\ell}):F]= \ell -1$, and let $\bar{r} \colon \gal{F} \to \mr{GL}_2(k)$ be a continuous representation unramified outside a finite set $\Sigma$ of finite places, which we assume to contain all places above $\ell$. Assume that $\bar{r}$ moreover satisfies the following:
\begin{enumerate}
\item For some subfield $k' \subset k$, 
\[
\mr{SL}_2(k') \subset \bar{r}(\gal{F}) \subset k^\times \cdot \mr{GL}_2(k');
\]
\item $\bar{r}$ is odd;
\item for each $v \vert \ell$, $\bar{r}|_{\gal{F_v}}$ is ordinary, satisfying
\[
\bar{r}|_{I_{F_v}} \sim \begin{pmatrix}
\chi_{1, v} & * \\
0 & \chi_{2, v} \\
\end{pmatrix},
\]
where $\left(\chi_{1,v}/\chi_{2,v}\right)|_{I_{F_v}}= {\kbar}^{r_v}$ for an integer $r_v \geq 2$ such that $\ell>r_v(h^\vee-1)+1$;
\end{enumerate}
We then choose a quadratic totally imaginary extension $\tF/F$ with the following properties:
\begin{itemize}
\item All elements of $\Sigma$ split in $\tF/F$.
\item $\tF$ is linearly disjoint from $F(\bar{r}, \zeta_{\ell})$ over $F$.
\end{itemize}
Then using $\tF/F$ we can define the L-group ${}^L G$ over $\Z$ as in \S \ref{outerodd} and consider the composite homomorphism
\[
\xymatrix{
\gal{\Sigma} \ar[r]^-{\bar{r}} \ar@/^2 pc/[rr]^{\br} & \mr{PGL}_2(k) \times \Gal(\tF/F) \ar[r]^-{\varphi} & {}^L G(k).
} 
\]
We additionally assume that for all primes $v \in \Sigma$ not dividing $\ell$, we can find liftable local deformation conditions $\mc{P}_v$ for $\br|_{\gal{\tF_{\tv}}}$ such that $\dim L_{\tv}= h^0(\gal{F_v}, \brgd)$.

Then there exists a finite set of places $Q$ disjoint from $\Sigma$ and a lift
\[
\xymatrix{
& {}^L G(\mc{O}) \ar[d] \\
\gal{\Sigma \cup Q} \ar@{-->}[ur]^{\rho} \ar[r]^{\br} & {}^L G(k)
}
\]
such that $\rho|_{\gal{\tF_{\tv}}}$ is of type $\mc{P}_v$ for all $v \in \Sigma$, and (having specified a place $\tilde{w} \vert w$ of $\tF$ for all $w \in Q$) $\rho|_{\gal{\tF_{\tilde{w}}}}$ is of Ramakrishna type for all $w \in Q$.
\end{thm}
\begin{rmk}
As with Theorem \ref{principalsl2lift}, the argument will apply to simple types $\mr{D}_{2n+1}$ and $\mr{A}_n$ once Lemma \ref{magmacalc} is established in those cases.
\end{rmk}
\proof
First, it is clear that we can find such an extension $\tF/F$: it suffices to choose a quadratic imaginary field $\Q(\sqrt{-D})$ in which all (rational) primes below $\Sigma$ are split, and which ramifies at some prime that is unramified in $F(\br, \zeta_{\ell})$; existence of such a $D$ follows from the Chinese Remainder Theorem.

We now rapidly verify the six conditions of the axiomatized lifting theorem, as enumerated at the start of \S \ref{formalout}. That the centralizer in $\fg^\vee$ of $\br|_{\gal{\tF, \Sigma}}$ is trivial follows as in Theorem \ref{principalsl2lift}, since $\tF$ is linearly disjoint from $F(\bar{r})$. Likewise condition (\ref{dualinvsout}) follows as before, it even sufficing to consider $\gal{\tF, \Sigma}$-invariants. Condition (\ref{locdefcondout}) is satisfied by assumption, and by taking an appropriate ordinary deformation condition at $v \vert \ell$. Oddness of $\bar{r}$ and Lemma \ref{oddout} together imply condition (\ref{splitcartanout}). 

The argument of Theorem \ref{principalsl2lift} also implies the cohomological vanishing statements of condition (\ref{cohvanishout}): to be precise, the argument there directly applies to the cohomology of $\Gal(K/\tF)$, but since $\ell$ is coprime to $[\tF:F]=2$ the slight strengthening here also holds. Condition (\ref{lindisjointout}) is also the identical argument (using the element $\sigma$ to be constructed in the verification of condition (\ref{avoidout})). For condition (\ref{avoidout}) we construct an element $\sigma \in \Gal(F(\ad^0(\bar{r}), \zeta_{\ell})/L)$, where $L$ is as before the intersection $F(\ad^0(\bar{r})) \cap F(\zeta_{\ell})$, exactly as in Theorem \ref{principalsl2lift}; we then note that since $\tF$ is linearly disjoint from $F(\ad^0 \bar{r}, \zeta_{\ell})$, we can in fact extend $\sigma$ to an element of $\gal{\tF \cdot L}$, and in particular we may regard it as an element of $\gal{\tF, \Sigma}$.
Finally, the group theory establishing the trickiest condition (\ref{avoidout}) was already checked for type $\mr{E}_6$ in the proof of Lemma \ref{magmacalc}: for the desired simple root $\alpha$, we can take any simple root not fixed by the outer automorphism of $\mr{E}_6$.
\endproof
\subsection{Deformations with monodromy group $\mr{E}_6$}
Finally in this section we complete the proof of Theorem \ref{mainintro} by treating the case of $\mr{E}_6$. 
\begin{thm}\label{E6theorem}
There is a density one set of rational primes $\Lambda$ such that for all $\ell \in \Lambda$ there exists a quadratic imaginary field $\tF/\Q$, an almost simple group $G/\Q$ of type $\mr{E}_6$, and an $\ell$-adic representation
\[
\rho_{\ell} \colon \gal{\Q} \to {}^L G(\Qlb)
\]
whose image is Zariski-dense in ${}^L G \cong G^\vee \rtimes \Z/2\Z$. After restriction to $\gal{\tF}$, the image of $\rho_{\ell}|_{\gal{\tF}}$ is Zariski-dense in $G^\vee$.
\end{thm}
\proof
We choose a (non-CM) weight 3 cuspidal eigenform that is new of some level $\Gamma_0(p) \cap \Gamma_1(q)$, exactly as in Theorem \ref{principalsl2lift}, and consider the associated residual representations
\[
\bar{r}_{f, \lambda} \colon \gal{\Q, \Sigma} \to \mr{GL}_2(\mathbb{F}_{\lambda})
\]
where $\Sigma = \{p, q, \ell\}$ ($\lambda \vert \ell$). For a density one set of $\ell$, and a quadratic imaginary field $\tF/\Q$ chosen as in Theorem \ref{principalsl2liftout}, we obtain a homomorphism
\[
\br \colon \gal{\Sigma} \to {}^L G(\mathbb{F}_{\lambda})
\]
satisfying all the hypotheses of Theorem \ref{principalsl2liftout}, where we take (after possibly enlarging $\mathbb{F}_{\lambda}$) Steinberg, minimal, and (sufficiently generic) ordinary deformation conditions at (the specified split place of $\tF$ above) $p$, $q$, and $\ell$, as in Theorem \ref{principalsl2lift}. Let $\rho$ denote the resulting lift to ${}^L G(\Qlb)$. Again by Lemma \ref{reductivemono} and Lemma \ref{avoidsl2}, now applied to $\rho|_{\gal{\tF, \Sigma}}$, we see that the Zariski closure of the image $\rho(\gal{\tF, \Sigma})$ is all of $G^\vee$. 
\endproof
\begin{rmk}
Recall that we had a great deal of flexibility in choosing the field $\tF$, and we acquire more by allowing the modular form $f$ (and, more precisely, its primes $p$ and $q$ of ramification) to vary. Some strengthening of Theorem \ref{E6theorem} is surely possible in which one tries to describe the fields $\tF$ for which the conclusion of the theorem can be shown to hold, but we do not pursue this here.

As with Theorem \ref{notE6theorem}, note that we can check computationally whether a given prime $\ell$ belongs to the density one set admitted in the theorem statement.
\end{rmk}
\bibliographystyle{amsalpha}
\bibliography{biblio.bib}

\providecommand{\bysame}{\leavevmode\hbox to3em{\hrulefill}\thinspace}
\providecommand{\MR}{\relax\ifhmode\unskip\space\fi MR }
\providecommand{\MRhref}[2]{%
  \href{http://www.ams.org/mathscinet-getitem?mr=#1}{#2}
}
\providecommand{\href}[2]{#2}
\begin{thebibliography}{BLGGT14}

\bibitem[And96]{andre:motivated}
Yves Andr{\'e}, \emph{Pour une th\'eorie inconditionnelle des motifs}, Inst.
  Hautes \'Etudes Sci. Publ. Math. (1996), no.~83, 5--49. \MR{1423019
  (98m:14022)}

\bibitem[AP15]{adibhatla-patrikis}
Rajender Adibhatla and Stefan Patrikis, \emph{in preparation}.

\bibitem[BLGGT14]{blggt:potaut}
Thomas Barnet-Lamb, Toby Gee, David Geraghty, and Richard Taylor,
  \emph{Potential automorphy and change of weight}, Ann. of Math. (2)
  \textbf{179} (2014), no.~2, 501--609. \MR{3152941}

\bibitem[Bor79]{borel:L}
A.~Borel, \emph{Automorphic {$L$}-functions}, Automorphic forms,
  representations and {$L$}-functions ({P}roc. {S}ympos. {P}ure {M}ath.,
  {O}regon {S}tate {U}niv., {C}orvallis, {O}re., 1977), {P}art 2, Proc. Sympos.
  Pure Math., XXXIII, Amer. Math. Soc., Providence, R.I., 1979, pp.~27--61.
  \MR{546608 (81m:10056)}

\bibitem[Bou68]{bourbaki:lie456}
N.~Bourbaki, \emph{\'{E}l\'ements de math\'ematique. {F}asc. {XXXIV}. {G}roupes
  et alg\`ebres de {L}ie. {C}hapitre {IV}: {G}roupes de {C}oxeter et syst\`emes
  de {T}its. {C}hapitre {V}: {G}roupes engendr\'es par des r\'eflexions.
  {C}hapitre {VI}: syst\`emes de racines}, Actualit\'es Scientifiques et
  Industrielles, No. 1337, Hermann, Paris, 1968. \MR{0240238 (39 \#1590)}

\bibitem[Car85]{carter:finitelie}
Roger~W. Carter, \emph{Finite groups of {L}ie type}, Pure and Applied
  Mathematics (New York), John Wiley \& Sons, Inc., New York, 1985, Conjugacy
  classes and complex characters, A Wiley-Interscience Publication. \MR{794307
  (87d:20060)}

\bibitem[CHT08]{clozel-harris-taylor}
Laurent Clozel, Michael Harris, and Richard Taylor, \emph{Automorphy for some
  {$l$}-adic lifts of automorphic mod {$l$} {G}alois representations}, Publ.
  Math. Inst. Hautes \'Etudes Sci. (2008), no.~108, 1--181, With Appendix A,
  summarizing unpublished work of Russ Mann, and Appendix B by Marie-France
  Vign{\'e}ras. \MR{2470687 (2010j:11082)}

\bibitem[Con14]{conrad:luminy}
Brian Conrad, \emph{Reductive group schemes}, Autour des sch\'{e}mas en
  groupes, Panoramas et Synth\`eses [Panoramas and Syntheses], vol. 42-43,
  Soci\'et\'e Math\'ematique de France, Paris, 2014, p.~458.

\bibitem[CPS75]{cline-parshall-scott}
Edward Cline, Brian Parshall, and Leonard Scott, \emph{Cohomology of finite
  groups of {L}ie type. {I}}, Inst. Hautes \'Etudes Sci. Publ. Math. (1975),
  no.~45, 169--191. \MR{MR0399283 (53 \#3134)}

\bibitem[DR10]{dettweiler-reiter:rigidG2}
Michael Dettweiler and Stefan Reiter, \emph{Rigid local systems and motives of
  type {$G_2$}}, Compos. Math. \textbf{146} (2010), no.~4, 929--963, With an
  appendix by Michael Dettweiler and Nicholas M. Katz. \MR{2660679
  (2011g:14042)}

\bibitem[Edi97]{edixhoven:serre}
B.~Edixhoven, \emph{Serre's conjecture}, Modular forms and \protect{F}ermat's
  last theorem (Boston, MA, 1995) (New York), Springer, 1997, pp.~209--242.

\bibitem[FG09]{frenkel-gross:rigid}
Edward Frenkel and Benedict Gross, \emph{A rigid irregular connection on the
  projective line}, Ann. of Math. (2) \textbf{170} (2009), no.~3, 1469--1512.
  \MR{2600880 (2012e:14020)}

\bibitem[Gee09]{gee:sato-tatewt3}
Toby Gee, \emph{The {S}ato-{T}ate conjecture for modular forms of weight 3},
  Doc. Math. \textbf{14} (2009), 771--800. \MR{2578803 (2011f:11056)}

\bibitem[Gro97]{gross:principalsl2}
Benedict~H. Gross, \emph{On the motive of {$G$} and the principal homomorphism
  {${\rm SL}_2\to\widehat G$}}, Asian J. Math. \textbf{1} (1997), no.~1,
  208--213. \MR{1480995 (99d:20077)}

\bibitem[HNY13]{heinloth-ngo-yun:kloosterman}
Jochen Heinloth, Bao-Ch{\^a}u Ng{\^o}, and Zhiwei Yun, \emph{Kloosterman
  sheaves for reductive groups}, Ann. of Math. (2) \textbf{177} (2013), no.~1,
  241--310. \MR{2999041}

\bibitem[HSBT10]{hsbt:dwork}
Michael Harris, Nick Shepherd-Barron, and Richard Taylor, \emph{A family of
  {C}alabi-{Y}au varieties and potential automorphy}, Ann. of Math. (2)
  \textbf{171} (2010), no.~2, 779--813. \MR{2630056 (2011g:11106)}

\bibitem[Kat96]{katz:rls}
Nicholas~M. Katz, \emph{Rigid local systems}, Annals of Mathematics Studies,
  vol. 139, Princeton University Press, Princeton, NJ, 1996. \MR{1366651
  (97e:14027)}

\bibitem[Kis07a]{kisin:durham}
Mark Kisin, \emph{Modularity for some geometric {G}alois representations},
  {$L$}-functions and {G}alois representations, London Math. Soc. Lecture Note
  Ser., vol. 320, Cambridge Univ. Press, Cambridge, 2007, With an appendix by
  Ofer Gabber, pp.~438--470. \MR{2392362 (2009j:11086)}

\bibitem[Kis07b]{kisin:cdm}
\bysame, \emph{Modularity of 2-dimensional {G}alois representations}, Current
  developments in mathematics, 2005, Int. Press, Somerville, MA, 2007,
  pp.~191--230. \MR{2459302 (2010a:11098)}

\bibitem[Kle68]{kleiman:algcycles}
S.~L. Kleiman, \emph{Algebraic cycles and the {W}eil conjectures}, Dix
  espos\'es sur la cohomologie des sch\'emas, North-Holland, Amsterdam, 1968,
  pp.~359--386. \MR{0292838 (45 \#1920)}

\bibitem[Kos59]{kostant:principalsl2}
Bertram Kostant, \emph{The principal three-dimensional subgroup and the {B}etti
  numbers of a complex simple {L}ie group}, Amer. J. Math. \textbf{81} (1959),
  973--1032. \MR{0114875 (22 \#5693)}

\bibitem[{LMF}13]{lmfdb}
The {LMFDB Collaboration}, \emph{The l-functions and modular forms database},
  \url{http://www.lmfdb.org}, 2013, [Online; accessed 15 February 2015].

\bibitem[Maz89]{mazur:deforming}
Barry Mazur, \emph{Deforming \protect{G}alois representations}, Galois groups
  over \protect{${\bf Q}$} (Berkeley, CA, 1987), Springer, New York, 1989,
  pp.~385--437.

\bibitem[Pat06]{stp:undergrad}
Stefan Patrikis, \emph{Lifting symplectic galois representations}, Harvard
  undergraduate thesis (2006).

\bibitem[Ram02]{ramakrishna02}
Ravi Ramakrishna, \emph{Deforming {G}alois representations and the conjectures
  of {S}erre and {F}ontaine-{M}azur}, Ann. of Math. (2) \textbf{156} (2002),
  no.~1, 115--154. \MR{MR1935843 (2003k:11092)}

\bibitem[Rib85]{ribet:modlimage2}
Kenneth~A. Ribet, \emph{On {$l$}-adic representations attached to modular
  forms. {II}}, Glasgow Math. J. \textbf{27} (1985), 185--194. \MR{819838
  (88a:11041)}

\bibitem[Sch68]{schlessinger:functors}
Michael Schlessinger, \emph{Functors of {A}rtin rings}, Trans. Amer. Math. Soc.
  \textbf{130} (1968), 208--222. \MR{0217093 (36 \#184)}

\bibitem[Sei91]{seitz:maximalinexceptional}
Gary~M. Seitz, \emph{Maximal subgroups of exceptional algebraic groups}, Mem.
  Amer. Math. Soc. \textbf{90} (1991), no.~441, iv+197. \MR{1048074
  (91g:20038)}

\bibitem[Ser94a]{serre:motivicgalois}
Jean-Pierre Serre, \emph{Propri\'et\'es conjecturales des groupes de {G}alois
  motiviques et des repr\'esentations {$l$}-adiques}, Motives ({S}eattle, {WA},
  1991), Proc. Sympos. Pure Math., vol.~55, Amer. Math. Soc., Providence, RI,
  1994, pp.~377--400. \MR{1265537 (95m:11059)}

\bibitem[Ser94b]{serre:semisimpletensor}
\bysame, \emph{Sur la semi-simplicit\'e des produits tensoriels de
  repr\'esentations de groupes}, Invent. Math. \textbf{116} (1994), no.~1-3,
  513--530. \MR{1253203 (94m:20091)}

\bibitem[Ser96]{serre:principalsl2}
\bysame, \emph{Exemples de plongements des groupes {${\rm PSL}_2({\bf F}_p)$}
  dans des groupes de {L}ie simples}, Invent. Math. \textbf{124} (1996),
  no.~1-3, 525--562. \MR{1369427 (97d:20056)}

\bibitem[SS70]{springer-steinberg:conjugacy}
T.~A. Springer and R.~Steinberg, \emph{Conjugacy classes}, Seminar on
  {A}lgebraic {G}roups and {R}elated {F}inite {G}roups ({T}he {I}nstitute for
  {A}dvanced {S}tudy, {P}rinceton, {N}.{J}., 1968/69), Lecture Notes in
  Mathematics, Vol. 131, Springer, Berlin, 1970, pp.~167--266. \MR{0268192 (42
  \#3091)}

\bibitem[SS97]{saxl-seitz:dynkin}
Jan Saxl and Gary~M. Seitz, \emph{Subgroups of algebraic groups containing
  regular unipotent elements}, J. London Math. Soc. (2) \textbf{55} (1997),
  no.~2, 370--386. \MR{1438641 (98m:20057)}

\bibitem[Ste68]{steinberg:chevalley}
Robert Steinberg, \emph{Lectures on {C}hevalley groups}, Yale University, New
  Haven, Conn., 1968, Notes prepared by John Faulkner and Robert Wilson.
  \MR{0466335 (57 \#6215)}

\bibitem[SW01]{skinner-wiles:basechange}
C.~M. Skinner and A.~J. Wiles, \emph{Base change and a problem of {S}erre},
  Duke Math. J. \textbf{107} (2001), no.~1, 15--25. \MR{1815248 (2002c:11058)}

\bibitem[Til96]{tilouine:defs}
Jacques Tilouine, \emph{Deformations of {G}alois representations and {H}ecke
  algebras}, Published for The Mehta Research Institute of Mathematics and
  Mathematical Physics, Allahabad; by Narosa Publishing House, New Delhi, 1996.
  \MR{1643682 (99i:11038)}

\bibitem[Wes04]{weston:unobstructed}
Tom Weston, \emph{Unobstructed modular deformation problems}, Amer. J. Math.
  \textbf{126} (2004), no.~6, 1237--1252. \MR{2102394 (2006c:11061)}

\bibitem[Wil88]{wiles:ordinary}
A.~Wiles, \emph{On ordinary {$\lambda$}-adic representations associated to
  modular forms}, Invent. Math. \textbf{94} (1988), no.~3, 529--573.
  \MR{MR969243 (89j:11051)}

\bibitem[Yun14]{yun:exceptional}
Zhiwei Yun, \emph{Motives with exceptional {G}alois groups and the inverse
  {G}alois problem}, Invent. Math. \textbf{196} (2014), no.~2, 267--337.

\end{thebibliography}
\end{document}